\magnification=1100
\baselineskip=14truept
\voffset=.75in
\hoffset=1truein
\hsize=4.5truein
\vsize=7.75truein
\parindent=.166666in
\pretolerance=500 \tolerance=1000 \brokenpenalty=5000
\def\anote#1#2#3{\smash{\kern#1in{\raise#2in\hbox{#3}}}%
  \nointerlineskip}     
\def\note#1{%
  \hfuzz=50pt%
  \vadjust{%
    \setbox1=\vtop{%
      \hsize 3cm\parindent=0pt\eightrm\baselineskip=9pt%
      \rightskip=4mm plus 4mm\raggedright#1%
      }%
    \hbox{\kern-4cm\smash{\box1}\hfil\par}%
    }%
  \hfuzz=0pt
  }
\def\note#1{\relax}

\newcount\equanumber
\equanumber=0
\newcount\sectionnumber
\sectionnumber=0
\newcount\subsectionnumber
\subsectionnumber=0
\newcount\snumber  
\snumber=0

\def\section#1{%
  \subsectionnumber=0%
  \snumber=0%
  \equanumber=0%
  \advance\sectionnumber by 1%
  \noindent{\bf \the\sectionnumber .~#1~}%
}%
\def\subsection#1{%
  \advance\subsectionnumber by 1%
  \snumber=0%
  \equanumber=0%
  \noindent{\bf \the\sectionnumber .\the\subsectionnumber .~#1~}%
}%
\def\prevs{\the\sectionnumber .\the\subsectionnumber .\the\snumber }
\long\def\Definition#1{\noindent{\bf Definition.}{\it#1}}
\long\def\Claim#1{\noindent{\bf Claim.}{\it#1}}
\long\def\Corollary#1{%
  \global\advance\snumber by 1%
  \bigskip
  \noindent{\bf Corollary~\prevs .}%
  \quad{\it#1}%
}%
\long\def\Lemma#1{%
  \global\advance\snumber by 1%
  \bigskip
  \noindent{\bf Lemma~\prevs .}%
  \quad{\it#1}%
}%
\def\Proof{\noindent{\bf Proof.~}}
\long\def\Proposition#1{%
  \advance\snumber by 1%
  \bigskip
  \noindent{\bf Proposition~\prevs .}%
  \quad{\it#1}%
}%
\long\def\Remark#1{%
  \bigskip
  \noindent{\bf Remark.~}#1%
}%
\long\def\Theorem#1{%
  \advance\snumber by 1%
  \bigskip
  \noindent{\bf Theorem~\prevs .}%
  \quad{\it#1}%
}%
\def\ifundefined#1{\expandafter\ifx\csname#1\endcsname\relax}
\def\labeldef#1{\global\expandafter\edef\csname#1\endcsname{\prevs}}
\def\labelref#1{\expandafter\csname#1\endcsname}
\def\label#1{\ifundefined{#1}\labeldef{#1}\note{$<$#1$>$}\else\labelref{#1}\fi}

\def\preveq{(\the\sectionnumber .\the\subsectionnumber .\the\equanumber)}
\def\neq{\global\advance\equanumber by 1\eqno{\preveq}}

\def\ifundefined#1{\expandafter\ifx\csname#1\endcsname\relax}

\def\equadef#1{\global\advance\equanumber by 1%
  \global\expandafter\edef\csname#1\endcsname{\preveq}%
  \setbox1=\hbox{\hskip .1in[#1]}\dp1=0pt\ht1=0pt\wd1=0pt%
  \preveq\box1}
\def\equadef#1{\global\advance\equanumber by 1%
  \global\expandafter\edef\csname#1\endcsname{\preveq}%
  \preveq}

\def\equaref#1{\expandafter\csname#1\endcsname}
 
\def\equa#1{%
  \ifundefined{#1}%
    \equadef{#1}%
  \else\equaref{#1}\fi}

\def\poorBold#1{\setbox1=\hbox{#1}\wd1=0pt\copy1\hskip.25pt\box1\hskip .25pt#1}

\font\funfont=cmb10%
\font\largecmr=cmr12 at 24pt%
\font\eightrm=cmr8%
\font\sixrm=cmr6%
\font\fourrm=cmr4%
\font\eighttt=cmtt8%
\font\eightbf=cmb8%
\font\elevenbb=msbm10%
\font\eigthbb=msbm8%
\font\sixbb=msbm6%
\font\bftt=cmbtt10%
\newfam\bbfam%
\textfont\bbfam=\elevenbb%
\scriptfont\bbfam=\eigthbb%
\scriptscriptfont\bbfam=\sixbb%
\font\tencmssi=cmssi10%
\font\eightcmssi=cmssi7%
\font\sixcmssi=cmssi6%
\newfam\ssfam%
\textfont\ssfam=\tencmssi%
\scriptfont\ssfam=\eightcmssi%
\scriptscriptfont\ssfam=\sixcmssi%
\def\ssi{\fam\ssfam\tencmssi}%

\font\tenmsam=msam10%
\font\eigthmsam=msam7%
\font\sixmsam=msam6%

\def\bb{\fam\bbfam\elevenbb}%

\def\hexdigit#1{\ifnum#1<10 \number#1\else%
  \ifnum#1=10 A\else\ifnum#1=11 B\else\ifnum#1=12 C\else%
  \ifnum#1=13 D\else\ifnum#1=14 E\else\ifnum#1=15 F\fi%
  \fi\fi\fi\fi\fi\fi}
\newfam\msamfam%
\textfont\msamfam=\tenmsam%
\scriptfont\msamfam=\eigthmsam%
\scriptscriptfont\msamfam=\sixmsam%
\mathchardef\leq"3\hexdigit\msamfam 36%
\mathchardef\geq"3\hexdigit\msamfam 3E%

\def\convM{\buildrel{\raise-1pt\hbox{\sixrm M}}\over *}
\def\convMss{\setbox1=\hbox{\fourrm M}\ht1=0pt\dp1=0pt%
  \buildrel{\raise -2pt\box1}\over *}
\def\<{\langle}
\def\>{\rangle}
\def\AmoGnk{{\cal A}_m\oGnk}
\def\cX{\<c,X\>}
\def\d{\,{\rm d}}
\def\diag{{\rm diag}}
\def\D{{\rm D}}
\long\def\DoNotPrint#1{\relax}
\def\ds{\displaystyle}
\def\equald{\buildrel{\rm d}\over=}
\def\finetune#1{#1}
\def\fixedref#1{#1}

\def\Id{{\rm Id}}
\def\limn{\lim_{n\to\infty}}
\def\limt{\lim_{t\to\infty}}
\def\limsupt{\limsup_{t\to\infty}}
\def\L{{\ssi L}}
\def\M{{\ssi M}}
\def\P#1{P\{\,#1\,\}}
\def\oB{\overline B}
\def\oDelta{\overline\Delta}
\def\oF{\overline F{}}

\def\oFik{\overline F_i^{(k)}}
\def\oFnk{\overline F_n^{(k)}}
\def\oG{\overline G{}}
\def\oGn{\oG_n}
\def\oGnk{\oGn^{(k)}}
\def\oH{\overline H{}}
\def\oK{\overline K{}}
\def\onu{\overline \nu}
\def\oW{\overline W{}}
\def\qed{~~~{\vrule height .9ex width .8ex depth -.1ex}}
\def\sign{\hbox{\rm sign}}
\def\sumn{\sum_{n\geq 1}}
\def\ss{\scriptstyle}
\def\T{\hbox{\tencmssi T}}
\def\U{\hbox{\tencmssi U}}

\def\HH{{\bb H\kern .5pt}}
\def\II{{\bb I}}

\def\NN{{\bb N}\kern .5pt}

\def\RR{{\bb R}}
\def\ZZ{{\bb Z}}

\def\calA{{\cal A}}
\def\calC{{\cal C}}
\def\calD{{\cal D}}

\def\calL{{\cal L}}
\def\calM{{\cal M}}
\def\calN{{\cal N}}
\def\calP{{\cal P}}

\long\def\donotprint#1{\relax}

%
%
%
%
\def\uncatcodespecials 
    {\def\do##1{\catcode`##1=12}\dospecials}%
{\catcode`\^^I=\active \gdef^^I{\ \ \ \ }
 \catcode`\`=\active\gdef`{\relax\lq}}
\def\setupverbatim 
    {\parindent=0pt\tt %
     \spaceskip=0pt \xspaceskip=0pt 
     \catcode`\^^I=\active %
     \catcode`\`=\active %
     \def\par{\leavevmode\endgraf}
     \obeylines \uncatcodespecials \obeyspaces %
     }%
{\obeyspaces \global\let =\ }
%
%
%
%
%
%

\pageno=1

\centerline{\bf ASYMPTOTIC EXPANSIONS}
\centerline{\bf FOR INFINITE WEIGHTED CONVOLUTIONS}
\centerline{\bf OF HEAVY TAIL DISTRIBUTIONS}
\centerline{\bf AND APPLICATIONS} 

\medskip
 
\centerline{Ph.\ Barbe and W.P.\ McCormick}
\centerline{CNRS, France, and University of Georgia}
 
{\narrower
\baselineskip=9pt\parindent=0pt\eightrm

\bigskip

{\eightbf Abstract.} We establish some asymptotic
expansions for infinite weighted convolution of distributions having regular
varying tails. Various applications to statistics and probability
are developed.

\bigskip

\noindent{\eightbf AMS 2000 Subject Classifications:}
Primary: 41A60, 60F99.
Secondary: 41A80, 44A35, 60E07, 60G50, 60K05, 60K25, 62E17, 62G32.

\bigskip
 
\noindent{\eightbf Keywords:} asymptotic expansion, regular variation,
convolution, tail area approximation, ARMA models, tail estimation, randomly 
stopped sums, infinitely divisible distributions, renewal theory.

}
\vfill\eject

{\eightrm
\obeylines\baselineskip=9pt\setbox1=\hbox{1. }
\newdimen\dd
\dd=\wd1
\def\quad{\hskip\dd}
{\eightbf Contents}
1. Introduction
\quad 1.1. Prolegomenom
\quad 1.2. Mathematical overview and heuristics
\quad 1.3. Asymptotic scales
2. Main result
\quad 2.1. The Laplace characters
\quad 2.2. Smoothly varying functions of finite order
\quad 2.3. Asymptotic expansion for infinite weighted convolution
3. Implementing the expansion
\quad 3.1. How many terms are in the expansion?
\quad 3.2. Practical implementation: from Laplace's characters to linear algebra
\quad 3.3. Two terms expansion and second order regular variation
\quad 3.4. Some open questions
4. Applications
\quad 4.1. ARMA models
\quad 4.2. Tail index estimation
\quad 4.3. Randomly weighted sums
\quad 4.4. Randomly stopped sums
\quad 4.5. Queueing theory
\quad 4.6. Branching process
\quad 4.7. Infinitely divisible distributions
\quad 4.8. Implicit transient renewal equation and iterative systems
5. Proof in the positive case
\quad 5.1. Decomposition of the convolution into integral and multiplication 
\quad\phantom{5.1. }operators
\quad 5.2. Organizing the proof
\quad 5.3. Regular variation and basic tail estimates
\quad 5.4. The fundamental estimate
\quad 5.5. Basic lemmas
\quad 5.6. Inductions
\quad 5.7. Conclusion
6. Removing the sign restriction on the random variables
\quad 6.1. Elementary properties of $\ss U_H$
\quad 6.2. Basic expansion of $\ss U_H$
\quad 6.3. A technical lemma
\quad 6.4. Conditional expansion and removing conditioning
7. Removing the sign restriction on the constants
\quad 7.1. Neglecting terms involving the multiplication operators
\quad 7.2. Substituting $\ss \oH^{(k)}$ and $\ss \oG^{(k)}$ by their expansions
8. Removing the smoothness restriction
Appendix. {\eighttt Maple} code
References
}

\vfill\eject
 
\section{Introduction.}\quad
The primary focus of this paper is to obtain precise understanding of the
distribution tail of linear and related stochastic processes based on
heavy tail innovations. In doing so, we will develop some new mathematical
objects which are tailored to efficiently write and compute asymptotic
expansions of these tails. Also, we will derive simple bounds of theoretical
importance for the error between the tail and its asymptotic expansion.

These tails and their expansions are of interest in a variety of contexts.
In the following subsection, we provide some typical examples to illustrate
their use.
The second subsection of this introduction overviews the new perspective
and techniques developed in this paper; this will be done at a heuristic
level, explaining the intuition and sketching the broad expanse wherein our
methods lie. The last subsection contains basic facts on asymptotic scales.

\bigskip
\subsection{Prolegomenom.}%
The first basic problem we will deal with is related to the tail
behavior of the marginal distribution of linear processes.  To be
specific, let $c=(c_i)_{i\in\ZZ}$ be a sequence of real constants, and
let $X=(X_i)_{i\in\ZZ}$ be a sequence of independent and identically
distributed random variables. Let $F$ be the common distribution
function of these $X_i$'s and write $\oF=1-F$ the tail
function. Assume that $F$ has a heavy tail, that is, for any positive
$\lambda$,
$$
  \lim_{t\to\infty} \oF(t\lambda)/\oF(t) = \lambda^{-\alpha}
  \eqno{\equa{IntroA}}
$$
for some positive and finite $\alpha$. Let $G_c$ be the distribution function
of the series $\sum_{i\in\ZZ} c_iX_i$. Set $\oG_c=1-G_c$.
For simplicity, assume in this introduction that the sequence $c$ 
is nonnegative. Define $C_\alpha=\sum_{i\in\ZZ} c_i^\alpha$. It is 
well known that under a mild additional condition,
$$
  \oG_c\sim C_\alpha \oF
$$
at infinity, that is $\lim_{t\to\infty} \oG_c(t)/\oF(t) = C_\alpha$.

It is then a natural question to investigate higher order expansions.
Under suitable conditions, we will obtain higher order asymptotic expansions
both for $G_c$ and its derivatives. In particular all ARMA processes fall
within the scope of our result, provided the innovation distribution satisfies
certain mild conditions beyond that of heavy tail.

This question of obtaining higher order terms has several
motivations. One is that very little is known on the marginal
distributions of linear processes, and since those are ubiquitous in
statistics, better knowledge and understanding of their
properties is desirable. There are also specific instances where
higher order asymptotics are needed, such as in the tail index
estimation problem which we will study in section 4.2. To achieve such
refined distributional results in this setting, we must first build
the mathematical language and theory needed for higher order
expansions.

A third motivation comes from a large number of applications to
related processes obtained by allowing the weights to be random. Whereas
many first order results are known, few second order results exist,
and no higher order results are available that we are aware of.  We
will obtain asymptotic expansions when the weights are random, not
necessarily independent or identically distributed, but still assuming
that the sequence $c$, now random, is independent of $X$.

This will allow us to derive tail expansions for compound sums, that is a
sum with a random number of summands. As a consequence, we will derive a tail
expansion for some infinitely divisible distributions. In the same vein,
having a result with random weights yields expansions for tails arising in 
transient renewal theory, and even in implicit transient renewal theory.
By the latter, we mean, for example, distributions defined implicitly as 
follows. Let $M$, $Q$ and $R$ be independent nonnegative random variables
and suppose that $R$ has the same distribution as $Q+MR$. Given the 
distributions of $M$ and $Q$, this defines that of $R$ implicitly. We will
obtain an asymptotic expansion for the distribution of $R$, when $Q$ has
heavy tail and $M$ has enough moments.

Furthermore, it will be clear that by the algebraic paradigm of this paper, 
it is quite a simple matter to write out what the expansions should be, and 
often not so difficult to prove that those expansions are valid under some
reasonable conditions.

In all the contexts mentioned, first order results exist under heavy
tail assumptions, or more generally in the subexponential framework.
By now, these results are well understood. We refer to the book by
Bingham, Goldie and Teugels (1989) for a comprehensive study of these
topics and Broniatowski and Fuchs (1995) for a different perspective
on first order results for sums of independent and identically
distributed random variables. Resnick (1986) is a good reference for
many first order results obtained through a point process argument.
In comparison, second order results are few. In the heavy-tail case
papers which provide second order results for the sum of two random
variables include Omey (1988), Geluk (1994), Geluk et al.~(1997) and
Geluk et al.~(2000). With regard to subexponential distributions,
second order work for a convolution of two such distributions may be
found in Cline (1986), Cline (1987) and Geluk and Pakes (1991). We
also mention that Geluk et al.~(2000) shows a second order result for
a sum of a finite number of independent and identically distributed
random variables in the special case of the underlying distribution
being a member of the Hall-Weissman class.  Barbe and McCormick (2005)
obtained second order result for the sum of a finite number of heavy
tailed random variables when the distributions also possess a mild
smoothness property. Using a different set of assumptions and for
a different purpose, Borovkov and Borovkov (2003) obtained higher
order expansions also for a finite number of summands.

In renewal theory, a little more is known. For instance, in recurrent
renewal theory, higher order results are already in Feller (1971). But
our technique is useful only in the transient case, for which far less
is known.  Chover et al.~(1973) provides a first order
result while Gr\"ubel (1987) gives a second order formula. In the
implicit case, Goldie (1991) obtained second order formulas. However,
our work complements nicely with Goldie's, in the sense that our results
typically apply when his fail. The techniques used by Goldie (1991)
build on those devised by Kesten (1973) and Grincevi\v cius
(1975), and make use of delicate Tauberian theory, while that of
Gr\"ubel (1987) is the so-called Banach algebra technique.  Our
technique is radically different from these.

Gr\"ubel's (1987) result on renewal measures allows him to obtain a second
order formula for the tail of some infinitely divisible distributions, using a
decomposition which we could trace back at least to Feller (1971). We will
use the same decomposition to obtain higher order results.

Since our approach to deriving these asymptotic expansions may appear
very algebraic to analytically oriented minds, we feel compelled at this point
to quote a result which can be proved in less than 5 minutes after 
understanding the statement of our main result (Theorem \fixedref{2.3.1})
and being familiar with our algebraic formalism. So assume that $F$ is the 
distribution function defined by
$\oF(t)=(a+b)^{-1} (at^{-\alpha} + bt^{-\beta})$ for $t$ more than $1$,
and with $1<\alpha<\beta$. Assume further that the $c_i$'s are nonnegative.
Write $C_r=\sum_{i\in\ZZ} c_i^r$, and 
$$
  \mu_{F,1}={1\over a+b} \Bigl({a\over \alpha-1} + {b \over \beta-1}\Bigr)
$$
for the first moment of $F$. Then we have the two terms expansion
$$\displaylines{\quad
  \oG_c(t)={a\over a+b} {C_\alpha\over t^\alpha}
  \hfill\cr\noalign{\vskip 3pt}\hfill
  {}+\left\{ 
   \matrix{\ds {bC_\beta\over a+b} {1\over t^\beta} \hfill
            &\ds \hbox{ if $\beta<\alpha+1$,}\cr
            \noalign{\vskip 3pt}
            {\ds 1\over\ds a+b} 
            \bigl( bC_\beta +a\mu_{F,1}\alpha 
                   (C_1C_\alpha-C_{\alpha+1}) \bigr) 
                   {\ds 1\over\ds t^{\alpha+1}} \hfill
            &\ds \hbox{ if $\beta=\alpha+1$,}\cr
            \ds{a\alpha\over a+b}\mu_{F,1}(C_1C_\alpha-C_{\alpha+1}) 
            {1\over t^{\alpha+1}} \hfill
            &\ds \hbox{ if $\alpha+1<\beta$.}\cr
           }
   \right\}
  \hfill\cr\quad
   {}+ o\Bigl({1\over t^{\alpha+1}}\Bigr)
  \hfill\cr}
$$
With our formalism and notations to be explained later, this complicated 
formula has the far more manageable form, where no cases need to be 
distinguished,
$$
  \oG_c\sim \sum_{i\in\ZZ}\L_{G_c\sharp\M_{c_i}F,1}\overline{\M_{c_i}F} \, . 
$$
Note that if $\beta$ is $\alpha+1$, the $c_i$'s and $b$ can be
chosen so that the above second order expansion yields only
$$
  \oG_c={a\over a+b} C_\alpha t^{-\alpha} + o(t^{-\alpha-1}) \, .
$$
It is then desirable to get an extra term. In section \fixedref{3.1}, we 
will obtain this term, generically, when $\alpha$ is more than $2$. 
We will also explain
why in some very exceptional cases, obtaining a second term is rather hopeless
(last example of section \fixedref{3.1}), while in other cases, obtaining
as many terms as one desires can be trivial (penultimate example in section
\fixedref{3.2}).

To conclude this preliminary discussion, and perhaps to motivate further
our investigation, we mention that there are many unsolved very basic
problems related to linear processes. For instance, for discrete
innovation, Davis and Rosenblatt (1991) show that all but trivial infinite
order linear processes have a continuous marginal distribution; 
but there is still no good criterion to know if the marginal
distribution of the process is absolutely continuous. Even in the most
basic cases, the behavior of the marginal distribution is amazingly
difficult to analyze.  Solomyak's (1995) breakthrough --- a generic
result on absolute continuity of the marginal distribution of first
order autoregressive models with Bernoulli innovations --- gives a
chilling reminder on how little is known in general.  So this paper
can be taken as a contribution to the understanding of those
distributions, in the continuous and heavy tail situation.

\bigskip

\subsection{Mathematical overview and heuristics.}%
We now outline the mathematical content of the key parts of this paper
and provide some intuition behind the main result.

In our opinion, the main reason why so few higher order results are
available in the problems we are interested in, is that too much emphasis
has been given to an analytical
perspective. Our first basic remark is that the convolution operation
$(F,G)\mapsto F\star G$ is bilinear and defines a semi-group. Bilinearity
is a notion in linear algebra; while semi-group is related to the group 
structure, which in turn 
suggests representation theory. So our view is that the asymptotic behavior
of the convolution semi-group should be analyzed by a linear representation
which captures only the tail behavior of the semi-group. In fact we will
obtain two linear representations, one which is suitable for writing the
expansions, the other one, derived from the first one, and requiring
more assumptions, suitable for practical computations.

To explain further our view on the subject, let us give a heuristic argument
on how to derive higher order expansions and build the mathematical language
to handle these expansions.

\medskip

\noindent{\funfont Heuristic.} 
Let $F$ and $G$ be two distribution functions. Their convolution is
$$
  \overline{F\star G}(t) = \int \oF(x-y)\d G(y) \, .
$$
The left hand side in this formula is symmetric, while the right hand side, as
far as the notation goes, is not. To symmetrize it, we split the integral
and integrate by parts to obtain
\hfuzz=2pt
$$
  \overline{F\star G}(t) = \int_{-\infty}^{t/2} \oF (t-x)\d G(x)
  + \int_{-\infty}^{t/2} \oG(t-x) \d F(x) + \oF(t/2)\oG(t/2) \, .
$$
\hfuzz=0pt

This formula involves the translation $\tau_x$ acting on 
functions and defined by $\tau_x h(t) =h(t+x)$. The translations 
form a semigroup since $\tau_x\tau_y=\tau_{x+y}$. The corresponding 
infinitesimal generator is the derivative $\D$, because for smooth 
functions
$$
  \lim_{\epsilon\to 0} {\tau_\epsilon-\tau_0\over\epsilon} h(t)
  = \lim_{\epsilon\to 0} {h(t+\epsilon)-h(t)\over\epsilon}
  = \D h(t) \, .
$$
So, one expects to be able to write $\tau_x=e^{x\D}$ --- this is nothing
but a neat way to write a Taylor expansion, which we could trace back to 
Delsarte (1938). Indeed, applying both sides to a smooth function $h$ yields
$$
  h(t+x)=\tau_x h(t) = e^{x\D} h(t) =\sum_{i\geq 0}{x^i\D^i \over i!} h(t)
  = \sum_{i\geq 0} {x^i\over i!} h^{(i)}(t) \, .
$$

So, formally,
$$
  \int_{-\infty}^{t/2} \oF(t-x) \d G(x)
  = \int_{-\infty}^{t/2} \tau_{-x} \d G(x) \oF (t)
  = \int_{-\infty}^{t/2} e^{-x\D} \d G(x) \oF (t) \, .
$$
As $t$ tends to infinity, the linear operator 
$\int_{-\infty}^{t/2} e^{-x\D}\d G(x)$ should
tend to the Laplace transform of $G$ at $\D$, say $L_G(\D)$. 
So, ideally, we have
$$
  \overline{F\star G} = L_G(\D)\oF + L_F(\D)\oG 
  + O(\oF\,\oG) \, .
  \eqno{\equa{HeuristicA}}
$$
Of course, for heavy tail distributions, many things are not well defined.
Certainly, the Laplace transform does not exist, the dis\-tri\-bu\-tion 
lacking sufficient integrability.

So, consider the Taylor polynomial of the Laplace transform. Write
$\mu_{G,k}$ for the $k$-th moment of $G$. Assume that these moments
are finite for $k$ at most equal to some $m$. The Taylor polynomial
of the Laplace transform gives rise to the differential operator
$$
  \L_{G,m}=\sum_{0\leq j\leq m} {(-1)^j\over j!} \mu_{G,j} \D^j \, .
$$

Assume that $\oF$ and $\oG$ are regularly varying of index $-\alpha$,
that is obey \IntroA. Take $m$ be an integer less than $\alpha$.
Throughout this paper, $\Id$ denotes the identity function of whatever
space is under consideration. When using the
identity on the real line, we write $\Id^k$ for the function $t\mapsto t^k$. 
So, $\Id^{-m}\oF$ is the function whose value
at $t$ is $t^{-m}\oF(t)$. We certainly would
like to have the $m$-th derivative $\oF^{(m)}$ regularly varying of index 
$-\alpha-m$. In that case $(L_G-\L_{G,m})\oF$ could be of order smaller
than the last term, that is
$\oF^{(m)}$, or equivalently, $\Id^{-m}\oF$, and therefore dominates
$\oF\,\oG$. So replacing $L_G(\D)$ and 
$L_F(\D)$ in {\HeuristicA} by $\L_{G,m}$ and $\L_{F,m}$
and taking $m$ less than $\alpha$ suggests
$$
  \overline{F\star G} = \L_{G,m}\oF + \L_{F,m}\oG 
  + o(\Id^{-m}\oF) \, .
  \eqno{\equa{HeuristicB}}
$$
The advantage of this substitution is that now everything is well defined,
in particular it does not assume that $F$ and $G$ have moments of arbitrary
large order.

By induction, if we have $n$ distribution functions $F_i$, $1\leq i\leq n$, 
we obtain
$$
  \overline{\star_{1\leq i\leq n} F_i} 
  = \sum_{1\leq i\leq n} 
  \L_{\ds\star_{\matrix{\ss 1\leq j\leq n\cr
                        \noalign{\vskip -3pt}
                        \ss j\not= i\hfill\cr}}
      \ss F_j,m}
  \oF_i + \hbox{ remainder.}
$$
Taking limit at infinity yields a formula, which, restricted to the context 
of linear processes, is in fact the main result of this paper.

The substantial task of transforming this heuristic argument into a
rigorous proof will be carried out in sections 5, 6, 7 and 8.

Going back to our heuristic, where is the representation which 
we announced? We will show that
in the proper ring, that of polynomials in $\D$ modulo the ideal generated
by $\D^{m+1}$, the composition of $\L_{F,m}$ with $\L_{G,m}$ gives
$\L_{F\star G,m}$. So the map $F\mapsto \L_{F,m}$ turns out to be a linear
representation of the convolution algebra in a ring of differential
operators. In this ring, the operators $\L_{F,m}$ will be invertible.

This gives us a very powerful formalism to state expansions.
But, arguably, it is not quite explicit. To obtain a powerful computational
machinery, we need to go from a ring of differential operators to a
vector space where we can manipulate finite dimensional
matrices. This is not always possible though, and it is related to a rather
intricate matter concerning asymptotic scales. Those will be introduced in the
next subsection and further examined in our context in subsections 3.1 and
3.2. 

Nevertheless, the basic idea is quite intuitive. For any positive real number
$c$, write $\M_cF$ to denote the distribution function of $cX_i$. More 
generally, when $c$ is positive and $h$ is a function, we write $\M_ch$ for 
the function whose value at $t$ is $h(t/c)$. Since 
$\D\overline{\M_cF}=c^{-1}\M_c\D\oF$, the expansion {\HeuristicB} with $m=1$ 
yields
$$\displaylines{\quad
  \overline{\M_{c_1}F\star \M_{c_2}F}
  =\M_{c_1}\oF +\M_{c_2}\oF 
  - c_1\mu_{F,1}c_2^{-1}\M_{c_2}\D\oF 
  \hfill\cr\hfill
  {}- c_2\mu_{F,1}c_1^{-1}\M_{c_1}\D\oF
  +\hbox{remainder.}\quad\equa{HeuristicC}\cr
  }
$$
This suggests that we could think of such expansion as a decomposition in
a basis comprised of $\oF$ and $\D\oF$, that is as a sort of projection onto 
the 2-dimensional vector space spanned by these functions.
But there is a problem to overcome, namely that the sum 
$\M_{c_1}\oF +\M_{c_2}\oF$
is not quite the usual first order result $(c_1^\alpha+c_2^\alpha)\oF$.
In doing such a replacement, we make an error, which may be larger than
the second order term. So, even in an asymptotic sense, $\M_c\oF$
may not be in the space spanned by $\oF$ and $\D\oF$. We will give examples
in sections 3.1 and 3.2 showing incidentally that this problem is not 
a failure of the expansion, but that of
the basis $\oF$, $\D\oF$. This suggests that some bases are better than
others, and of course raises the question of what is a good basis for our
purpose.

Since we will now be dealing with issues related to asymptotic analysis,
we will switch to a terminology used in that field. Specifically, we will refer
to an asymptotic scale for that which we called a basis in the previous 
paragraph. Asymptotic 
scales will be introduced formally in the next subsection.

In studying linear processes, two operations are involved. One is the
convolution obtained from summing independent random variables, the other
is the scaling $\M_c$ obtained by rescaling the variables. We have
seen that as far as the convolution goes, we have a representation using
the differential operators $\L_{F,m}$, which themselves involve $\D$. This
suggests that the good asymptotic scales are in some sense stable by
differentiation and rescaling. In which sense? Since we are interested
in asymptotic expansions, and we want to preserve the linear aspect of 
the convolution, they should be stable in the sense that
derivative and scaling of functions in this scale should have an asymptotic
expansion in the same scale. The simplest example of such functions are the
powers $t^{-\alpha-n}$, $n\in\NN$. We will see many others.  We call such
stable scale a $\star$-asymptotic scale, and that will be defined
formally in subsection 3.2.

So, assume that $e=(e_i)_{i\in I}$ is a set 
of functions defining our $\star$-asymptotic scale. Then there exists a
matrix $\calD=(\calD_{i,j})_{i,j\in I}$, such that for any $i$ in $I$, the
derivative
$\D e_i$ has asymptotic expansion $\sum_{j\in I }\calD_{i,j}e_j$; and
for any positive $c$, there exists a matrix $\calM_c$, such that $\M_ce_i$ has
asymptotic expansion $\sum_{j\in I}{(\calM_c)}_{i,j} e_j$.
In particular, defining the matrix 
$$
  \calL_G=\sum_{0\leq j\leq m} {(-1)^j\over j!} \mu_{G,j} \calD^j \, ,
$$
we see that $\L_{G,m} e_i$ has asymptotic expansion 
$\sum_{j\in I}({\calL_G})_{j,i}e_j$. Since the differential operators 
$\L_{G,m}$ are representations of the convolution semi-group, so are the 
matrices $\calL_{G}$. This is our second representation of the convolution 
semigroup, which is now made up of finite dimensional matrices. Note that
compared to the first representation using differential operators, we assume
the existence of a $\star$-asymptotic scale in which functions can be expanded.
We will see that for usual distributions this is not a restriction at all,
and it is even very natural and desirable. Also, the representation $\L_{F,m}$
has dimension $m+1$, while that given by $\calL_F$ will often be of higher
dimension.

Now, if $\oF$ has an asymptotic expansion $\sum_{i\in I}p_{\oF,i}e_i$, this
expansion is encoded in the vector $p_{\oF}=(p_{\oF,i})_{i\in I}$. Then, the
asymptotic expansion of $\overline{\M_cF}$ is encoded in the vector 
$\calM_c p_{\oF}$, which is simply a product of a matrix by a vector.
Thus one can hope that because $\L_{G,m}$ is a linear operator, 
$\L_{G,m}\oF$ has an asymptotic expansion encoded by the vector
is $\calL_Gp_{\oF}$. In particular, formula {\HeuristicC} can be written as
$$
  p_{\overline{\M_{c_1}F\star\M_{c_2}F}}
  = (\calL_{\M_{c_1}F}\calM_{c_2}+\calL_{\M_{c_2}F}\calM_{c_1})p_{\oF} \, .
$$
This means that we can calculate the asymptotic expansion of $\M_{c_1}F\star
\M_{c_2}F$ in a $\star$-asymptotic scale using multiplication and addition
of matrices and vectors.
The advantage of using this matrix representation is that we are now in the
realm of linear algebra, where computers can be used both to do numerical
and formal work. This is the key to efficient derivations of expansions
in practical cases. One last point worth mentioning with regard to this
algebraic approach is that the complexity of argument to obtain a higher
order expansion is essentially the same as that to obtain a second order
expansion. The effect of higher order is to raise the dimension of the space
one works with; but, for second order and beyond, the spaces used will all have
dimension at least 2. Only a first order result for which a 1-dimensional
space suffices can attain a significant reduction in effort. In the same vein, 
we note that expansions on densities and their derivatives are obtained with 
little extra effort. Furthermore, these results are useful, for example, in 
discussing the Von Mises condition for infinite order weighted averages
--- see subsection \fixedref{4.2}.

\medskip

Of course, all this heuristic elaboration is rather formal at this point, and
clearly the realization of this program requires us to define the proper
analytical setup --- this will be done in subsection 2.1 --- as well as
the algebraic machinery --- to be done in part in subsection 2.2 --- 
simply to be able to state the proper theorem in subsection 2.3. 

\medskip

Up to now we have mainly presented the algebraic formalism. It is clear that 
in order to justify rigorously these heuristic arguments, one has to connect 
this formalism to analysis. As mentioned earlier, this is mostly done 
in sections 5, 6 and 7. At this point, it is not useful to outline the proof;
such outline is given in section 5.2. But we mention that the proof is
by induction on the number of summands, and then a limiting argument
as the number of summands goes to infinity. Of course, any limiting argument
in this type of asymptotics requires a good error bound between the original
functions and their approximations, while induction somewhat requires
tractable bounds. Those two requirements, sharpness and simplicity, tend
to be in opposition.
The key to obtain both features will be a functional analytic approach.
We will write the convolution in terms of some operators, and the study
of those operators will give us good bounds.

The many steps alluded to in our heuristic arguments may suggest that 
particularly severe restrictions will be needed on the distribution function 
of the $X_i$'s and the constants  $c$. However, this is not the case and
the main result is rather sharp. Before going any further,
we mention that our asymptotic expansions involve derivatives and
moments. It will be clear from the proof that some regularity of
the underlying distribution must be assumed in order to obtain
even a two terms expansion. The regularity assumption has a great
impact on the form of the expansion. The assumption chosen hereafter,
namely differentiability up to a certain order, seems the most natural
one for applications to concrete distributions. In the same vein,
Barbe and McCormick (2005) show that when $\alpha$ is $1$ and the 
first moment is infinite, 
the truncated first moment is involved in a two terms expansion for 
finite convolution; moreover, the expansion takes a very different
form if $\alpha$ is less than $1$. It is quite clear from 
Barbe and McCormick (2005) how some results of the present paper can
be formally modified when one does not have enough moments. The overall
philosophy of this work is to prove that in a natural setting,
limit as the number of summands tends to infinity and differentiation
can be permuted with asymptotic expansions.
It is however a daunting task --- if not hopeless --- to prove a general 
theorem which could cover all possible cases. In particular, it
is useful to remember that even for the convolution of two heavy
tail distributions, there is no known second order formula that
covers every possible case; however, we will provide some formulas or 
techniques which cover most distributions of interest in applications,
if not all.

\bigskip

\subsection{Asymptotic scales.}
To conclude this introductory section, we present a few facts on terminology
and notations 
related to asymptotic expansions. We refer to Olver (1974) for a complete
account.

Let $N$ be either a positive integer or $+\infty$.
A family of functions $e_i$, $0\leq i<N$,
is called an asymptotic scale if whenever $i+1$ is less than $N$, the
asymptotic relation $e_{i+1}=o(e_i)$ --- i.e. $e_{i+1}(t)=o\bigl( e_i(t)\bigr)$
as $t$ tends to infinity --- holds. We say that a function $f$ has an
$(n+1)$-terms asymptotic expansion in the scale $(e_i)_{0\leq i<N}$ if there
exist real numbers $f_i$, $0\leq i\leq n$ such that
$$
  f(t)=\sum_{0\leq i\leq n} f_i e_i(t) + o\bigl( e_n(t)\bigr)
$$
as $t$ tends to infinity. We then write
$$
  f \sim \sum_{0\leq i\leq n} f_i e_i \, .
$$
When $n$ vanishes, this notation agrees with the more usual one when one writes
$f\sim g$ to mean that $f/g$ tends to $1$ at infinity.
Sometimes, it will be more convenient to index asymptotic scale by an
ordered set, and an obvious variation of the definition
will be understood.

However, this 
definition is insufficient for our purpose. The next step, following
Olver (1974), is to introduce the notion of generalized asymptotic expansion
with respect to an asymptotic scale $e_i$, $0\leq i < N$. If there are 
functions $\phi_s$, $0\leq s\leq n$ such that for any nonnegative integer
$k$ at most $n$,
$$
  f(t)=\sum_{0\leq i\leq k} \phi_i(t) + o\bigl( e_k(t)\bigr) \, ,
$$
we say that $\sum_{0\leq i\leq n} \phi_i(t)$ is a ($n+1$-terms) generalized 
expansion with respect to the scale $e_i$, $0\leq i <N$, and the developing
functions $\phi_i$. We write
$$
  f\sim \sum_{0\leq i\leq n} \phi_i \qquad (e_n) \, ,
$$
or simply
$$
  f\sim \sum_{0\leq i\leq n} \phi_i 
$$
if the scale is understood. Note that there is no a priori understanding
that $e_k=o(\phi_i)$ for any $i$ between $0$ and $n$; in theory, it could be
that the only information conveyed in the expansion is that $f=o(e_k)$.

It is important to note that once the asymptotic scale is chosen, the 
asymptotic expansion of a given function in that scale, if it exists, is 
unique. However, a generalized expansion may not be unique, even if one
fixes the functions $\phi_i$.

One should also be aware that the dependence of an asymptotic expansion
on the asymptotic scale may be an important matter. For instance,
write $\overline\Phi$ for the tail of the standard normal distribution and
consider the asymptotic scale $e_i=\overline\Phi^i$, $i\in\NN^*$. 
Define the function $f(t)=\overline\Phi(t) + \overline\Phi^2(t)$.
One has the asymptotic expansion
$$
  f\sim e_1 +e_2 \, .
$$
Clearly, this approximation cannot be improved.
On the other hand, if one decides to use the asymptotic scale $\tilde e_i(t)
=\overline\Phi(t)/t^i$, $i=0,1,2,\ldots$,
one has for any nonnegative $n$
$$
  f\sim \tilde e_0\qquad (\tilde e_n) \, .
$$
In this asymptotic expansion the part of $f$ coming from $e_2$ is not taken 
into consideration, despite the fact that in the
chosen asymptotic scale, the expansion has as many terms as one wishes.
Finally, if one decides to use the asymptotic scale 
$\overline e_i(t)=e^{-t^2/2}/t^i\sqrt{2\pi}$, $i\geq 1$, then
$$
  f\sim \overline e_1+\sum_{i\geq 1} (-1)^i {(2i-1)!\over 2^{i-1}(i-1)!} 
  \overline e_{2i+1} 
  \, .
$$
Again, in the asymptotic scale $\overline e_i$, $i\geq 1$, this asymptotic
expansion has as many terms as one wishes though the term in $e_2$ in $f$ is
not taken into consideration. The term $e_2$ is effectively $0$ at every
level of scale with respect to the asymptotic scale $(\tilde e_i)_{i\geq 0}$
or the asymptotic scale $(\overline e_i)_{i\geq 1}$.
Moreover, this last asymptotic expansion
is a divergent series, which gives a useless approximation if one does
not truncate it.

We highly recommend the short first chapter of Olver's (1974) book,
where many simple examples are discussed; very useful, not so useful and
totally useless expansions are shown; and the pros and cons of asymptotic
expansions are well presented.

\bigskip


\section{Main result.} 

To state our main result, we need to introduce
two objects: the class of distribution functions we are dealing with and some
differential operators which play a role in writing the asymptotic
expansions. This is the purpose of the next two subsections. The main
result of this paper is given in the third subsection. 

Let us now introduce some notation and conventions.

We write $\II\{\cdot\}$ to denote the indicator function of a set $\{\cdot\}$.

Throughout this paper we let $\M_c$ denote the multiplication operator on
distribution functions corresponding to the multiplication of random 
variables by $c$. Hence, if $F$ is a distribution function and $X$ is 
distributed according to $F$, 
$$
  \M_c F(t) = P\{\, cX\leq t\,\}
  =\cases{ F(t/c) & if $c>0$, \cr
           \oF(t/c-) & if $c<0$, \cr
           \II\{\, t\geq 0\,\} & if $c=0$.\cr}
$$
In the sequel, all distribution functions will be assumed continuous,
thereby obviating the need for left limit notation.
Note that for $c$ positive, $\overline{\M_c F}=\M_c\oF$, while for $c$
negative, $\overline{\M_c F}=\M_{-c}\overline{\M_{-1}F}= F(\cdot/c)$. For a
general function $h$ and positive $c$, we define also $\M_c h$ to be 
$h(\cdot/c)$.

\bigskip

\subsection{The Laplace characters.}%
Let $\D^k$ be the differential operator defined by $\D^kh(x)=h^{(k)}(x)$ for 
any $k$ times differentiable function $h$. As usual, we set $\D^0$ to be the 
identity and $\D^1=\D$.

Recall that a linear differential operator with constant coefficients
is a finite sum $\sum_{0\leq i\leq m}p_i\D^i$. This is a polynomial in
$\D$, and the order of the operator is the degree of the polynomial.
Thus, $p_m$ is nonzero if and only if the previous differential
operator is of order $m$. We write $\RR_m[\,\D\,]$ the set of all linear
differential operators with constant coefficients and order at most $m$.
It is naturally endowed with a vector space structure, corresponding
to that of the polynomials. We will be mostly interested in a
ring structure though. In $\RR_m[\,\D\,]$, we define a composition as
follows. If
$$
  p=\sum_{0\leq i\leq m} p_i\D^i
  \hbox{\qquad and\qquad}
  q=\sum_{0\leq i\leq m} q_i\D^i\, ,
$$
are in $\RR_m[\,\D\,]$, we set
$$
  p\circ q=\sum_{0\leq i\leq m} \Bigl( \sum_{0\leq j\leq i}p_jq_{i-j}\Bigr)\D^i
  \, .
$$
In other words, we use for the ring structure of $\RR_m[\,\D\,]$ that of the 
quotient ring of the polynomials in $\D$ modulo the ideal generated by
$\D^{m+1}$.

In the following, we use the notation $\mu_{F,k}$ for the
moment of order $k$ of a distribution function $F$, that is $\int x^k \d F(x)$.
In particular, $\mu_{F,0}$ is $1$.

\bigskip

\Definition
{ Let $F$ be a distribution function and $m$ be an integer such that $F$ has 
  a finite $m$-th absolute moment. The $m$-th Laplace character of $F$ is 
  the element of $\RR_m[\,\D\,]$ defined by
  $$
    \L_{F,m}=\sum_{0\leq k\leq m} {(-1)^k\over k!} \,\mu_{F,k}\,\D^k \, .
  $$%
}%
When $m$ is clear from the context, we may simply write $\L_F$.
Observe that $\L_{\delta_0}$ is the identity.
The origin of the term `Laplace character' will be explained after our next
proposition.

\Proposition%
{
  \label{LaplaceRepresentation}%
  If $F$ and $G$ are two distribution functions with finite absolute moment
  of order $m$, then $\L_{F,m}\circ \L_{G,m}=\L_{F\star G,m}$.
}

\bigskip

\Proof
The binomial formula shows that
$$
  \sum_{0\leq i\leq k} {\mu_{F,i}\over i!} {\mu_{G,k-i}\over (k-i)!}
  = {\mu_{F\star G,k}\over k!} \, .
$$
The conclusion follows in an easy way from the definition of the
composition.\hfill$\qed$

\bigskip

The previous proposition asserts that the Laplace characters
are morphisms of semigroups. Since the maps $F\mapsto \mu_{F,k}$
are linear, the Laplace characters are representations of
convolution algebras in $\RR_m[\,\D\,]$. One may wonder if other
representations exist, that depend on the distribution through some other set 
of characteristics, say for instance, moments of fractional
order or quantiles. In fact, the answer is negative. By a recent result
of Mattner (2004), the mapping of a distribution with $m$ moments to its
first $m$ cumulants is universal among all continuous homomorphisms
of the convolution algebra of distributions with $m$ moments into Hausdorff
topological groups. This implies in particular that one can express the
Laplace characters in terms of cumulants.

We can now explain our choice of the name `Laplace
character'. Formally, $\L_{F,m}$ is the $m$-th Taylor polynomial of
the Laplace transform of $F$ where we substitute the operator $\D$ for
the variable.  Hence the `Laplace' part of the name. Equivalently, we
could define $\L_{K,m}$ to be $Ee^{-Y\D}$ in $\RR_m[\,\D\,]$, where
$Y$ has distribution function $K$ --- note that with respect to
Mattner's (2004) result, it is clear from this definition that the
Laplace characters can be expressed in terms of cumulants.  Moreover,
Proposition {\LaplaceRepresentation} shows that the Laplace character
obeys a property very similar to the Laplace transform, namely that
this trivializes the convolution into a product of polynomials (modulo
an ideal here).  This property is very similar to that defining the
characters of a group, hence the name. There is an other reason for
the use of the word `character' which has to do with the fact that
$\L_{F,m}$ can be viewed as a formal sum of differential operators,
which reminded us about the Chern character in algebraic topology. And
we also observe that the representation of convolution algebra is
pretty much the same as the representation of modules realized by the
Chern character, our module being over the field of real numbers.

One needs to be careful with one subtlety, namely that our composition
in $\RR_m[\,\D\,]$ is not the usual composition of operators acting on
functions. If $T_1$ and $T_2$ are two operators, we write $T_2T_1$ for
their usual composition, that is $(T_2T_1)h=T_2(T_1h)$ for any
smooth function $h$. We then see that if $i+j>m$, we have $\D^i\D^j=
\D^{i+j}$ while $\D^i\circ \D^j=0$. However, $\D^i\D^j$ is $\D^i\circ \D^j$
modulo an operator in the ideal generated by $\D^{m+1}$. This is in
fact important (see subsection \fixedref{5.2}).
We could avoid this subtlety by defining operators acting on an infinite
sequence of tuples of the form $(\oF_i,\oF_i^{(1)},\ldots, \oF_i^{(m)})$,
$i\geq 1$. It is really a matter of taste, and we feel that the framework
proposed here is somewhat more convenient for our purposes.

For computational purposes, it is of interest to note that operators
in $(\RR_m[\,\D\,],\circ)$ which are not in the ideal generated by $\D$
--- in other words, those which, as a polynomial in $\D$, have a nonzero
constant term --- can be inverted. In particular, the Laplace characters
can be inverted. To write an explicit formula for the inverse, recall
that an ordered partition of length $k$ of an integer $n$ is a
$k$-tuple of positive integers $p=(p_1,\ldots , p_k)$ such that
$p_1\geq p_2\geq \cdots \geq p_k>0$ and $p_1+\cdots + p_k=n$. For such 
partition of length $k$, it is convenient for what follows to agree
on the notation $p_i=0$ for all $i$ larger than $k$. We write $\calP(m,n)$ 
for the set of all ordered partitions of length at most $m$ of $n$.
For such partition $p$, we write $\Delta p$
for $(p_1-p_2,p_2-p_3,\ldots ,p_m-p_{m+1})$. Also for a tuple 
$k=(k_1,\ldots ,k_m)$ and an integer $q$ we write 
${q\choose k}={q \choose k_1\ldots k_m}$ the multinomial type coefficient 
$q!/(k_1!\,\cdots\, k_m!)$.

\Proposition%
{\label{InverseLaplaceChar}%
  In $\RR_m[\,\D\,]$, we have the inversion formula
  $$
    \L_{F,m}^{-1}
    =\sum_{0\leq n\leq m} \Bigl( \sum_{p\in\calP(m,n)}
    (-1)^{n+p_1}{p_1\choose \Delta p} \prod_{1\leq k\leq m}
    \Bigl({\mu_{F,k}\over k!}\Bigr)^{(\Delta p)_k} \Bigr) \D^n \, .
  $$
}

\Proof In this proof we write $i=(i_1,\ldots ,i_m)$ a generic tuple
of nonnegative integers, and $|i|=i_1+\cdots +i_m$. In the space of power
series in $t$,
$$\displaylines{\qquad
  \Bigl(1+\sum_{1\leq k\leq m} {(-1)^k\over k!} \mu_{F,k} t^k\Bigr)^{-1}
  \hfill\cr\hfill
  \eqalign{
  {}={}& \sum_{j\geq 0} (-1)^j\Big( \sum_{1\leq k\leq m} {(-1)^k\over k!}
         \mu_{F,k}t^k\Bigr)^j \cr
  {}={}& \sum_{j\geq 0} (-1)^j \sum_{|i|=j}{j\choose i}\prod_{1\leq k\leq m}
         \Bigl( {(-1)^k\over k!}\mu_{F,k}t^k\Bigr)^{i_k} \, . \cr
  }
  \qquad\cr}
$$
Consequently, in $\RR_m[\,\D\,]$, the inverse of $\L_{F,m}$ is
$$
  \sum_{0\leq j\leq m} (-1)^j \sum_i {j\choose i} 
  \biggl(\prod_{1\leq k\leq m}\Bigl( {(-1)^k\over k!}\mu_{F,k}\Bigr)^{i_k} 
  \biggr)
  \D^{i_1+2i_2+\cdots +mi_m}
  \, ,
$$
where the sum over $i$ is over all tuples $(i_1,\ldots , i_m)$ such 
that $|i|=j$ and $i_1+2i_2+\cdots+mi_m\leq m$.
Set $p_{m+1}=0$ and $p_l=i_l+p_{l+1}$ for $l=m,m-1,\ldots , 1$. Thus,
$p_m=i_m$, $p_{m-1}=i_m+i_{m-1}$, $\ldots$, $p_1=i_m+\cdots + i_1$.
Set $p=(p_1,\ldots , p_m)$. We see that $\sum_{1\leq k\leq m}ki_k=n$
if and only if $p$ is an ordered partition of $n$, and the correspondance
between $p$ and $i$ is one-to-one. Moreover, $|i|=p_1$. Therefore
$$
  \L_{F,m}^{-1}
  =\sum_{n\leq m} \Bigl( \sum_{p\in\calP(m,n)}(-1)^{p_1}
  {p_1\choose\Delta p}\prod_{1\leq k\leq m} 
  \Bigl( {(-1)^k\over k!}\mu_{F,k}\Bigr)^{(\Delta p)_k} \Bigr) \D^n \, ,
$$
which is the inversion formula.\hfill$\qed$

\bigskip

To use Proposition \InverseLaplaceChar, one may need to generate
all the partitions of a given integer. We refer to
Nijenhuis and Wilf (1978) or Stanton and White (1986) for algorithms
related to that matter.

Our next result is an other inversion formula for the Laplace characters.
It is not as explicit as the previous one, because it is written in
the ring $(\RR_m[\,\D\,],\circ)$. However, we will make use of it. To write 
this inversion formula, for any differential operator $T$ in 
$(\RR_m[\,\D\,],\circ)$ and any nonnegative integer $k$, we define 
$T^{\circ k}$ inductively by $T^{\circ 0}=\Id$ and 
$T^{\circ k}=T\circ T^{\circ (k-1)}$. 

\Proposition
{\label{InverseLaplaceCharB}%
In $(\RR_m[\,\D\,],\circ)$, we have the inversion formula
  $$\L_{F,m}^{-1}
    = \sum_{0\leq k\leq m} (\Id-\L_{F,m})^{\circ k}
    = \sum_{0\leq j\leq m} (-1)^j \L_{F,m}^{\circ j} \sum_{j\leq k\leq m} 
      {k\choose j} \, .
  $$
}

\Proof Since $\L_{F,m}-\Id$ is in the ideal generated by $\D$, it is nilpotent
of nullity at most $m+1$ in $(\RR_m[D],\circ)$. Consequently,
$$
  \L_{F,m}^{-1} 
  = \bigl( \Id -(\Id -\L_{F,m})\bigr)^{-1}
  = \sum_{0\leq k\leq m} (\Id -\L_{F,m})^{\circ k} \, .
$$
Then, the second equality follows from the binomial formula
$$
  (\Id-\L_{F,m})^{\circ k} 
  = \sum_{0\leq j\leq k} {k\choose j} (-1)^j \L_{F,m}^{\circ j} \, .
  \eqno{\qed}
$$

It is also clear that there are other ways of thinking of the Laplace
characters and their inverses that can be useful in applications.
For instance we can think of $\L_{F,m}$ as the $m$-th Taylor polynomial
of the Laplace transform $\calL_F$ (if it exists) of $F$ evaluated at $\D$. 
Then $\L_{F,m}^{-1}$ is simply obtained by taking the $m$-th Taylor 
polynomial of $1/\calL_F$ and evaluating it in the variable $\D$.

\bigskip

There are other algebraic operations on Laplace characters which are of 
interest. For instance, we will use the Mellin-Stieltjes convolution between 
two distribution functions $F$ and $G$ on the nonnegative half line. It is 
writen $F\convM G$, and defined by
$$
  F\convM G(t)=\int_0^\infty F(t/x) \d G(x) \, .
$$
If $X$ and $Y$ are two independent random variables with respective 
distributions $F$ and $G$, then $F\convM G$ is the distribution of the product
$XY$. Since $E(XY)^j=EX^j EY^j$, the Laplace character of $F\convM G$
is obtained by multiplying the Laplace characters of $F$ and $G$ 
coefficient-wise, that is
$$
  \L_{F\convMss G,m}
  =\sum_{0\leq j\leq m} {(-1)^j\over j!} \mu_{F,j}\mu_{G,j} \D^j \, .
$$

Yet, another algebraic operation that one could consider, and which we will
not use in this paper, is that of differentiating formally the Laplace
characters as polynomials in $\D$. Assuming that $\mu_{F,k}$ exists and does
not vanish, writing $\d G(x)=\mu_{F,k}^{-1}x^k \d F(x)$,
we see that $\L_{G,m}$ is $(-1)^k/\mu_{F,k}$ times the $k$-th derivative with
respect to the variable $\D$ of $\L_{F,m+k}$. 

Similarly, if $F$ is a distribution on the nonnegative half line and 
$G(t)=\mu_{F,1}^{-1}\int_0^t \oF (x) \d x$, then 
$$
  \L_{G,m} = -\mu_{F,1}^{-1}(\L_{F,m+1}-\Id)/\D \, .
$$

The interest of such algebraic formulas is to produce a form of operational 
calculus to derive tail expansions. This will be clear in section \fixedref{3}.

\medskip

We conclude this subsection with two lemmas. The first one is related to
Proposition {\LaplaceRepresentation} while the second one describes the 
behavior of the Laplace characters under some multiplications.

\Lemma
{\label{LaplaceBinomial}%
  If $H$ and $K$ are two distribution functions with finite $m$-th absolute 
  moments, then
  $$
    \L_{K\star H,m}
    =\sum_{0\leq j\leq m} {(-1)^j\over j!}\mu_{H,j}\L_{K,m-j}\D^j 
    \, .
  $$
}

\Proof The right hand side of the equality is
$$
  \sum_{0\leq j\leq m}{(-1)^j\over j!}\mu_{H,j} 
  \sum_{0\leq l\leq m-j} {(-1)^l\over l!}
  \mu_{K,l} \D^{j+l} \, .
$$
Setting $s=l+j$, it is
$$\displaylines{\qquad
  \sum_{s,j\geq 0} 
  \II\{\, j\leq s\leq m\,\} {(-1)^s\over s!}{s\choose j}
  \mu_{H,j}\mu_{K,s-j} \D^s
  \hfill\cr\hfill
  {}=\sum_{0\leq s\leq m} {(-1)^s\over s!} \mu_{K\star H,s} \D^s \, ,
  \qquad\cr}
$$
which gives the conclusion.\hfill$\qed$

\Lemma
{\label{LaplaceCharMul}%
  If $K$ is a distribution function whose $m$-th moment is finite, then
  $$
    \L_{\M_\lambda K,m}\M_\lambda = \M_\lambda \L_{K,m} \, . 
  $$
}

\Proof It follows from the equalities
$$\eqalign{
  \L_{\M_\lambda K,m}(\M_\lambda h)
  &{}=\sum_{0\leq j\leq m} {(-1)^j\over j!} \mu_{\M_\lambda K,j} 
   \D^j(\M_\lambda h) \cr
  &{}=\sum_{0\leq j\leq m} {(-1)^j\over j!} \lambda^j\mu_{K,j}\lambda^{-j}
    \M_\lambda\D^j h \, , \cr
  }
$$
and the equality of the last term with $\M_\lambda \L_{K,m}h$.\hfill$\qed$

\bigskip

\subsection{Smoothly varying functions of finite order.}%
Among the regularly varying functions, the normalized ones will be of 
importance
to us. Following Bingham, Goldie and Teugels (1989, \S 1.3.2), we say that
a function $g$ defined on $[\, a,\infty)$ is a normalized regularly
varying function with index $\rho$ if it has the representation
$$
  g(t)=t^\rho c\exp\int_a^t {\epsilon (u)\over u} \d u \, ,
  \eqno{\equa{KaramataRepresentation}}
$$
where $\epsilon (\cdot )$ is a function converging to $0$ at infinity.

For any real number $x$ and any positive integer $k$, we write $(x)_k$
for $x(x-1)\cdots (x-k+1)$. We also set $(x)_0=1$. If $h$ is a function,
$h^{(k)}$ is the $k$-th derivative of $h$ if it exists, with the usual
convention $h^{(0)}=h$.

Recall that a function $h$ defined in some neighborhood
of infinity is smoothly varying with index $-\alpha$ if it is
ultimately infinitely differentiable and
$$
  \limt t^k h^{(k)}(t)/h(t) = (-\alpha )_k 
  \neq
$$
for every integer $k$ (see Bingham, Goldie and Teugels, 1989, \S 1.8).
The set of all smoothly varying functions with index $-\alpha$ is
written $SR_{-\alpha}$. The set of smoothly varying functions of 
a given finite order which we will define contains $SR_{-\alpha}$
and can be thought of as a Sobolev space in the framework of
regular variation. Before proceeding, it is useful to remark that for $k=1$
relation {\preveq} forces $h$ to be normalized 
regularly varying of index $-\alpha$
(see Bingham, Goldie and Teugels, 1989, \S 1.8). This and {\preveq}
forces $h^{(k)}$ to be regularly varying of index $-\alpha-k$.

\bigskip

\Definition
{ Let $m$ be a positive integer. A function $h$ is smoothly varying
  of index $-\alpha$ and order $m$ if it is ultimately $m$-times
  continuously differentiable and $h^{(m)}$ is regularly
  varying of index $-\alpha-m$.
  We write $SR_{-\alpha,m}$ for the set of all such functions.
}

\bigskip

From Karamata's theorem, we see that a function $h$
in $SR_{-\alpha,m}$ satisfies {\preveq} for any $0\leq k\leq m$.
Note that if $m_1$ is at most $m_2$ then $SR_{-\alpha,m_2}$ is included
in $SR_{-\alpha,m_1}$.

Also, if $h$ belongs to $SR_{-\alpha,m}$, then all $h^{(k)}$'s for $k=0,1,
\ldots , m-1$ are normalized regularly varying. Consequently, for any such $k$,
the function $t\mapsto t^\sigma |h^{(k)}(t)|$ is ultimately decreasing 
(respectively increasing) when $\sigma<\alpha+k$ 
(respectively $\sigma>\alpha+k$) --- see Bingham, Goldie and Teugels, 1989,
Theorem 1.5.5.

To define $SR_{-\alpha,s}$ for a real number $s$, we write $\delta_x$ for the
point mass at $x$. Thus if $h$ is a function, $\delta_x h=\int h\d\delta_x
=h(x)$. We then define the operator
$$
  \Delta_{t,x}^r=\sign (x) {\delta_{t(1-x)}-\delta_t\over |x|^r\delta_t} \, .
$$
In other words, for any function $g$ we set
$$
  \Delta_{t,x}^r g = \sign(x) {g\bigl(t(1-x)\bigr)-g(t)\over |x|^r g(t)} \, .
$$

\Definition
{ Let $s$ be a positive real number. Write $s=m+r$ where $m$ is the
  integer part of $s$ and $r$ is in $[\, 0,1)$. A function $h$ is
  smoothly varying of index $-\alpha$ and order $s$ if it belongs
  to $SR_{-\alpha,m}$ and
  $$
    \lim_{\delta\to 0}\limsupt \sup_{0<|x|\leq\delta} |\Delta_{t,x}^r h^{(m)}|
    = 0 \, .
    \eqno{\equa{Deltar}}
  $$
}

\noindent Note that $m=m+0$, and that if $h$ belongs to $SR_{-\alpha,m}$ then 
{\Deltar} holds for $r=0$ thanks to the uniform convergence on compact subsets
of $(0,\infty)$ theorem (Bingham, Goldie, Teugels,1989, Theorem 1.2.1)

When $s$ is smaller than $1$, the set $SR_{-\alpha,s}$ is closely 
related to that of all
the functions satisfying the Lipschitz condition $[\,\D_s\,]$ of Borovkov
and Borovkov (2002), while for positive $m$ it is related to their
condition $[\,\D_m\,]$.

Our next result shows that these spaces are nested.

\Proposition
{
  If $r$ is at most $s$, then $SR_{-\alpha,s}\subset SR_{-\alpha,r}$.
}

\bigskip

\Proof If $r$ or $s$ is an integer, the result is obvious, hence we assume
that both are not integers.
Write $r=m+\rho$ and $s=n+\sigma$ with $m$ and $n$ integers
and $\rho$, $\sigma$ positive and less than $1$. Let $h$ be a function in 
$SR_{-\alpha,s}$. If $n$ is at most $r$, then $m$ and $n$ are equal, and 
$\rho$ is less than $\sigma$. It is 
then plain that $h$ is in $SR_{-\alpha,r}$ since  the function 
$|x|^{\sigma-\rho}$ is bounded in any neighborhood of the origin.

Assume that $n$ is larger than $r$. Then $n$ is at least $m+1$ and
$$
  \limt t^{m+1} h^{(m+1)}(t)/h(t) = (-\alpha)_{m+1} \, .
$$
Since $m+1$ is at least $1$, the function $h^{(m+1)}$ is regularly varying 
of index $-\alpha-m-1$. We write
$$
  h^{(m)}\bigl( t(1-x)\bigr) - h^{(m)}(t)
  = -\int_{1-x}^1 t h^{(m+1)}(tu) \d u
$$
to obtain
$$
  |\Delta_{t,x}^\rho h^{(m)}|
  \leq |x|^{1-\rho} \sup_{1-x\leq u\leq 1} |th^{(m+1)}(tu)/h^{(m)}(t)| \, .
$$
Using the uniform convergence Theorem (Bingham, Goldie and Teugels, 1989,
Theorem 1.2.1) and that $\rho$ is less than $1$, we conclude that
$$
  \lim_{\delta\to 0} \limsupt\sup_{0<|x|\leq\delta} |\Delta_{t,x}^\rho h^{(m)}| = 0 \, ,
$$
which shows that $h$ belongs to $SR_{-\alpha,r}$.\hfill$\qed$

\bigskip

Before going further, let us make a digression about condition 
{\Deltar}. One may wonder about the analogous condition when $r$ is 
$1$. In fact, as we will explicate next, the condition {\Deltar} with $r=1$
has bearing on several issues, e.g.\ understanding the spaces 
$SR_{-\alpha,\omega}$, the monotone density theorem which is a classical
result in the theory of regular variation, the class of asymptotically
smooth functions introduced in Barbe and McCormick (2005) as well as the
conditions $[\,\D\,]$ in Borovkov and Borovkov (2003).

For clarity of the arguments, if a function $g$ is differentiable, then
$$
  \lim_{x\to 0}\Delta_{t,x}^1 g = -tg'(t)/g(t) \, .
$$
Therefore, if $\beta$ is a positive number and if $g'$ is regularly 
varying with index $-\beta-1$,
$$
  \lim_{t\to\infty} \lim_{x\to 0} \Delta_{t,x}^1 g = \beta \, .
$$
For $r=1$, it is then natural to introduce the condition
$$
  \lim_{\delta\to 0} \limsup_{t\to\infty} \sup_{0<|x|<\delta}
  |\Delta_{t,x}^1 g-\beta| = 0 \, .
  \eqno{\equa{DeltaOne}}
$$

\Proposition{\label{RVequiv}
  The following are equivalent:

  \noindent (i) condition \DeltaOne.

  \noindent (ii) $g$ is normalized regularly varying with index $-\beta$.

  \noindent (iii) $g$ is ultimately absolutely continuous and a version of
  its Radon-Nykodym derivative is regularly varying with index $-\beta-1$.
}

\bigskip

\Proof We write $\dot g$ a Radon-Nykodym derivative of $g$ when it exists.

We start by proving that (iii) implies (ii). If (iii) holds, by Karamata
theorem (see Bingham et al., 1989, Theorem 1.5.11), $\Id\, \dot g/g\sim -\beta$
at infinity. Integrating $\dot g/g$ gives (ii).

Next, we prove that (ii) implies (iii). Consider the Karamata representation
of $g$,
$$
  g(t) = t^{-\beta} c\exp \int_{t_0}^t {\epsilon(u)\over u} \d u \, .
$$
Such a function is ultimately absolutely continuous and a Radon-Nykodym 
derivative is
$$
  \dot g(t) = {g(t)\over t} \Bigl( -\beta + {\epsilon(t)\over\beta}\Bigr) \, .
$$
This implies (iii).

We now prove that (ii) implies (i). Again, the Karamata representation of 
$g$ yields for $t$ large enough and $x$ small enough,
$$
  \Delta_{t,x}^1 g
  = {(1-x)^{-\beta}-1\over x} + (1-x)^{-\beta} {1\over x} 
  \Bigl( \exp\int_t^{t(1-x)} {\epsilon(u)\over u} \d u -1\Bigr) \, .
$$
Let $\eta$ be a positive number and let $t_1$ be large enough so that $g$ is
absolutely continuous on $(t_1,\infty)$ and the absolute value of 
$\epsilon (\cdot)$ is at most $\eta$ on $(t_1,\infty)$. Then, for any $t$
and $t(1-x)$ more than $t_1$,
$$
  -\eta|\log (1-x)| 
  \leq \int_t^{t(1-x)} {\epsilon(u)\over u} \d u
  \leq \eta |\log(1-x)| \, .
$$
For $\delta$ small enough, this implies
$$
  \limsup_{t\to\infty} \sup_{0<|x|<\delta}\, \Bigl| {1\over x} 
  \Bigl( \exp\int_t^{t(1-x)}{\epsilon (u)\over u} \d u -1\Bigr)\Bigr| 
  \leq 2\eta \, .
$$
On the other hand,
$$
  \lim_{\delta\to 0} \limsup_{t\to\infty} \sup_{0<|x|<\delta}\,
  \Bigl| {(1-x)^{-\beta}-1\over x} -\beta\Bigr| = 0 \, .
$$
Therefore, {\DeltaOne} and (i) hold.

We conclude the proof of the Proposition by showing that (i) implies
(ii). We first note that if {\DeltaOne} holds, then $g$ must be
ultimately continuous; otherwise, the supremum in $x$ in {\DeltaOne}
would be infinite for some arbitrary large $t$'s, precluding
{\DeltaOne} to hold. Also, $g$ does not vanish ultimately. So, without
loss of generality, we assume that $g$ is ultimately positive. Let
$\dot g_U$ and $\dot g_L$ be the upper and lower derivatives of $g$,
that is
$$\eqalign{
  \dot g_U(t) &{}=\limsup_{\epsilon\to 0} {g(t+\epsilon)-g(t)\over\epsilon}
               \, , \cr
  \dot g_L(t) &{}=\liminf_{\epsilon\to 0} {g(t+\epsilon)-g(t)\over\epsilon}
               \, . \cr
}
$$
Condition {\DeltaOne} implies that $\Id\, \dot g_U/g$ and $\Id\, \dot g_L/g$ 
have limit $-\beta$ at infinity. Set $f(t)=t^\beta g(t)$. Since $g$ is 
ultimately continuous,
$$
  \dot f_U(t)=\beta t^{\beta-1}g(t) + t^\beta \dot g_U(t) \, ,
$$
and an analogous relation holds for the lower derivative of $f$. This implies
$$
  t\dot f_U(t)/f(t) = \beta + t\dot g_U(t)/g(t)
$$
has limit $0$ at infinity, and similarly for the upper derivative. By Bojanic
and Karamata (1963) --- see Bingham et al.\ 1989, exercice 1.11.8
--- and Theorem 1.5.5 in Bingham et al.\ (1989), this implies that $f$
is normalized slowly varying.\hfill$\qed$

\bigskip

Note that Proposition {\RVequiv} implies that if a function $h$ belongs
to $SR_{-\alpha,m}$ and {\DeltaOne} holds with $h^{(m)}$ in place of $g$ and 
$\beta=\alpha+m$, then the existence of the derivative $h^{(m+1)}$ 
guarantees that $h$ belongs to $SR_{-\alpha,m+1}$. However, without assuming
that $h^{(m+1)}$ exists, {\DeltaOne} would only give existence and 
regular variation of a Radon-Nykodym derivative of $h^{(m)}$.

The equivalence between (ii) and (iii) in Proposition {\RVequiv} shows that
the class of normalized regularly varying functions is a very natural one
when dealing with regular variation of a function and its derivative.

It also follows from Proposition {\RVequiv} that the class of asymptotically
smooth distributions in Barbe and McCormick (2005) is in fact those which
are normalized regularly varying; this was suggested to
us by Jaap Geluk and prompted Proposition {\RVequiv}. It also shows that 
condition [D$_m$] of Borovkov and Borovkov (2003) relates to normalized
regular variation.

Proposition {\RVequiv} is also interesting with respect to the monotone
density theorem (Bingham et al., 1989, Theorem 1.7.2). Indeed, this theorem 
asserts that if a regularly
varying function $g$ with index $-\beta$  has ultimately monotone derivative 
$g'$, then this derivative is regularly varying of index $-\beta-1$; and
moreover, by Karamata's theorem, $\Id\, g'/g\sim-\beta$. This implies that $g$
is normalized regularly varying, and therefore satisfies (ii) of Proposition
{\RVequiv}. In other words, the monotone density theorem can be viewed as a 
particular case of the implication of (iii) by (ii) in Proposition {\RVequiv}.
This concludes our digression on condition {\Deltar}

\medskip

The importance of the classes $SR_{-\alpha,m}$ in this paper stems 
from the nice behavior of their functions with respect to differentiation,
global 
Potter type bounds, and to Taylor formula used asymptotically. 
In order to elaborate on this assertion, recall first that Potter's 
bounds (see, Bingham, Goldie and Teugels, 1989, 
Theorem 1.5.6) assert that if $g$ is a function which is regularly varying
with index $\rho$, then for any $A$ larger than $1$ and any 
$\delta$ positive, there
exists a $t_0$ such that for any $\lambda$ at least $1$ and any
$t$ more than $t_0$,
$$
  A^{-1}\lambda^{\rho-\delta}\leq g(\lambda t)/g(t) \leq A \lambda^{\rho+\delta}
  \, .
$$

This standard Potter inequality already yields an upper bound
on the decay of the scaled derivatives of smoothly varying functions
of fixed order.

\Lemma{%
  \label{BoundDkMF}
  Let $h$ be a smoothly varying function of index $-\alpha$ and 
  order $s$ more than $1$. Let $\epsilon$ be a positive number.
  There exists $t_1$ such that for any $c$ positive at most $1$ and 
  any $t$ at least $t_1$ and any integer $k$ at most $s$,
  $$
    |(\M_ch)^{(k)}(t)|
    \leq 2|(-\alpha)_k| c^{\alpha-\epsilon} t^{-k}|h(t)| \, .
  $$
}

\Proof
We have $(\M_ch)^{(k)}(t)= c^{-k} h^{(k)}(t/c)$.
Since $h$ is in $SR_{-\alpha,s}$ with $s$ at least $k$, there exists
$t_0$ such that for any $t$ at least $t_0$, 
$$
  |t^k h^{(k)}(t)/h(t)|\leq \sqrt 2 \, |(-\alpha)_k| \, .
$$
Since $c$ is positive and at most $1$, 
the inequality $t\geq t_0$ implies $t/c\geq t_0$.
Therefore, for $t$ at least $t_0$,
$$
  c^{-k}|h^{(k)} (t/c)|
  \leq \sqrt 2|(-\alpha)_k| t^{-k} |h(t/c)| \, .
$$
Now, if $t$ is larger than some $t_1$ independent of $c$, Potter's bounds give
$|h(t/c)|\leq \sqrt2 c^{\alpha-\epsilon}|h(t)|$. This implies the 
result.\hfill$\qed$

\bigskip

The following lemma shows that for normalized regularly varying functions,
in particular for smoothly varying functions of fixed order at least 
$1$
Potter's bounds can be improved by taking $A=1$.
Although this may seem a minor feature, this improvement is essential for
our purpose.

\Lemma{
  \label{Potter}
  Let $g$ be a normalized regularly varying function with index
  $\rho$. Let $\delta$ be a positive number. There exists $t_2$ such
  that for any $\lambda$ larger than $1$ and any $t$ at least $t_2$,
  $$
    \lambda^{\rho-\delta} \leq g(t\lambda )/g(t) \leq \lambda^{\rho+\delta} \, .
  $$
}

\Proof Let $\epsilon(\cdot )$ be as in {\KaramataRepresentation}. 
Let $t_2$ be at least
$a$ and such that $\sup_{t\geq t_2}|\epsilon (t)|\leq \delta$. Then,
$$
  {g(\lambda t)\over g(t)}
  = \lambda^\rho \exp\Bigl( \int_t^{t\lambda} {\epsilon (u)\over u} \d u\Bigr)
  \cases{ \leq \lambda^\rho\exp\bigl(\delta\int_t^{t\lambda}u^{-1}\d u\bigr) 
          &\cr
          \noalign{\vskip 2pt}
          \geq \lambda^\rho\exp\bigl(-\delta\int_t^{t\lambda}u^{-1}\d u\bigr)
          \, , &\cr
        }
$$
which yields the result.\hfill$\qed$

\bigskip

To connect Taylor formula and asymptotic expansions, we set
$$
  \oDelta_{\tau,\delta}^q(h)
  =\sup_{t\geq\tau}\sup_{0<|x|\leq\delta} |\Delta_{t,x}^q h| \, .
$$

\Proposition
{ 
  \label{Taylor}%
  Let $r$ be in $[\, 0,1\,]$. If
  $h$ is $m$ times differentiable, then for any positive $t$ and 
  any $u$,
  $$
    \Bigl|h(t+u)-\sum_{0\leq j\leq m} {u^j\over j!} h^{(j)}(t)\Bigr|
    \leq {|u|^{m+r}\over t^r} {|h^{(m)}(t)|\over m!}\, 
    \oDelta_{t,|u|/t}^r(h^{(m)})
    \, .
  $$
}

\Proof The Taylor-McLaurin formula yields
$$
  h(t+u)
  =\sum_{0\leq j\leq m} {u^j\over j!} h^{(j)}(t)
  + {u^m\over m!} \bigl( h^{(m)}(t+\theta_{t,u}u)-h^{(m)}(t)\bigr)
$$
with $\theta_{t,u}$ between $0$ and $1$. The equality
$$
  |h^{(m)}(t+\theta_{t,u}u)-h^{(m)}(t)|
  = \theta_{t,u}^r t^{-r}  |u|^r |h^{(m)}(t)| |\Delta_{t,-\theta_{t,u}u/t}^r h^{(m)}| 
$$
implies the result.\hfill$\qed$

\bigskip

Observe that Lemma {\BoundDkMF} asserts that if $h$ is smoothly varying
of order $s=m+r$, then $h^{(m)}(t)$ is 
of order $t^{-m}h(t)$. Thus, the remainder term in Proposition {\Taylor}
is expected to be asymptotically rather small and certainly of smaller
order of magnitude than any of the $h^{(j)}(t)$ involved in the inequality
of the Proposition. In other words, for functions
in $SR_{-\alpha,s}$ we can relate the local character of Taylor's formula
to an asymptotic expansion of the translation $u\mapsto t+u$ acting on those
functions. The Proposition implies that if $h$ is in $SR_{-\alpha,s}$, with
$s$ larger than $m$, then $\sum_{0\leq j\leq m} (u^j/j!) h^{(j)}(t)$
is an asymptotic expansion of $h(t+u)$ as $t$ tends to infinity.
Using the translation $\tau_x$ defined on functions by $\tau_x h(t)=h(t+x)$,
Proposition {\Taylor} asserts that $\tau_u$ is approximately the Laplace
character $\L_{\delta_u,m}$ when read on the tail of smoothly varying 
functions.

Before moving to the next section, we make a remark concerning the
behavior of the operator $\Delta_{t,x}^r$ with respect to composition
by differentiation and multiplications $\M_c$. This remark does not
have much intrinsic values, but we will need it during our proof.

\Lemma{%
  \label{DeltaDM}%
  For any positive $c$, the equality 
  $\Delta_{t,x}^r \D^j \M_c=\Delta_{t/c,x}^r\D^j$ holds.
}

\bigskip

\Proof
We have 
$$
  \D^j\M_ch(t) = {\d^jh(t/c)\over \d t^j} = c^{-j} h^{(j)}(t/c) \, .
$$
Therefore,
$$
  \Delta_{t,x}^r \D^j \M_c h
  = {h^{(j)}\bigl( (t/c)(1-x)\bigr) -h^{(j)}(t/c)\over |x|^r h^{(j)}(t/c)}
  = \Delta_{t/c,x}^r \D^j h \, .
  \eqno{\qed}
$$

\bigskip

\subsection{Asymptotic expansion for infinite convolution.}%
To state our main result, recall that the $\ell_p$ norm of a sequence
$c=(c_i)_{i\in\ZZ}$ is
$$
  |c|_p=\Bigl(\sum_{i\in\ZZ} |c_i|^p\Bigr)^{1/p} \, .
$$
This defines a norm only if $p$ is at least $1$. Nevertheless, we will
still use the same notation with the same meaning when $p$ is less than $1$.
When $p$ is infinite, $|c|_\infty$ is simply the supremum of the $|c_i|$'s.
Given three nonnegative numbers $\alpha$, $\gamma$ and $\omega$, we define
$$
  N_{\alpha,\gamma,\omega}(c)
  =|c|_{\gamma ({\ss\alpha\over\ss\alpha+\omega}\wedge {\ss 1\over\ss 2})}
  \vee 2^{\alpha/(\alpha+\omega)}|c|_\infty
  \, .
$$
This may or may not be a norm, according to whether or not $\gamma \alpha/
(\alpha+\omega)$ and $\gamma/2$ are at least $1$. For the values of $\gamma$
that we will use, that is $\gamma$ positive and less than $1$, this is not a 
norm.

The next theorem is the main result of the paper. It establishes a 
generalized expansion for some infinite weighted convolutions with respect
to the asymptotic scale $\Id^{-i}\oF$, $i\geq 0$. To state
this result, we need further notation.

If $F$, $G$, $H$ are three distribution functions, with $H=F\star G$, we 
write $F=H\natural G$, the division of $H$ by $G$. 

Note that if $c_i$ is negative, then the upper tail of $c_iX_i$ is driven by
the lower tail of $X_i$. This induces some complications in stating a result
for the tail of $\<c,X\>$ because one needs stronger assumptions when
the signs of the constants and the random variable can be arbirary.
This explains the formulation of the next result. It does not
cover all the possible variations, but seems to give what is needed in most
applications. Other cases often may be obtained by simple changes in the proof.

For two functions $f$ and $g$, we write $f\asymp g$ to mean that the
ratio $f/g$ is ultimately bounded away from $0$ and infinity.

In the following theorem, we write $G_c$ for the distribution 
of $\sum_{i\in\ZZ} c_i X_i$. Note that $G_c\natural \M_{c_i}F$ is the
distribution of that infinite series with the $i$-th term removed.

\Theorem
{\label{MainTheorem}%
Let $\omega$ be at least $1$. Let $F$ be a distribution function with $\oF$ in 
$SR_{-\alpha,\omega}$. Let 
$m$ and $k$ be two integers such that $m$ is smaller than $\alpha$, and $m+k$ 
is smaller than $\omega$. Furthermore, 
if some $c_i$'s are negative 
assume either that $\overline{\M_{-1}F}$ is also in $SR_{-\alpha,\omega}$ and
$\overline{\M_{-1}F}\asymp \oF$ or that $F$ vanishes in a neighborhood of
$-\infty$.

Let $\gamma$ be a positive number less than $\omega-m-k$ and $1$.
Then, there exists a function $\eta (\cdot )$ converging
to $0$ at infinity and a real number $t_0$ such that for any $t$ at least
$t_0$, for any sequence $c$ with $N_{\alpha,\gamma,\omega}(c)\leq 1$,
$$
  \Bigl|\oG_c^{(k)}(t)
  -\sum_{i\in\ZZ}\L_{G_c\natural \M_{c_i}F,m}(\overline{\M_{c_i}F})^{(k)} (t)
  \Bigr|
   \leq t^{-m-k}\oF(t)\eta(t) \, .
$$
In particular for any sequence $c$ with $N_{\alpha,\gamma,\omega}(c)$ finite,
$$
  \oG_c^{(k)}
  \sim \sum_{i\in\ZZ}\L_{G_c\natural \M_{c_i}F,m}(\overline{\M_{c_i}F})^{(k)}
  \qquad (\Id^{-m-k}\oF)\, .
$$
}

This theorem calls for several remarks.

\Remark By allowing a certain uniformity in the Potter type bounds and in the
asymptotic smoothness assumption, a degree of uniformity in our result
with respect to the underlying distribution $F$ can be achieved. However,
we choose not to develop such a refinement since it would entail a greater
level of technicality which may be distractive to the main aim of the paper.

\Remark When some $c_i$'s are negative, the assumption 
$\oF\asymp \overline{\M_{-1}F}$ is not necessary. It is
quite clear from the proof how to modify the statement if we assume 
that $\overline{\M_{-1}F}$ is in some
$SR_{-\beta,\omega'}$. However, the statement is far nicer when either both
tails are comparable, or when $F$ is supported by the nonnegative half line.

\Remark When the $c_i$'s are nonnegative, one does not really need to 
assume $\overline{\M_{-1}F}=O(\oF)$. Certainly, the result holds provided 
only that the $m$-th absolute moment of the distribution is finite. 
This technical point appears in Lemma \fixedref{6.3.1.}

\Remark If either $k$ or $m$ are positive and their sum is less than
$\omega$, then $\omega$ is more than $1$, and the first sentence in Theorem
{\MainTheorem} is not needed. Assuming that $\omega$ is at least $1$ is
only needed when both $k$ and $m$ vanish; it is assumed only to
ensure that $\oF$ and possibly $\overline{M_{-1}F}$ are normalized regularly
varying.

\Remark The analoguous statement on the lower tail holds as well.

\Remark Since the Laplace characters are differential operators with constant
coefficients, they commute with the derivative $\D$. Consequently, Theorem
{\MainTheorem} implies that for a certain class of distributions and sequences,
taking limit with
respect to the number of nonvanishing terms in the sequence $c$ and 
differentiation can be permuted in an asymptotic expansion. It is quite
striking that this can be done with some uniformity with respect to a rather
large class of sequences.

\Remark It will be clear in the sequel that the 
formulation of Theorem {\MainTheorem} is well suited for applications. 
One may still wonder if one can relax its assumptions on $F$. The example
of first order autoregressive processes suggests that some refinements may
be given, but not in a fundamental way. Indeed, consider the sequence
$c_i=a^i$ for nonnegative integers $i$, and $c_i=0$ for negative $i$. Thus,
$G_c$ is the distribution of $Y=\sum_{i\geq 0} a^i X_i$. For any integer
$k$, set
$$
  Z_{k,i}=\sum_{0\leq j<k} a^j X_{ki+j}\, .
$$
For fixed $k$, as $i$ runs through the nonnegative integers, these
random variables are independent and equidistributed. Moreover,
$Y=\sum_{i\geq 0} a^{ki}Z_{k,i}$. Hence, if the distribution function
of $Z_{k,i}$ fulfills the assumptions of Theorem {\MainTheorem}, one
can obtain the tail expansion of $G_c$, even though the distribution
of $X_i$ may not satisfy the assumptions of Theorem
{\MainTheorem}. This may happen for instance if the distribution of
$X_i$ is a mixture containing point masses. Note however that the tail
expansion will be expressed in term of the distribution of $Z_{k,i}$
and not of $X_i$. One would then need to relate the derivative of the
distribution function of $Z_{k,i}$ to that of $X_i$. In the problem at
hand, if $Z_{k,i}$ has a density for $k$ large enough, then by
increasing $k$, one can assume that this density is smooth. However,
if no matter how large $k$ is $Z_{k,i}$ does not have a density, then
$Y$ may still have one. As mentioned in the introduction, Solomyak's
(1995) work is a cruel reminder on how difficult it is to study the
marginal distribution of such a basic time series model.

For more general linear processes, one could imagine to use similar blocking
techniques. One would then need to extend Theorem {\MainTheorem} in replacing
the $\M_{c_i}F$'s by some more general $F_i$'s, each $F_i$ being in 
$SR_{-\alpha,\omega}$. Formally, the expression is obvious to obtain. A proof
along the line of Theorem {\MainTheorem} is possible under some assumptions
on the $F_i$'s. However, it is not clear that such a refinement presents much
interest in applications. There is in fact a great similarity with the central
limit theorem for density, where one does not need
the summands to have a density, but the characteristic function to be in some
$L^r$ space with some $r$ at least $1$ (see Feller, 1971, \S XV.5). Yet, 
in most real life applications --- if not all --- the local central limit
theorem is applied when the summands have a well behaved density. The same
seems true about the conditions of Theorem {\MainTheorem}.

\bigskip


\section{Implementing the expansion.}

In this section we illustrate Theorem {\MainTheorem} and its implementation
through examples. In the first subsection we discuss the
number of terms in the expansion given by Theorem \MainTheorem. In the
second subsection, we explain how for standard statistical distributions,
matrix identification of the differential and multiplication operators
provides a computationally efficient way to derive tail expansions.
In the third subsection, we discuss an open problem related to
second order formula for weighted convolutions under the assumption of
regular variation with remainder.

In what follows, if the sequence $c$ is nonnegative, we write $C_p$ for 
the series $\sum_{i\in\ZZ}c_i^p$. When $C_{p+q}$ as well as $C_p$ 
and $C_q$ are finite, we write $C_{p;q}$ for $C_{p+q}-C_pC_q$.

\bigskip

\subsection{How many terms are in the expansion?}%
Let us consider the expansion given in Theorem {\MainTheorem} for the tail
of the distribution function, that is when $k$ vanishes. It asserts
$$
  \oG_c=\sum_{i\in\ZZ} \L_{G_c\natural \M_{c_i}F,m}\overline{\M_{c_i}F}
  +  o\bigl(\Id^{-m}\oF\bigr) \, .
  \eqno{\equa{GAsymptm}}
$$
It seems that this formula provides an $m+1$-term expansion. However,
one should remember that the number of terms in an expansion depends
crucially on the asymptotic scale chosen. In the case at hand, we will
show that, for natural scales, this formula may give more or less than
$m$ terms. We suspect that this fact explains the failure of previous
purely analytical attempts to find only a two terms expansion in the general
infinite order case. The algebraic flavor of the Laplace characters will 
be developed further in the next subsection, and will be the key to efficiently
derive asymptotic expansions on specific distributions.

When $m$ is $1$, formula {\GAsymptm} is
$$
  \oG_c=\sum_{i\in\ZZ} (\overline{\M_{c_i}F} 
  -\mu_{G\natural \M_{c_i}F,1}\overline{\M_{c_i}F}{}')
  + o\bigl(\Id^{-1}\oF\bigr) \, .
  \eqno{\equa{GAsymptOne}}
$$
For simplicity of the discussion, assume that the $c_i$'s are nonnegative,
and let us ask: how many terms does {\GAsymptOne} provide?

A naive answer is $2$, for the first term is 
$\sum_{i\in\ZZ}\overline{\M_{c_i}F}$, while the second one is
$-\sum_{i\in\ZZ}\mu_{G\natural \M_{c_i}F,1}\overline{\M_{c_i}F}{}'$.
The following four examples show that the truth is more complex.

\medskip

\noindent{\it Example 1.} Take $\oF(t)=t^{-\alpha} (1+1/\log t)$ for $t$
large enough. Then
$$\eqalign{
  \overline{\M_cF}(t)
  &{}={c^\alpha\over t^\alpha} \Bigl( 1+{1\over \log t-\log c}\Bigr) \cr
  &{}={c^\alpha\over t^\alpha} + {c^\alpha\over t^\alpha\log t}
   + {c^\alpha\over t^\alpha\log t} {\log c\over \log t-\log c} \, . \cr
  }
$$
Moreover, $\overline{\M_cF}{}'(t)\sim -\alpha c^\alpha/t^{\alpha+1}$ 
as $t$ tends to infinity. Applying formula {\GAsymptOne}, we obtain
$$
  \oG_c(t) = {C_\alpha\over t^\alpha} + {C_\alpha\over t^\alpha\log t}
  \bigl(1+o(1)\bigr) \, .
$$
The interesting point is that this two terms expansion is given by the
sum $\sum_{i\in\ZZ}\overline{\M_{c_i}F}$; it does not involve the term
$\sum_{i\in\ZZ}\mu_{G\natural \M_{c_i}F,1}\overline{M_{c_i}F}{}'$.

\medskip

\noindent{\it Example 2.} Take $\oF(t)=t^{-\alpha}+t^{-\alpha-2}$ for $t$
large enough, and assume that $F$ is supported inside the nonnegative 
half line. This last assumption ensures that $\mu_{F,1}$ does not vanish. Now,
$$
  \oF'(t) = -\alpha t^{-\alpha-1} - (\alpha+2)t^{-\alpha-3} \, .
$$
Consequently, {\GAsymptOne} yields
$$\eqalign{
  \oG_c(t)
  &{}= C_\alpha t^{-\alpha} - \alpha t^{-\alpha-1} \sum_{i\in\ZZ} (C_1-c_i)
   c_i^\alpha\mu_{F,1} + o(t^{-\alpha-1}) \cr
  &{}= C_\alpha t^{-\alpha} +\alpha C_{\alpha; 1} \mu_{F,1}
     t^{-\alpha-1}
     + o(t^{-\alpha-1})\, .\cr
  }
$$
Therefore, if $C_1C_\alpha\not= C_{\alpha+1}$, the first term is given by
a contribution from $\sum_{i\in\ZZ} \overline{\M_{c_i}F}$, while the second one
comes from $\sum_{i\in\ZZ} \overline{\M_{c_i}F}{}'$.

\medskip

\noindent{\it Example 3.} Take $\oF(t)=t^{-\alpha}-t^{-\alpha-3}$ for large
$t$, and assume now that $F$ is symmetric. Hence $\mu_{F,1}$ vanishes. Formula
{\GAsymptOne} yields
$$
  \oG_c(t) = C_\alpha t^{-\alpha} + o(t^{-\alpha-1}) \, .
$$
In some sense, the formula fails to give a two terms expansion. So we 
consider formula {\GAsymptm} when $m$ is $2$. Since $\mu_{F,1}$ vanishes,
it yields
$$
  \oG_c(t) 
  = \sum_{i\in\ZZ} \overline{\M_{c_i}F}(t) + {1\over 2} \sum_{i\in\ZZ}
  \mu_{G\natural \M_{c_i}F,2} \overline{\M_{c_i}F}^{(2)}(t)
  + o(t^{-\alpha-2}) \, .
  \eqno{\equa{GAsymptTwo}}
$$
Also, since $\mu_{F,1}$ is zero,
$$
  \mu_{G\natural \M_{c_i}F,2} 
  = E\Bigl( \sum_{\matrix{\noalign{\vskip -2pt}\ss j\in\ZZ\cr
                          \noalign{\vskip -5pt}\ss j\not= i\cr}} 
            c_jX_j \Bigr)^2
  = (C_2-c_i^2)\mu_{F,2} \, .
$$
Consequently,
$$\eqalign{
  \oG_c(t)
  &{}= C_\alpha t^{-\alpha} + {1\over 2} \sum_{i\in\ZZ} (C_2-c_i^2)\mu_{F,2}
    \alpha (\alpha+1) c_i^\alpha t^{-\alpha-2} + o(t^{-\alpha-2}) \cr
  &{}= C_\alpha t^{-\alpha} - {\alpha (\alpha+1)\over 2} 
    C_{\alpha;2}\mu_{F,2} t^{-\alpha-2} + o(t^{-\alpha-2}) \, .
   \cr
  }
$$
We see that the second term in this formula comes in fact from the 
third term in \GAsymptm. Notice that the second term in the expansion
for $\oF$ played no role.

\medskip

\noindent{\it Example 4.} Let us modify example 2 in writing 
$\oF(t)= t^{-\alpha}+at^{-\alpha-2}$ for large $t$. And let us take
$F$ to be symmetric. Formula {\GAsymptm} with $m$ equal $2$ yields
\hfuzz=4pt
$$\eqalign{
  \oG_c(t)
  &{}= C_\alpha t^{-\alpha} + a C_{\alpha+2}t^{-\alpha-2}
    + {\alpha(\alpha+1)\over 2}\mu_{F,2} \sum_{i\in\ZZ} (C_2-c_i^2) c_i^\alpha
    t^{-\alpha-2} \cr
  \noalign{\vskip -10pt}
  &\hbox{\hskip 3in}{}+ o(t^{-\alpha-2}) \cr
  &{}= C_\alpha t^{-\alpha} + t^{-\alpha-2} 
   \bigl(a C_{\alpha+2} - 0.5\alpha(\alpha+1)\mu_{F,2} C_{\alpha;2}\bigr)
   + o(t^{-\alpha-2}) \, . \cr
  }
$$
\hfuzz=0pt
For any fixed sequence $c=(c_i)_{i\in\ZZ}$ of nonnegative numbers, we can 
find $a$ such that
$$
  a C_{\alpha+2} -0.5\alpha(\alpha+1)\mu_{F,2}C_{\alpha;2}=0 \, .
$$
In this case, we need to include at least one more term, taking at least
$m=3$ in formula {\GAsymptm} if we want two nonzero terms.

\medskip

A first point of these examples is that obtaining simply a two terms
expansion formula is rather hopeless in general. Clearly, some tricky
cancelation may occur, depending on the particular sequence $c$
considered as well as the coefficients in the expansion  of $\oF$.

A second point is that when some cancellations occur, we may need to 
add terms; at least if we can, since the condition $m$ less than $\alpha$
caps the number of terms we can obtain. In our examples, this was easy,
but in other cases, this may involve far more complicated calculations.
Hence, we need more effective ways to implement the expansion, and this
will be discussed in the next subsection.

A third point is that it raises an interesting question: Is it possible
that for any fixed $m$, as large as we want, can we find a distribution
for which formula {\GAsymptm} provides only a $1$ term expansion? We address
this question in the rest of this subsection. The answer is yes, though it
is not a generic case. To make this clear, let us examine the following
statement.

\bigskip

\noindent{\bf Statement.} {\quad\it%
  Generically, {\GAsymptm} with $m$ equal $2$ provides at least a two
  terms expansion, even if $\mu_{F,1}$ vanishes.
}

\bigskip

When $m$ is $2$, the right hand side in {\GAsymptm}, up to the 
$o(\cdot)$-term, is
$$
  \sum_{i\in\ZZ} \overline{\M_{c_i}F}
  -\sum_{i\in\ZZ} \mu_{G\natural \M_{c_i}F,1} \overline{\M_{c_i}F}{}'
  + {1\over 2} \sum_{i\in\ZZ} 
  \mu_{G\natural\M_{c_i}F,2}\overline{\M_{c_i}F}{}''
  \, .
$$
Let $\sigma_F^2$ be the variance of $F$. Since
$$\eqalignno{
  \mu_{G\natural\M_{c_i}F,1} &{}= (C_1-c_i)\mu_{F,1} \cr
  \mu_{G\natural\M_{c_i}F,2} &{}=(C_2-c_i^2)\sigma_F^2
                              + (C_1-c_i)^2\mu_{F,1}^2 \cr
  }
$$
and $\overline{M_{c_i}F}{}''(t) \sim c_i^\alpha \oF''(t)$, as $t$ tends to 
infinity, the formula reads
$$\displaylines{
  \sum_{i\in\ZZ} \overline{\M_{c_i}F} 
  -\sum_{i\in\ZZ} (C_1-c_i)\mu_{F,1}\overline{\M_{c_i}F}{}'
  \hfill\cr\hfill
  {}-{1\over 2} \Bigl( C_{\alpha;2}\sigma_F^2
  + (C_1C_{\alpha;1}-C_{\alpha+1;1})\mu_{F,1}^2\Bigr) \oF''(t)
  +  o\bigl(\oF(t)/t^2\bigr) \, .\cr 
  }
$$
Generically, there is no cancellation, and even if $\mu_{F,1}$ vanishes,
the term in $\oF''(t)$ yields a nonzero term which (generically) does not
cancel with $\sum_{i\in\ZZ}\overline{\M_{c_i}F}$. So, the formula yields
at least two terms.

\medskip

The companion of the positive statement is as follows.

\bigskip

\noindent{\bf Statement.}{\qquad\it%
  Consider the formula {\GAsymptm} for both $m$ and the sequence 
  $c=(c_i)_{i\in\ZZ}$ fixed. Then there exists a distribution (depending
  on $m$ and the sequence $c$) for which the formula yields only a $1$-term
  expansion.
}

\bigskip

We will prove this statement in the next section, once we have an effective
way of calculating expansions.

\bigskip

\subsection{Practical implementation: from Laplace's characters to 
linear algebra.}%
The few examples that we gave in the previous subsection were calculated
easily by hand. In more complicated cases as well as for theoretical 
reasons (for instance to prove the statement at the end of the previous
subsection), we need more effective ways to derive expansions. In Barbe
and McCormick (2005), we introduced an algebraic technique, a tail calculus,
which allows one to derive asymptotic expansions for finite weighted 
convolutions of distribution functions such that 
$\oF(t)\sim t^{-\alpha}P(1/t)$ for some 
polynomial $P$. An extension to infinite order moving averages was carried out
for this class in Barbe and McCormick (200?).
Distributions for which this calculus is applicable include the Pareto, 
Student and the Generalized Pareto. However, it does not apply to other
well used distributions in the heavy tail literature. So we will develop
a tail calculus for those. It may be possible to derive an abstract
tail calculus which covers many more examples. However, we have not been
able to find a good formalism for that purpose. Consequently we will say 
few things on the general idea, and mostly develop two examples.

\medskip

\noindent{\it Generalities on tail calculus.} Assume that we start with a
distribution function $F$ such that $\oF$ and $\overline{\M_{-1}F}$ 
have asymptotic expansions
in an asymptotic scale $e=(e_i)_{i\in I}$ say. That is, for some
coefficients $p_{\oF,i}$ and $p_{\overline{\M_{-1}F},i}$,
$$
  \oF\sim \sum_{i\in I} p_{\oF,i}e_i 
  \qquad\hbox{and}\qquad 
  \overline{\M_{-1}F} \sim\sum_{i\in I}p_{\overline{\M_{-1}F},i}e_i \, .
$$ 
We write $p_{\oF}$ and $p_{\overline{\M_{-1}F}}$ for the vectors 
$(p_{\oF,i})_{i\in I}$ and $(p_{\overline{\M_{-1}F},i})_{i\in I}$.
Note that since $e$ is an asymptotic scale and both $F$ and $\oF$ are
nonnegative, the first nonzero term of $p_{F}$ and $p_{\overline{\M_{-1}F}}$
is positive.
Writing $|I|$ the cardinality of $I$, we can think of these expansions as 
a coding of the tails in the vector space $\RR^{2|I|}$ via the map
$$
  \rho \, : \, F\mapsto \rho(F) 
  = (p_{\oF},p_{\overline{\M_{-1}F}}) \, .
$$
The expansion in Theorem {\MainTheorem} involves multiplication by constants.
Both the derivative $\D$ and the multiplication operators $\M_c$ act 
componentwise on the scale $e$ by $\D e=(\D e_i)_{i\in I}$ and
$\M_c e=(\M_c e_i)_{i\in I}$. For every $t$, we can take linear
combinations of the components of the vector $e(t)$ by left multiplying
it by a matrix. For vector valued functions, we extend the expansion
notation $\sim$ as acting componentwise. So, a vector valued function
$f=(f_1,\ldots , f_p)$ has asymptotic expansion $g=(g_1,\ldots , g_p)$
if $f_i\sim g_i$ for each $i=1,\ldots , p$.

\bigskip

\Definition
{ An asymptotic scale $e=(e_i)_{i\in I}$ is a $\star$-asymptotic scale if
  there exist matrices $\calD$ and $\calM_c$ such that $\D e\sim\calD e$
  and for any $c$ positive, $\M_ce\sim\calM_c e$.
}

\bigskip

In particular, in an asymptotic sense, a $\star$-asymptotic scale is closed
under differentiation and multiplication of the variable by a positive 
constant.

For a $\star$-asymptotic scale, the following diagram commutes.

{

\def\mapright#1{\hskip-3pt\smash{\mathop{\hbox to 30pt{\rightarrowfill}}%
\limits^{#1}}\hskip-3pt}
\def\mapleft#1{\hskip-3pt\smash{\mathop{\hbox to 30pt{\leftarrowfill}}%
\limits^{#1}}\hskip-3pt}
\def\mapdown#1{\Big\downarrow\rlap{$\vcenter{\hbox{$\scriptstyle#1$}}$}}
$$
\matrix{
\M_cF&\mapleft{\M_c}&F&\mapright{\D}&F'\cr
\noalign{\vskip 2pt}
\mapdown{\rho}&  &\mapdown{\rho}&  &\mapdown{\rho}\cr
(p_{\overline{\M_cF}},p_{\overline{\M_{-c}F}})&\mapleft{\calM_c}%
&(p_{\oF},p_{\overline{\M_{-1}F}})&\mapright{\calD}&%
(p_{\oF'},p_{\overline{\M_{-1}F}{}'})\cr
}
$$
}

In particular, the identification of tails with the finite vector space
$\RR^{2|I|}$ implies that the operators $\calD$ and $\calM$ involved in the
lower part of the diagram are simply matrices. This allows one to identify the
Laplace character $L_{K,m}$ with the matrix
$$
  \calL_K=\sum_{0\leq j\leq m} {(-1)^j\over j!}\mu_{K,j}\calD^j \, .
$$
Note that this matrix depends on the order $m$ of the Laplace character
as well as on the $\star$-asymptotic scale chosen. The notation does not
show this dependence. Theorem {\MainTheorem} implies that the following 
diagram commutes.

{

\def\mapright#1{\hskip-3pt\smash{\mathop{\hbox to 30pt{\rightarrowfill}}%
\limits^{#1}}\hskip-3pt}
\def\mapleft#1{\hskip-3pt\smash{\mathop{\hbox to 30pt{\leftarrowfill}}%
\limits^{#1}}\hskip-3pt}
\def\mapdown#1{\Big\downarrow\rlap{$\vcenter{\hbox{$\scriptstyle#1$}}$}}
$$
\matrix{
(F,G)&\mapright{\star}&F\star G\cr
\noalign{\vskip 2pt}
\mapdown{}&  &\mapdown{}&\cr
(p_{\oF},p_{\oG},\calL_F,\calL_G)%
&\mapright{}%
&(\calL_Fp_{\oG} + \calL_G p_{\oF},\calL_F\calL_G)\cr
}
$$
}
In particular, the map $F\mapsto\calL_F$ is a linear representation of the
convolution algebra, whose dimension is $|I|$.

Practically, the diagram says that on a $\star$-asymptotic scale, asymptotic
expansions for weighted convolution can be done by manipulating finite
dimensional vectors and matrices. This is the key to effective computation.

In particular, using the existence of inverse for the Laplace characters, 
the asymptotic expansion of Theorem {\MainTheorem} is given by
$$\eqalign{
  p_{\oG_c^{(k)}}
  &{} = \sum_{i\in\ZZ} \calL_{G_c}\calL_{\M_{c_i}F}^{-1} 
    \calD^k\calM_{c_i}\bigl( \II\{\, c_i>0\,\} p_{\oF} 
    + \II\{\, c_i<0\,\} p_{\overline{\M_{-1}F}} \bigr)\cr
  &{}= \calL_{G_c} \sum_{i\in\ZZ} \calL_{\M_{c_i}F}^{-1} 
    \calD^k\calM_{c_i}\bigl( \II\{\, c_i>0\,\} p_{\oF} 
    + \II\{\, c_i<0\,\} p_{\overline{\M_{-1}F}} \bigr)\, .
   \cr
  }
$$
This expresses the asymptotic expansion as sums and products of
matrices.  As we did in Barbe and McCormick (2005), this method of
calculation is suitable for computer implementation.  It allows one to
automatize asymptotic expansion for this type of weighted heavy tail
convolutions.

\medskip

In the next two examples, we will work out the details and show how to 
obtain expansions with as many terms as the integrability 
condition $m<\alpha$ and computer memory allow.

\bigskip

\noindent{\it Distributions with asymptotic expansion of the form 
$\sum_{i\geq 0}a_it^{-\alpha_i}$.}\
The Hall-Weissman (1997) distributions are defined by 
$$
  \oF(t)=at^{-\alpha}+ bt^{-\beta} \, ,
$$
with $\alpha<\beta$ and for $t$ at least some $t_{a,b}$. Obviously, 
they have an expansion of the form $\sum_{i\geq 0}a_it^{\alpha_i}$.

The Burr distributions with positive paramters $\tau$ and $\gamma$ are 
defined by
$$
  \oF(t) = (1+t^\tau/\beta)^{-\gamma}\, , \rlap{$\qquad t\geq 0\, .$}
$$
When $\gamma$ is $1$, it is also called the log-logistic distribution.
A tail expansion is easily derived. Indeed, $\oF(t)=\beta^\gamma
t^{-\gamma\tau}(1+\beta t^{-\tau})^{-\gamma}$, and the expansion
$$
  \oF(t) \sim \beta^\gamma \sum_{k\geq 0}
  (-1)^k \beta^k {\Gamma (\gamma+k)\over k!\Gamma (\gamma)} 
  t^{-\tau\gamma-k\tau}
  \eqno{\equa{BurrExpansion}}
$$
follows. This expansion is of the form $\sum_{i\geq 0} a_i t^{-\alpha_i}$.

The Fr\'echet distribution is defined by $\oF(t)=1-\exp (-t^{-\alpha})$. Its
asymptotic expansion is given by
$$
  \oF(t)\sim \sum_{k\geq 1} {(-1)^{k+1}\over k!} t^{-\alpha k} \, .
$$
Again, this expansion has the form $\sum_{i\geq 0} a_i t^{-\alpha_i}$.

So, it is of some interest here to consider an increasing sequence 
$0<\alpha_0<\ldots < \alpha_p$ of real numbers, and investigate expansions
for weighted convolutions of distribution funcitons having an asymptotic 
expansion in the
scale $t^{-\alpha_i}$, $0\leq i\leq p$. For a more comprehensive list of 
heavy tailed distributions with examples of their use in modeling, 
particularly for insurance risk data, we refer to Beirlant et al.~(1996),
Embrechts et al.~(1997) and Rolski et al.~(1999). 

Before proceeding, we need to note several things.
There is of course no loss of generality in assuming that the coefficient of
$t^{-\alpha_0}$ in the tail expansion of $\oF$ does not vanish. Hence, we
can assume that $\oF(t)\asymp t^{-\alpha_0}$. We will apply 
Theorem {\MainTheorem} with $k=0$. This requires
$m<\alpha_0$; and so, up to truncating the sequence, we can also assume that
$\alpha_p\leq \alpha_0+m<2\alpha_0$. The next thing to notice is that the
asymptotic scale $(t^{-\alpha_i})_{0\leq i\leq p}$ may not be a 
$\star$-asymptotic scale. Indeed, it is closed under 
differentiation if and only if for any integer $0\leq i\leq p$ and any 
nonnegative integer $k$ for which $\alpha_i+k\leq \alpha_p$, there exists
a $j$ such that $\alpha_i+k=\alpha_j$. If this is not the case, we enlarge
our collection of $\alpha_j$'s by adding the missing ones recursively.
Thus, we assume without loss of generality that 
$(t^{-\alpha_i})_{0\leq i\leq p}$ is a $\star$-asymptotic scale.

Now, let us consider the functions $e_i(t)=t^{-\alpha_i}$, 
$0\leq i\leq p$. They
form the basis of a finite dimensional vector space, isomorphic
to $\RR^{p+1}$. Since $\D e_i=-\alpha_i e_j$ with $j$ defined by 
$\alpha_j=\alpha_i+1$, the derivative is identified to the 
matrix $\calD$  defined by
$$
  \calD e_i
  =\cases{ -\alpha_i e_j & if $\alpha_j=\alpha_i+1\leq \alpha_p$,\cr
           0             & otherwise.\cr}
$$
Since $\calD$ maps $e_i$ into the space spanned by $e_{i+1},\ldots , e_p$,
it is a nilpotent matrix.

Next, since $\M_c e_i=c^{\alpha_i}e_i$, the matrix
$\calM_c$ is the diagonal matrix $\diag(c^{\alpha_i})_{0\leq i\leq p}$.

\medskip

To show how this works practically, we consider the Burr distribution 
with $\gamma=10$ and $\tau=3/2$ say (other values would work just as well).
Then $\alpha_0=\gamma\tau=15$. Then, Theorem {\MainTheorem} allows for
a $14$ terms expansion of the weighted convolution. However, we limit
the expansion by chosing $m=4$ in Theorem {\MainTheorem},
though it will be clear at the end that to do so is motivated not by
any difficulty in obtaining the terms in such high order expansion, but
rather by concern with the space and display of such expansion.
The {\tt Maple} code that we used to implement the formal calculations 
is given in the appendix.

So we choose $m$ to be $4$. Display {\BurrExpansion} shows that $\oF$ has an
asymptotic expansion in the scale $t^{-15-3i/2}$, $i\geq 0$. Since we want
$4$ terms and $t^{-4}\oF(t)\asymp t^{-19}$, we can restrict the range
of $i$'s to $15+3i/2\leq 19$, that is $i\leq 2$. The expansion of $\oF$ is
then
$$
  \oF(t)= \beta^{10} t^{-15} - 10\beta^{11}t^{-33/2} + 55\beta^{12}t^{-18}
  + o(t^{-19})\, .
  \eqno{\equa{BurrExplicit}}
$$
So we need to consider the functions $t^{-15}$, $t^{-15-3/2}$ and $t^{-18}$.
This family is not closed under differentiation modulo terms in $o(t^{-19})$.
To close it, we need to add the derivative of $t^{-15}$, that is up to a
multiplicative constant, $t^{-16}$. But then we need to add the derivative
of $t^{-16}$ as well, that is to add $t^{-17}$. And similarly, we must add
$t^{-19}$. Next we also need to add the derivative of $t^{-15-3/2}$, which is
up to scale $t^{-18-1/2}$. This process leads to the $\star$-asymptotic scale
$(e_i)_{0\leq i\leq 7}$ corresponding to the following $\alpha_i$'s:
$$
  \alpha_0=15 <16 < 16+1/2 < 17 < 17+1/2 < 18 < 18+1/2 < 19 \, .
$$
Since we have eight $\alpha_i$'s, we need to work in the space $\RR^8$.

Display {\BurrExplicit} shows that
$$
  \oF = \beta^{10} e_0 - 10\beta^{11}e_2 + 55\beta^{12} e_5 + o(e_7)\, .
$$
So, the vector $p_{\oF}$ in $\RR^8$ is simply 
$$
  p_{\oF}=(\beta^{10},0,-10\beta^{11},0,0,55\beta^{12},0,0) \, .
$$

Since $\D e_0=-15 e_1$ and $\D e_i=-\alpha_ie_{i+2}$ for $i$ at least $1$, 
the matrix that corresponds to the differentiation is defined by

\setbox1=\hbox{\largecmr 0}\wd1=0pt\ht1=0pt\dp1=0pt
$$\calD=\pmatrix{
0   &     &       &     &       &       &   &   \cr
-15 & 0   &       &     &       &       &   &   \cr
0   & 0   & 0     &     &\copy1 &       &   &   \cr
    & -16 & 0     & 0   &       &       &   &   \cr
    &     & -16.5 & 0   & 0     &       &   &   \cr
    &     &       & -17 & 0     &  0    &   &   \cr
    &\box1&       &     & -17.5 &  0    & 0 &   \cr
    &     &       &     &       & -18   & 0 & 0 \cr
}
$$
Moreover, since $\M_ce_k=c^{\alpha_k}e_k$, the multiplication operator is
represented by
$$
  \calM_c = \diag( c^{\alpha_k})_{0\leq k\leq 7} \, .
$$

A Laplace character $\L_{K,4}$ is identified to the matrix
$$
  \calL_K 
  = \sum_{0\leq j\leq 4} {(-1)^j\over j!} \mu_{K,j}\calD^j \, .
$$

\noindent This shows that $\calL_{\M_cF}$ is the matrix
\hfuzz=7pt
\setbox1=\hbox{\largecmr 0}\ht1=0pt\dp1=0pt\wd1=0pt
$$
\pmatrix{
1 & & & & & & \cr
\noalign{\vskip 3pt}
15c\mu_1 & 1 & & & & & & \cr
\noalign{\vskip 3pt}
0 & 0 & 1 & &\box1 & & & \cr
\noalign{\vskip 3pt}
120c^2\mu_2 & 16c\mu_1 & 0 & 1 & & & & \cr
\noalign{\vskip 3pt}
0 & 0 & {\ds 33\over\ds 2}c\mu_1 & 0 & 1 & & & \cr
\noalign{\vskip 3pt}
680c^3\mu_3 & 136c^2\mu_2 & 0 & 17c\mu_1 & 0 & 1 & & \cr
\noalign{\vskip 3pt}
0 & 0 & {\ds 1155\over\ds 8}c^2\mu_2 & 0 & {\ds 35\over\ds 2}c\mu_1 & 0 & 1 &\cr
\noalign{\vskip 3pt}
3060c^4\mu_4 & 816c^3\mu_3 & 0 & 153c^2\mu_2 & 0 & 18c\mu_1 & 0 & 1 \cr
}
$$
\hfuzz=0pt

Finally, $\oG_c$ has an expansion whose coefficients are given by
$$
  \calL_{G_c}\sum_{i\in\ZZ} \calL_{\M_{c_i}F}^{-1}\calM_{c_i} p_{\oF} 
  \, .
$$
We formally compute the vector $\calL_{\M_{c_i}F}^{-1}\calM_{c_i}p_{\oF}$.
Assuming that the $c_i$'s are positive, summing the vectors 
$\calL_{\M_{c_i}F}^{-1}\calM_{c_i}p_{\oF}$ leads to expressions involving
$C_s=\sum_{i\in\ZZ} c_i^s$. The matrix $\calL_G$ involves moments of $G_c$
which in turn are expressed in terms of moments of $F$ and other expressions
involving the terms $C_s$.

We can then compute the matrix for the Laplace character associated to $G$.
We then do a formal matrix multiplication. To write the
final result recall the notation $C_{p;q}$ for $C_{p+q}-C_pC_q$.
We also write $\kappa_3$ for the third central moment of $F$, that is
for $E(X-\mu_{F,1})^3$; and $\kappa_4$ for the fourth one, 
$E(X-\mu_{F,1})^4$.
We then proved that $P\{\, \<c,X\> > t\,\}$ has asymptotic expansion
$\sum_{0\leq i\leq 7} q_ie_i$ with
$$\eqalignno{
q_0={}&\beta^{10}C_{15} \, , \cr 
q_1={}&-15\beta^{10}\mu_1 C_{15;1} \, , \cr
q_2={}&-10\beta^{11} C_{33/2} \, ,\cr
q_3={}&120\beta^{10}\mu_1^2 (C_{17;1}-C_1C_{15;1})
       - 120\beta^{10}\sigma^2 C_{15;2} \, , \cr
q_4={}&165\beta^{11}\mu_1 C_{33/2;1} \, , \cr
q_5={}&{}-680\beta^{10}\mu_1^3(C_{17:1}-2C_1C_{16;1}+C_1^2C_{15;1})
     \cr
      &{} + 2040\beta^{10}\sigma^2\mu_1 (C_{17;1}-C_2C_{15;1}
       - 680\beta^{10}\kappa_3 C_{15;3}
       + 55\beta^{12}C_{18} \, , \cr
q_6={}&{5775\over 4}\beta^{11}\bigl( -\mu_1^2 (C_{35/2;1}-C_1C_{33/2;1})
      +\sigma^2C_{33/2;2}\bigr) \cr
q_7={}&3060\beta^{10}\mu_1^4 (C_{18;1}-3C_1C_{17;1}-3C_1^2C_{16;1}
       - C_1^3C_{15;1}) \cr
      &{}-18360\beta^{10}\sigma^2\mu_1^2 (C_{18;1}-C_1C_{17;1}-C_2C_{16;1}
       -C_1C_2C_{15;1}) \cr
      &{}-9180\beta^{10}\sigma^4 (2C_{17;2}-C_{15}C_{2;2}) \cr
      &{}+12240\beta^{10}\kappa_3\mu_1 (C_{18;1}-C_3C_{15;1}) \cr
      &{}-990\beta^{12}\mu_1 C_{18;1}
       -3060\beta^{10}\kappa_4 C_{15;4} \, . \cr
}
$$
It should be clear now that Theorem {\MainTheorem} is not a pure 
abstraction and that the algebraic nature of the Laplace character 
is what makes such a computation possible.

\medskip

\noindent{\it The log-gamma distribution.} Recall that $X$ has a log-gamma
distribution with parameter $(\lambda,\alpha)$ if its density is given by
$$
  f(x)= {\alpha^\lambda\over \Gamma (\lambda)}
  (\log x)^{\lambda-1}x^{-\alpha-1}\, , \quad
  x\geq 1 \, .
$$
This is equivalent to saying that $X=\exp(Z/\alpha)$ where $Z$ has a standard
Gamma distribution with parameter $\lambda$. Successive integrations by parts
show that
$$
  P\{\, Z>t \,\} 
  \sim  {e^{-t}\over \Gamma (\lambda)} \sum_{k\geq 0} (\lambda-1)_k 
  t^{\lambda-1-k} \, . 
$$
The distribution function $F$ of $X$ is given 
by $P\{\, Z\leq \alpha \log x\,\}$. Thus, it has expansion
$$
  \oF(t) \sim {1\over \Gamma(\lambda) t^\alpha} \sum_{k\geq 0} (\lambda-1)_k
  \alpha^{\lambda-1-k} (\log t)^{\lambda-1-k} \, .
$$
This is an expansion in the scale $e_k(t)=t^{-\alpha}(\log t)^{\lambda-1-k}$,
$k\geq 0$. 

Since 
$$
  e_k'(t)=t^{-\alpha-1}(\log t)^{\lambda-1-k}
  \Bigl( -\alpha + {\lambda-1-k\over \log t}\Bigr) \, ,
$$
we see that $e_k'=o(e_i)$ at infinity for any $i$. Consequently, any finite
family $(e_i)_{0\leq i\leq p}$ is closed under differentiation up to
$o(e_{p+k})$, $k\geq 1$. Moreover $\D e_i=o(e_p)$ implies that the matrix
corresponding to derivative, $\calD$, is $0$.

To determine the matrices $\calM_c$, we write
$$\eqalign{
  e_j(t/c)
  &{}= c^\alpha t^{-\alpha} (\log t)^{\lambda-1-j} 
    \Bigl( 1-{\log c\over \log t}\Bigr)^{\lambda-1-j} \cr
  &\sim \sum_{k\geq 0} {(\lambda-1-j)_k\over k!} (-1)^k c^\alpha (\log c)^k
    e_{j+k}(t) \, . \cr
  }
$$
Hence, the matrix $\calM_c$ is the lower triangular matrix defined by
$$
  (\calM_c)_{i,j} = 
  \cases{ {\ds (\lambda-1-j)_{i-j}\over\ds (i-j)!}(-1)^{i-j} c^\alpha 
          (\log c)^{i-j}
          & if $i\geq j\geq 0$, \cr
          0
          & otherwise. \cr}
$$
(Note that we indexed our matrix so that the first row and column are labeled
$0$.) In particular $(e_k)_{k\geq 0}$ is a $\star$-asymptotic scale.
Since $\calD$ vanishes, the Laplace character of $F$ is the identity. So the 
formula in Theorem {\MainTheorem} yields for any $\alpha$ more than $1$,
$$
  G_c(t)\sim \sum_{i\in\ZZ} \calM_{c_i}p_{\oF}
$$
with 
$$
  p_{\oF}={\alpha^{\lambda-1}\over \Gamma(\lambda)} 
  \Bigl(1,{\lambda-1\over\alpha},{(\lambda-1)_2\over \alpha^2},\ldots , 
  {(\lambda-1)_p \over\alpha^p} \Bigr) \, .
$$
Writing $(C\log C)_{r,s}$ for $\sum_{i\in\ZZ} c_i^r(\log c_i)^s$, we obtain
for instance the asymptotic expansion
$$
  \oG_c\sim \sum_{j\geq 0} q_je_j
$$
with
$$
  q_j={\alpha^{\lambda-1}\over\Gamma (\lambda)}\sum_{0\leq i\leq j} 
  {(\lambda-1)_{j-i}\over\alpha^{j-i}} {(\lambda-1-j+i)_i\over i!} (-1)^i 
  (C\log C)_{\alpha,i} \, .
$$

The interesting feature of this example is that for any $m$, the
identification $\calL_F$ of the $m$-th Laplace character of $F$ is the
identity. A consequence is that in this asymptotic scale, increasing
$m$ in Theorem {\MainTheorem} does not yield extra terms compared to
taking $m$ to be $0$. This is entirely due to the choice of the
asymptotic scale. However, the expansion given in Theorem
{\MainTheorem} provides in fact a better estimate. The reason is that
the asymptotic scale in Theorem {\MainTheorem} is adapted to the
distribution function $F$ and its derivatives. It is not a scale given
a priori on which $\oF$ is expanded.

This situation is very similar to the fact that the normal tail
$$
  \overline\Phi (t) = \int_t^\infty {e^{-x^2/2}\over \sqrt{2\pi}} \d x
$$
has a one term expansion in the scale $\overline\Phi^k$, $k\geq 1$, say,
but it has an expansion given by a divergent series in the asymptotic scale
$t^{-k}e^{-t^2/2}$, $k\geq 0$.

\medskip

\noindent{\it A degenerate case.} For the log-gamma distribution, we
saw that we can find an asymptotic scale such that Theorem {\MainTheorem} 
with $m$ equal $0$ provides as many terms as we like. In the example to
be developed now, we show that the reverse situtation may occur; that is,
due to rather exceptional cancellations, there are distributions 
such that Theorem {\MainTheorem} provides
only $1$ term. This was the statement ending 
subsection 3.1. For simplicity we will show this only when the $c_i$'s
are nonnegative. It is conceptually easy to adapt the proof if this does 
not hold.

So, let us first choose $m$, as large as we want. Let $\alpha$ be larger than
$m$, and let $(c_i)_{i\in\ZZ}$ be a sequence of nonnegative numbers such 
that $N_{\alpha,\gamma,\omega}(c)$ is finite. We consider the 
$\star$-asymptotic scale $e_i(t)=t^{-\alpha-i}$, $0\leq i\leq m$.
As we have seen in the first example of this subsection, the derivative
$\D$ is identified with the matrix $\calD$ whose entries are
$$
  \calD_{i,j}=\cases{ -\alpha-j & if $i=j+1$,\cr 0 & otherwise.\cr}
$$
Moreover, for any $c$ positive,
$$
  \calM_c =\diag(c^\alpha,c^{\alpha+1},\ldots , c^{\alpha+m}) \, .
$$
As before, the $m$-th Laplace character of $F$ is identified with the matrix
$$
  \calL_{F,m}=\sum_{0\leq j\leq m}{(-1)^j\over j!}\mu_{F,j}\calD^j \, .
$$
Let $p=(p_0,\ldots , p_m)$ be the coefficients of the expansion of $\oF$
in the asymptotic scale $e_i$, $0\leq i\leq m$. Then, the coefficients of 
$\oG_c$ in this asymptotic scale are given by
$$
  q=\calL_{G_c,m} \sum_{i\in\ZZ} \calL_{M_{c_i}F,m}^{-1}\calM_{c_i}p \, .
$$
So, it suffices to prove that we can find $F$ such that $q=e_0$. Note that
there is no loss of generality in assuming all the first $m$ moments of $F$ 
fixed, because fixing these moments does not put any restriction on the 
vector $p$. Hence the first $m$ moments of $G_c$ are fixed.
In other words both the Laplace characters of $F$ and $G_c$ can be taken
fixed. We already know that $\calL_{G_c,m}$ is invertible for it is a lower
triangular matrix with all diagonal terms equal to $1$. Along the lines of
Proposition {\InverseLaplaceCharB} we obtain an expression for the inverse
of $\calL_{\M_{c_i}F}$. Specifically, we define the nilpotent matrix  
$\calN=\Id-\calL_{M_{c_i}F,m}$ and have
$$
  \calL_{M_{c_i}F,m}^{-1}=(\Id -\calN)^{-1} 
   = \sum_{0\leq j\leq m} \calN^j
$$
is also lower triangular, with its diagonal elements all equal to $1$.
Consequently, the diagonal elements of $\calL_{\M_{c_i}F,m}^{-1}\calM_{c_i}$
are those of $\calM_{c_i}$. Hence the matrix 
$\sum_{i\in\ZZ}\calL_{\M_{c_i}F,m}^{-1}\calM_{c_i}$ is lower 
triangular, with diagonal elements 
$(C_\alpha,C_{\alpha+1}, \ldots , C_{\alpha+m})$. It is invertible. Therefore 
it is possible to find $p$ with positive $e_0$ component such that 
$q$ is $e_0$. This shows the
existence of a distribution function $F$ such that the formula in 
Theorem {\MainTheorem} yields only a $1$ term expansion.

This discussion is particularly relevant to the problem of ascertaining second
order regular variation for infinite order moving averages. If an $m$-term
expansion shows that, for example, $\oG_c\sim \Id^{-\alpha}+\Id^{-\alpha-m}$,
then $\oG_c$ is second-order regularly varying; but this would not be revealed
until an $m$-terms expansion was calculated. Thus, we see that second order
regular variation for linear processes is not a second order expansion 
question; rather, it is a higher order expansion question. We discuss this 
problem in detail in the next section.

\bigskip

\subsection{Two terms expansion and second order regular variation.}%
Motivated by probabilist and statistical applications, consider the
following problem: If $\oF$ is regularly varying of index $-\alpha$
with remainder, what is an asymptotic equivalent of $\oG_c-C_\alpha\oF$?

In this section, we will first explain this problem and the terminology used,
and then show that Theorem {\MainTheorem} sheds an interesting light on 
the matter.

Recall that $\oF$ is regularly varying with index $-\alpha$ and remainder,
if there exist functions $k(\cdot)$ and $g(\cdot)$, with $g$ tending to $0$ at
infinity, such that
$$
  {\oF(\lambda t)\over \oF(t)} - \lambda^{-\alpha} 
  \sim \lambda^{-\alpha} k(\lambda) g(t)
  \eqno{\equa{RVRemainder}}
$$
as $t$ tends to infinity --- see, Bingham, Goldie and Teugels (1989, \S 3.12).
Note that this relation does not change if we multiply $g$ and divide $k$ by
the same constant. If this relation holds, then $g$ must be regularly 
varying with nonpositive index $\rho$ and, up to a possible multiplication
of $g$ by a constant, necessarily
$$
  k(\lambda ) = \cases{ (\lambda^\rho-1)/\rho & if $\rho<0$,\cr
                        \log\lambda               & if $\rho=0$.\cr}
$$
We then write $\oF\in 2RV(-\alpha,g)$. Unless otherwise specified, we assume
for simplicity that the $c_i$'s are nonnegative.
The problem mentioned can be rephrased as: assuming {\RVRemainder}, we know 
that $\oF$ is regularly varying with index $-\alpha$, hence that
$$
  \lim_{t\to\infty} \oG_c(t)/\oF(t) = C_\alpha \, .
$$
What is the exact rate of convergence in this limit? One would indeed think
that adding one term to the regular variation as in {\RVRemainder} brings
one more term in the asymptotic expansion for $\oG_c$. The examples
in subsection 3.1 show that this belief is not correct in general.

Since the problem is motivated by applications in time series where it is
natural to suppose that $F$ is centered, we assume that $\mu_{F,1}$ 
vanishes --- the following discussion can easily be modified if this first
moment does not vanish. In order to apply Theorem {\MainTheorem} with $k=0$ and
$m=2$, we assume that $\oF$ is smoothly varying of order larger than $2$ and
that $\alpha$ is larger than $2$. We calculate
the second Laplace character $\L_{G_c\natural \M_{c_i}F,2}=\Id
+ (1/2)\mu_{F,2}(C_2-c_i^2)\D^2$. Then 
Theorem {\MainTheorem}, with $m=2$ and $k=0$, and the fact that 
$\oF''(t)\sim \alpha (\alpha+1)t^{-2}\oF(t)$ at infinity, yield
$$
  \oG_c(t) = \sum_{i\in\ZZ}\overline{\M_{c_i}F(t)}
  - {\alpha(\alpha+1)\over 2} C_{\alpha;2}\mu_{F,2} t^{-2}\oF(t)
  + o\bigl(t^{-2}\oF(t)\bigr) \, 
$$
as $t$ tends to infinity. Consequently,
$$
  {\oG_c\over\oF}(t) -C_\alpha
  =\sum_{i\in\ZZ} \Bigl( {\oF(t/c_i)\over \oF(t)}-c_i^\alpha\Bigr)
  - {\alpha (\alpha+1)\over 2} C_{\alpha;2}\mu_{F,2} t^{-2}
  + o(t^{-2}) \, .
$$
The global Potter type bounds of Theorem 3.1.3 in Bingham, Goldie and Teugels
(1989) and {\RVRemainder} imply
$$\displaylines{\quad
  {\oG_c\over \oF}(t) -C_\alpha
  = g(t) \sum_{i\in\ZZ} c_i^\alpha k(1/c_i)  \bigl(1+o(1)\bigr)
  \hfill\cr\hfill
  - {\alpha(\alpha+1)\over 2} C_{\alpha;2} \mu_{F,2}t^{-2}
  + o(t^{-2})\, .
  \quad\cr}
$$
Consequently, one sees that the second order term depends on the behavior of
$t^2g(t)$ at infinity. If $\lim_{t\to\infty} t^2g(t)=\infty$, then
$$
  {\oG_c\over \oF}(t) -C_\alpha
  \sim g(t)\sum_{i\in\ZZ} c_i^\alpha k(1/c_i)
$$
while if $\lim_{t\to\infty}t^2g(t)=0$ then
$$
  {\oG_c\over\oF}(t)-C_\alpha \sim -{\alpha (\alpha+1)\over 2} 
  C_{\alpha;2} \mu_{F,2}t^{-2} \, ,
$$
with the usual convention that if the constant $C_{\alpha;2}$ in the right
hand side vanishes, then the right hand side should be read as $o(t^{-2})$.
If $g(t)\sim a t^{-2}$, then $\rho$ is $-2$ and we obtain
$$
  {\oG_c\over\oF}(t)-C_\alpha \sim 
  {1\over 2t^2} \Bigl( -a(C_{\alpha+2}-C_\alpha) -\alpha (\alpha+1)
  C_{\alpha;2}\mu_{F,2}\Bigr) \, .
$$
If the constant in the second order terms vanishes (which generically
does not happen), then if we can, we need to add one more term when applying
Theorem {\MainTheorem}. But one should be careful that the rate of convergence
of $t^2g(t)$ to its limit may cancel the extra term added. If cancellation
occurs, then more terms are needed, and so on. It is therefore not clear
that second order regular variation provides the right framework for studying
second order expansions of $\oG_c$. In particular, for some exceptional
sequences of constants and some distributions, higher order regular variation
will be needed to obtain the exact second order. We also would like to point 
out that smooth variation of finite order is far easier to check than 
second order regular variation, and that it holds for most --- if not all ---
heavy tail distributions used in practical applications.

We conclude this section with a somewhat more general result to illustrate
how weights of arbitrary signs appear in the expansion. If $X$ is a random
variable with distribution function $F$, we write $F_*$ the distribution 
function of $|X|$. Thus, on the nonnegative half line, 
$\oF_*=\overline{\M_{-1}F}+\oF$. In the statistical literature dealing with 
regular variation for the upper and lower tails of distributions, it is 
customary to replace {\RVRemainder} by a second order regular variation
assumption on $\oF_*$, that is, with the same notation as in {\RVRemainder},
$$
  {\oF_*(\lambda t)\over \oF_*(t)}-\lambda^{-\alpha}
  \sim \lambda^{-\alpha} k(\lambda)g(t) \, ,
$$
and a tail balancing condition with remainder, that is for some nonnegative
$p$ at most $1$,
$$
  \oF = p\oF_*+o(\oF_*g)
$$
at infinity. Set $q=1-p$ and
$$
  \kappa_c(\lambda)=\sum_{i\in\ZZ} (\lambda/|c_i|)^{-\alpha} 
  k(\lambda/|c_i|)\bigl(p\II\{\, c_i>0\,\}+q\II\{\, c_i<0\,\}\bigr) \, .
$$
Define the constants
$$
  C_{+,\alpha}
  =\sum_{\matrix{\noalign{\vskip -2pt}\ss i\in\ZZ\cr
                 \noalign{\vskip -3pt}\ss c_i>0\cr}
         }c_i^{\alpha}
  \qquad \hbox{ and }
  C_{-,\alpha}
  =\sum_{\matrix{\noalign{\vskip -2pt}\ss i\in\ZZ\cr
                 \noalign{\vskip -3pt}\ss c_i<0\cr}
         }|c_i|^{\alpha} \, .
$$
Furthermore, set
$$
  C_{*,\alpha}=pC_{+,\alpha}+qC_{-,\alpha}
$$
and
$$
  C_{*,\alpha;1}
  =p(C_{+,\alpha+1}-C_1C_{+,\alpha}) -q(C_{-,\alpha+1}+C_1C_{-,\alpha}) \, .
$$
Assume that $\alpha$ is larger than $1$ and that $F$ satisfies the assumptions
of Theorem \MainTheorem with $m=1$ and $k=0$. Finally, assume also that
$a=\lim_{t\to\infty}tg(t)$ exists, possibly infinite. Then,
$$\displaylines{\quad
  {\oG_c(\lambda t)\over \oF_*(t)} - C_{*,\alpha}\lambda^{-\alpha}
  \hfill\cr\hfill
  \sim \cases{ \bigl( a \kappa_c(\lambda)-\lambda^{-\alpha-1}
                      \alpha\mu_{F,1}C_{*,\alpha;1}\bigr) t^{-1} 
               & if $a$ is finite,\cr
               \kappa_c(\lambda)g(t) & if $a$ is infinite.\cr}
  \cr}
$$
We remark that the above results extend Theorem 3.2.III in Geluk et al.~(1997).

\bigskip

\subsection{Some open questions.}%
Formula {\GAsymptTwo} provides generically a two terms expansion when the
mean $\mu_{F,1}$ vanishes. For this expansion to be valid, we need a variance
in order to define the second Laplace character. Consequently, we do not know
what a two terms expansion is when the distribution pertaining to $F$ is
centered and has infinite variance; that is essentially in the range $\alpha$
between $1$ and $2$.

When $\alpha$ is less than $1$, the mean does not exists, and we only have
the classical equivalence $\oG_c\sim\sum_{i\in\ZZ}\overline{\M_{c_i}F}$. Some
examples in subsection 3.1 and the log-gamma distribution studied in 
subsection 3.2 suggest that in some instances this equivalent may still provide
two terms or more. 

When $\alpha$ is less than $1$ and $F$ is concentrated on the 
nonnegative half line with smoothly varying tail of index
$\alpha$ and order more than $1$, Theorems 2.5 and 2.6 of Barbe and McCormick 
(2005) imply a result for some finite convolutions. Specifically, define
$$
  I(\alpha)= 2\int_0^{1/2} \bigl( (1-y)^{-\alpha}-1\bigr) \alpha y^{-\alpha-1}
  \d y + 2^{2\alpha}-2^{\alpha+1} \, .
$$
Then, if the $c_i$'s are nonnegative constants,
$$
  \overline{\star_{1\leq i\leq n}\M_{c_i}F} 
  = \sum_{1\leq i\leq n} \overline{\M_{c_i}F} + {I(\alpha)\over 2} 
  \Bigl( \sum_{1\leq i,j\leq n} c_i^\alpha c_j^\alpha -\sum_{1\leq i\leq n}
  c_i^{2\alpha}\Bigr) \oF^2 \, .
$$
This suggests that under suitable conditions, the tail of the series
$\sum_{i\in\ZZ} c_iX_i$ should behave like
$$
  \sum_{i\in\ZZ}\overline{\M_{c_i}F} 
  + {I(\alpha)\over 2} (C_\alpha^2-C_{2\alpha})\oF^2 +  o(\oF^2) \, .
  \eqno{\equa{OpenA}}
$$
The techniques used in Barbe and McCormick (2005) combined with those of the
current paper can certainly be used to prove \OpenA. However, even in the
case of a finite number of summands, when $\alpha$ is less than $1$, we do 
not know a good formalism to remove
the restriction that the support of the $\M_{c_i}F$'s should be in the 
nonnegative half line.

We believe that the techniques and formalism developed in the proof of
Theorem {\MainTheorem} of this paper are generally useful for extending a 
result from finite convolutions to infinite ones. Unfortunately, at the present
time, finite convolution are still difficult to work with.

\bigskip


\section{Applications.}

In this section we develop some applications of Theorem {\MainTheorem}
and the tail calculus explained in section 3.2. As the section goes,
these applications leave room to more and more questions; the last
subsection, on implicit renewal equation, only touches on a subject
which deserves further consideration.

\bigskip


\subsection{ARMA models.}%
ARMA models are among the most used models in statistical analysis
of time series. Yet, very few facts are known on their distributions. The
purpose of this subsection is to show that in some circumstances, 
Theorem {\MainTheorem} provides some basic information on the 
marginal distribution. We refer to Brockwell and Davis (1991)
for the basic probabilistic and statistical aspects of these models.

To fix notation, we define the backward shift operator $B$ on
sequences as follows. A sequence $x=(x_i)_{i\in\ZZ}$ is mapped under
$B$ to the sequence whose $i$-th element is $x_{i-1}$. As usual, $B^0$
is the identity, $B^{-1}$ is the inverse of $B$, and $B^k=B
B^{k-1}$. It then makes sense to consider polynomials in $B$. Having
two polynomials $\Theta$ and $\Phi$, and a sequence
$\epsilon=(\epsilon_i)_{i\in\ZZ}$ of independent and identically
distributed random variables, an ARMA process $X=(X_i)_{i\in\ZZ}$ is
defined by the relation
$$
  \Theta (B)X=\Phi(B)\epsilon \, .
$$
For this process to be defined, we assume that $\Theta$ has all its roots
outside the closed unit disk of the complex plane. We can define
$\Theta (B)^{-1}$ by a series expansion. Then, under a mild integrability
condition on $\epsilon$, we obtain 
$X=\Theta (B)^{-1}\Phi (B)\epsilon$. Having a series expansion yields a
representation of $X$ as an infinite order moving average
$$
  X_i=\sum_{j\in\NN} c_j\epsilon_{i-j} \, .
$$
In general it is unknown how to calculate the marginal distribution of
the process, that is the distribution of $X_0$. Since the roots of
$\Theta$ are outside the closed unit disk, the sequence $c_j$
decreases exponentially fast; see e.g.\ Brockwell et al.~(1991,\S
3.1). Therefore, this sequence has finite
$N_{\alpha,\gamma,\omega}$-norm, whatever $\alpha$, $\gamma$ and
$\omega$ are in the positive half line. Theorem {\MainTheorem} yields
immediately an expansion for the tail of the marginal distribution of
the process provided the distribution of the innovations $\epsilon_i$ has
tail smoothly varying of sufficiently large order and index.

To discuss further, let us write $G_c$ the marginal distribution of
the process. In order to make the expansion explicit, we need to
evaluate the Laplace character $\L_{G_c\natural \M_{c_i}F,m}$. This
requires computing the moments of $G_c\natural \M_{c_i}F$, or,
equivalently, $E(X_0-c_i\epsilon_{-i})^k$ for various integers $k$. We
do not know any way to obtain a nice formula for those moments in
terms of the polynomials $\Theta$ and $\Phi$. In general, one needs to
rely on numerical methods.

There are however two cases for which explicit calculations may be
performed, namely for AR($1$) and MA($q$) models, that is when either
$\Theta (B)=\Id -a B$ and $\Phi(B)=\Id$, or $\Theta (B)=\Id$ and $\Phi$
is an arbitrary polynomial of degree $q$. For instance, in many examples
of section 3, the expansions can be expressed with coefficients involving
the quantity $C_p=\sum_{i\in\ZZ} c_i^p$. For an AR($1$), with autocorrelation
$a$ positive and less than $1$,
$$
  C_p=\sum_{i\geq 0} a^{ip} = (1-a^p)^{-1} \, ,
$$
while for a MA($q$) with $\Phi (B)=\sum_{0\leq i\leq q} \Phi_i B^i$,
$$
  C_p=\sum_{0\leq i\leq q} \Phi_i^p \, .
$$
By applying these formulas to the expression obtained in section 3.3, we
obtain two terms expansions of the marginal distribution of these processes
with an assumption of second order regular variation for the distribution
of the innovations $\epsilon$. This will be used in the next subsection where
we develop an application to statistical inference for heavy tail data.

\bigskip


\subsection{Tail index estimation.}%
In heavy tail analysis, a critical parameter to estimate is the index of
regular variation when the marginal distribution of the data has a regularly
varying tail. In this example we are concerned with observations that follow
a causal linear process which we denote $(Y_i)_{i\in\NN^*}$. Thus, we assume
that there is a sequence of real number $c_i$'s and a sequence of independent
and identically distributed random variables $X_i$'s, so that
$$
  Y_i=\sum_{j\geq 0} c_j X_{i-j}\, , \qquad i\in\NN^* \, .
$$
By the first order result on tail behavior, we have that, if the distribution
function of the innovations $X_i$'s has regularly varying tails with index
$-\alpha$, then the same is true for the marginal distribution function 
of the stationary linear process. The statistical problem is
to estimate and find a confidence interval for $\alpha$, based on sample data.

The usual semiparametric procedure is to use the Hill (1973)
estimator, defined as follows. Let $|Y|_{i,n}$ be $i$-th largest value among
$|Y_1|,\ldots , |Y_n|$, so that $|Y|_{n,n}\leq \cdots \leq |Y|_{1,n}$.
Let $k_n$ be an integer between $1$ and $n$. The Hill estimator
is
$$
  \alpha_n = k_n \Bigl( \sum_{1\leq i\leq k_n} 
  \log {|Y|_{i,n}\over |Y|_{k_n+1,n}}\Bigr)^{-1} \, .
$$
For this estimate to be consistent, it is necessary to take points coming
from the tail of the distribution, and of course we need enough of them.
This is expressed by the conditions
$$
  \lim_{n\to\infty} k_n/n=0
  \qquad\hbox{ and }\qquad
  \lim_{n\to\infty} k_n=\infty \, .
$$
To derive a nondegenerate limiting distribution, and therefore obtain 
asymptotic
confidence intervals and tests, further hypotheses need to be imposed. These
were obtained by Resnick and St\u aric\u a (1997), and until the present paper
could not be verified except maybe in some very special situations. To explain
what the problem is, we need to list their conditions, and this requires
further notation as well some rather technical consideration.

Let $F_*$ be the distribution function of $|X_i|$. It is assumed 
that $F_*$ is second
order regularly varying; that is, with the same notation as in section 3.3,
$F_*$ belongs to $2RV(-\alpha,g)$ for some regularly varying function $g$,
i.e.~{\RVRemainder} holds with $\oF_*$ in place of $\oF$.
Next, it is assumed a tail balancing condition holds,
namely that for some $p$ in the closed unit interval and $q=1-p$,
$$
  \oF = p \oF_* +o(g\oF_*)
  \eqno{\equa{CBa}}
$$
and 
$$
  \overline{\M_{-1}F}=q \oF_* +o(g\oF_*)
  \eqno{\equa{CBb}}
$$
at infinity. This ensures that both $\oF$ and $\overline{\M_{-1}F}$ are
second order regularly varying, with same index $-\alpha$ and same auxiliary
rate function $g$. The next assumption made is that $F$ has a density $F'$ 
which is Lipschitz in mean, that is there exists a positive $\kappa$,
for which
$$
  \int_\RR |F'(x)-F'(x+y)| \d x \leq \kappa y \, .
  \eqno{\equa{CDF}}
$$
This condition ensures that the linear process is strong mixing.

Note that so far all the conditions are on the unknown distribution of the 
innovations. In statistical analysis of heavy tailed time series, one takes
as one's model assumptions that the innovation distribution satisfies certain
properties such as we have listed above. Model assumptions are simply assumed
to hold. However, conditions on the marginal distribution $G$ may be worrisome.
They may be redundant, i.e.~one derivable from the assumptions on $F$, or they
may even be inconsistent, i.e.~the hypotheses for a theorem may apply only to
an empty class of models. For our particular tail estimation problem, we will 
elucidate this issue next.

The next condition required is that the marginal distribution function $G$ 
of the process satisfies a Von Mises condition, that is has a density $G'$ and
$$
  \lim_{t\to\infty} {t G'(t)\over \oG(t)}=\alpha \, .
  \eqno{\equa{CGa}}
$$
Let $G_*$ be the distribution function of $|Y_i|$. It is further 
assumed that $\oG_*$ is second order regularly varying of index $-\alpha$, 
using the notation introduced after {\RVRemainder}, 
$$
  \oG_*\in 2RV(-\alpha,g_{G_*})
  \eqno{\equa{CGb}}
$$
for some regularly varying $g_{G_*}$.

Yet, an other assumption is that
$$
  \lim_{n_\to\infty} \sqrt{k_n}\, g_{G_*}\hskip-3pt\circ 
  G_*^\leftarrow(1-k_n/n)=  0 \, .
  \eqno{\equa{CGc}}
$$
Assumptions \CGa, {\CGb} and {\CGc} are rather problematic since they do not
involve the distribution function of the innovation, and may potentially 
put rather stringent conditions on $F$ or $k_n$. We will investigate that 
matter after stating Resnick and St\u aric\u a's result.

The last assumption is on the sequence $k_n$, and presents no difficulty,
since this sequence is chosen by whomever uses the estimator. It is assumed
that
$$
  \limsup_{n\to\infty} n^{2/3}/k_n < \infty
  \qquad\hbox{ or }\qquad
  \liminf_{n\to\infty} n^{2/3}/k_n > 0 \, .
  \eqno{\equa{Cknb}}
$$
Set
$$
  \lambda={1\over \alpha^2} \Bigl( 1+ 2{\sum_{j\geq 1}\sum_{k\geq 0}
  |c_k|^\alpha \wedge |c_{j+k}|^\alpha\over \sum_{k\geq 0} |c_k|^\alpha}\Bigr)
  \, .
$$

\Theorem{\hskip-5pt{\rm(Resnick, St\u aric\u a, 1997).}%
  \quad Assume that there exists some $u$ more than $1$ and some positive
  constant $A$ such that $|c_i|\leq A u^{-i}$ for any nonnegative $i$ and that
  $E|X_1|^d$ is finite for some positive $d$ less than $1$.
  
  \noindent (i) If conditions \CDF, {\CGa} and {\Cknb} hold, then
  $$
    \sqrt{k_n} \Bigl( \alpha_n^{-1}- {n\over k_n} \int_1^\infty
    P\bigl\{\, |Y_1|\geq x G_*^\leftarrow(1-k_n/n)\,\bigr\} {\d x\over x}\Bigr)
  $$
  has a centered normal limiting distribution, with variance $\lambda$.

  \noindent (ii) If, furthermore, {\CGb} and {\CGc} hold, then
  $\sqrt{k_n} (\alpha_n^{-1}-\alpha^{-1})$ has a normal limiting distribution
  with mean $0$ and variance $\lambda$.
}

\bigskip

Let us now show that Theorem {\MainTheorem} yields rather simple
and explicit conditions which ensure that {\CGa}, {\CGb} and {\CGc}
hold, and therefore makes Resnick and St\u aric\u a's theorem far easier
to use. Before stating the result we introduce some conditions on the
moving average weights. Recall at the end of section \fixedref{3.2},
we noted that in general higher-order expansions are required to
determine the auxiliary function for an infinite order weighted
average. In order that the derived auxiliary function be determined
by only the second order information for the innovation distribution,
certain restrictions must be met. They are encompassed in the
following conditions.

Similarly to $C_r$, we define the notation
$$
  |C|_r=\sum_{i\in\ZZ} |c_i|^r \, .
$$
Note that $|C|_r=|c|_r^r$.
Recall that $k(\cdot)$ is the function appearing in the definition of second
order regular variation for $\oF_*$ as in \RVRemainder. Furthermore, the
parameter $\rho$ which appear in the function $k$ is the index of
regular variation of the function $g$.
We will use the conditions
$$
  |C|_\alpha +\rho \sum_{i\in\ZZ} |c_i|^\alpha k(1/|c_i|) \not= 0 \, ,
  \eqno{\equa{TailCondA}}
$$
and, with $a$ a real number to be fixed later,
$$\displaylines{\quad
  a |C|_\alpha + a\rho \sum_{i\in\ZZ} |c_i|^\alpha k(1/|c_i|) 
  \hfill\cr\noalign{\vskip -3pt}\hfill
  {}+ \alpha\rho
  \mu_{F,1} (p-q) 
  \bigl( C_1\sum_{i\in\ZZ} |c_i|^\alpha\sign(c_i)-|C|_{\alpha+1}\bigr)
  \not= 0\, ,
  \quad\equa{TailCondB}\cr}
$$
as well as
$$
  a |C|_\alpha + a\rho \sum_{i\in\ZZ} |c_i|^\alpha k(1/|c_i|) 
  - \alpha(\alpha+1)\mu_{F,2} 
  \bigl( C_2|C|_\alpha-|C|_{\alpha+2} \bigr)\not= 0 \, .
  \eqno{\equa{TailCondC}}
$$
We comment on those conditions after the following statement.

\Proposition{%
  \label{PropTail}%
  Assume that $\oF$ and $\overline{\M_{-1}F}$ are smoothly varying of index
  $-\alpha$ and order more than $1$, and that they belong to $2RV(-\alpha,g)$. 
  Assume also that $F'$ ultimately exists and is continuous. Then,

  \smallskip

  \noindent (i) $G$ obeys the Von Mises condition \CGa;

  \smallskip

  \noindent (ii) If $F''$ exists and is Lebesgue integrable, then
  {\CDF} holds;

  \smallskip

  Next, let $\xi=1$ if $\mu_{F,1}$ does not vanish, and $\xi=2$ 
  otherwise. Assume furthermore that $\oF$ and $\overline{\M_{-1}F}$ are 
  of order more than $\xi$ and that $a=\lim_{t\to\infty} t^\xi g(t)$
  exists, possibly infinite. 

  \smallskip

  \noindent (iii) The function $\oG_*$ is second order regularly
  varying in any of the following three cases:
 
  \noindent\hskip1cm case 1. $a=+\infty$ and {\TailCondA} holds;

  \noindent\hskip1cm case 2. $a$ is finite, $\mu_{F,1}$ does not vanish, 
  and {\TailCondB} holds;

  \noindent\hskip1cm case 3. $a$ is finite, $\mu_{F,1}$ vanishes and
  {\TailCondC} holds.
 
  \smallskip

  \noindent (iv) Moreover,
  $$
    g_{G_*}\asymp\cases{ g          & if $a\not= 0$ \cr
                         \Id^{-\xi} & if $a=0$ \cr }
  $$  
  and
  
  \smallskip

  \noindent (v) condition {\CGc} is equivalent to
  $$
    \cases{ \lim_{n\to\infty} \sqrt{k_n}g\circ F_*^\leftarrow (1-k_n/n)=0
            & if $a\not=0$, \cr
            \lim_{n\to\infty} \sqrt{k_n} F_*^\leftarrow (1-k_n/n)^{-\xi}=0
            & if $a=0$. \cr}
  $$
}

\finetune{\vskip-15pt}

\Remark When $\mu_{F,1}$ does not vanish, the result is obtained by an 
application of Theorem {\MainTheorem} with $m=1$; otherwise
it requires $m=2$. The all but intuitive conditions {\TailCondA}--{\TailCondC} 
are not there for technical reasons. If those conditions are not
met as specified in the Proposition, then the last statement of the 
Proposition --- the equivalence with {\CGc} --- does not hold. We would
need to use Theorem {\MainTheorem} with a higher $m$ or higher order
regular variation, which would lead in some cases (but not always) to a
different result. Note that
generically, the Proposition covers all the cases; however, there are some
exceptional cases where it fails. This is the very same phenomenon as that
observed in subsection 3.1, and commented further toward the end of the proof
of Proposition \PropTail. 
This unfortunate fact commands caution when
estimating tail index in time series.

\Remark To fix the ideas, if $\oF(t)\asymp t^{-\alpha}$ and
$g(t)\asymp t^{-\beta}$ with $\beta$ positive and less than $\xi$. In
this situation, $a$ is infinite.  Then {\CGc} is equivalent to
$k_n=o(n^{2\beta/(\alpha+2\beta)})$. One sees that the smaller $\beta$
is, the smaller $k_n$ should be. In crafting a good estimator, one needs
to be mindful of such restrictions when using the Hill estimator in a time
series context.

\bigskip

\Proof (i) We first set
$$
  C_{+,\alpha}=\sum_{i\in\ZZ \, ;\, c_i>0}
  c_i^\alpha
  \qquad\hbox{ and }
  C_{-,\alpha}=\sum_{i\in\ZZ \, ;\, c_i<0}
  (-c_i)^\alpha
$$
Recall the classical first order equivalence
$$
  \oG
  \sim \sum_{i\in\ZZ} \overline{\M_{c_i}F}
  \sim (pC_{+,\alpha}+qC_{-,\alpha}) \oF_* \, .
$$
Applying Theorem {\MainTheorem} with $k=1$ and $m=0$ yields
$$
  \oG' 
  \sim \sum_{i\in\ZZ} \overline{\M_{c_i}F}{}'
  \sim (pC_{+,\alpha}+qC_{-,\alpha}) (-\alpha) \Id^{-1}\oF_* \, .
$$
We then deduce $\CGa$.

\noindent (ii) follows from the fundamental theorem of calculus as well 
as Fubini's theorem, writing
$$\eqalign{
  \int |F'(x)-F'(x+y)| \d x 
  &{}\leq \int\int |F''(x+u)|\d x\, \II\{\, 0\leq u\leq y\,\}\d u \cr
  &{}\leq y |F''|_1 \, . \cr
  }
$$

To prove the other statements, we first derive a two terms expansion for
$\overline{G_*}$. We write $G_c$ for $G$. Thus, $\M_{-1}G=G_{-c}$.

We first assume that $\mu_{F,1}$ does not vanish.
Since $\oG_*$ coincides with $\oG+\overline{\M_{-1}G}$ on the positive
half line, Theorem {\MainTheorem} and the tail balance conditions 
{\CBa}, {\CBb} imply
$$\eqalign{
  \oG_*
  & {}= \sum_{i\in\ZZ} \L_{G_c\natural \M_{c_i}F,1} \overline{\M_{c_i}F}
    +\sum_{i\in\ZZ} \L_{G_{-c}\natural \M_{-c_i}F,1} \overline{\M_{-c_i}F} 
    + o(\Id^{-1}\oF_*) \cr
  & {}=\sum_{i\in\ZZ} (\overline{\M_{c_i}F} + \overline{\M_{-c_i}F})
    - \mu_{F,1} \sum_{i\in\ZZ} (C_1-c_i) \D 
    (\overline{\M_{c_i}F}-\overline{\M_{-c_i}F}) \cr
  & \qquad{}+ o(\Id^{-1}\oF_*) \, .  \cr
  }
$$
We note that on the positive half line
$$
  \overline{\M_{c_i}F} + \overline{\M_{-c_i}F}
  =\overline{\M_{|c_i|}F_*}\, .
$$
Moreover, because $\oF$ and $\overline{\M_{-1}F}$ are smoothly varying of order
more than $1$, we also have
$$
  \D\overline{\M_{c_i}F} \sim -\alpha \Id^{-1}\overline{\M_{c_i}F} \, .
$$
We then obtain
$$\displaylines{\quad
  \oG_*=\sum_{i\in\ZZ} \overline{\M_{|c_i|}F_*}
  \hfill\cr\hfill
  {}+ \alpha\mu_{F,1}
  \sum_{i\in\ZZ} (C_1-c_i)(p-q)\sign (c_i)\Id^{-1}\overline{\M_{|c_i|}F_*} 
  + o(g) \oF_* \, .
  \cr}
$$
Consequently, we obtain the two terms expansion
$$\displaylines{\quad
  {\oG_*\over\oF_*}
  = |C|_\alpha +\sum_{i\in\ZZ} |c_i|^\alpha k(1/|c_i|) g
  \hfill\cr\hfill
  {}+ \alpha\mu_{F,1} (p-q)\Bigl(C_1\sum_{i\in\ZZ}\sign(c_i)|c_i|^\alpha 
  -|C|_{\alpha+1}\Bigr)
  \Id^{-1}
  + o(\Id^{-1}\vee g) \, .
  \cr}
$$

Let us now consider the case where $\mu_{F,1}$ vanishes.  
Again, Theorem {\MainTheorem}
and the tail balance conditions {\CBa}, {\CBb} imply
$$\eqalign{
  \oG_*
  & {}= \sum_{i\in\ZZ} \L_{G_c\natural \M_{c_i}F,2} 
    (\overline{\M_{c_i}F}+\overline{\M_{-c_i}F}) + o(\Id^{-2}\oF_*) \cr
  & {}=\sum_{i\in\ZZ} \overline{\M_{c_i}F} + \overline{\M_{-c_i}F}
    + {\mu_{F,2}\over 2} \sum_{i\in\ZZ} (C_2-c_i^2) \D^2 
    (\overline{\M_{c_i}F}+\overline{\M_{-c_i}F}) \cr
  & \qquad{}+ o(\Id^{-2}\oF_*) \, .  \cr
  }
$$
Because $\oF$ and $\overline{\M_{-1}F}$ are smoothly varying
of order more than $2$, we also have
$$
  \D^2\overline{\M_{c_i}F} 
  \sim \alpha (\alpha+1) \Id^{-2} \overline{\M_{c_i}F} \, .
$$
Then, up to $o(\Id^{-2}\oF_*)$, the tail $\oG_*$ is
$$\displaylines{\qquad
    \sum_{i\in\ZZ} \overline{\M_{|c_i|}F_*}  
    + {\mu_{F,2}\over 2} \alpha (\alpha+1) \sum_{i\in\ZZ}(C_2-c_i^2)
    \Id^{-2}\overline{\M_{|c_i|}F_*}
  \hfill\cr\hfill
    {}= \sum_{i\in\ZZ} \overline{\M_{|c_i|}F_*}
    + {\mu_{F,2}\over 2} \alpha (\alpha+1) (C_2|C|_\alpha-|C|_{\alpha+2})
    \Id^{-2} \oF_* \, .
  \qquad\cr}
$$
Therefore, we have the two terms expansion
$$\displaylines{\quad
  {\oG_*\over\oF_*}
  = |C|_\alpha + \sum_{i\in\ZZ} |c_i|^\alpha k(1/|c_i|) g
  \hfill\cr\hfill
  + {\alpha (\alpha+1)\over 2} \mu_{F,2} (C_2|C|_\alpha-|C|_{\alpha+2})\Id^{-2}
  + o(\Id^{-2}\vee g) \, .
  \quad\cr}
$$

In either the case $\mu_{F,1}$ vanishes or does not vanish, we obtained 
an expansion of the form
$$
  {\oG_*\over \oF_*} = U + Vg + W\Id^{-\xi} + o(\Id^{-\xi}\vee g)\, .
$$
This implies
$$\displaylines{\qquad
  {\oG_*(\lambda t)\over \oG_*(t)} 
  = \lambda^{-\alpha} + \lambda^{-\alpha} g(t)k(\lambda) (1+\rho U^{-1}V)
  \hfill\cr\hfill
  + \lambda^{-\alpha} U^{-1} W (\lambda^{-\xi}-1) t^{-\xi}
  + o\bigl(t^{-\xi}\vee g(t)\bigr) \, .
  \qquad\cr}
$$

We now prove the result in case 1. Indeed, 
if $a=\lim_{t\to\infty} t^\xi g(t)$ is infinite, we have
$$
  {\oG_*(\lambda t)\over \oG_*(t)} 
  = \lambda^{-\alpha} + \lambda^{-\alpha} g(t)k(\lambda) (1+\rho U^{-1}V)
  + o(g) \, .
$$
If $1+\rho U^{-1}V$ does not vanish, which is
condition {\TailCondA}, this implies $\oG_*$ is second order 
regularly varying with auxiliary function proportional to $g$.

The stated equivalence with {\CGc} follows in this case from the fact that
$\oG_*^\leftarrow(1-u)\asymp\oF_*^\leftarrow(1-u)$ as $u$ tends to $0$.

This proves statements (iii), (iv) and (v) in case 1.

Note that as we have seen in section
3.1, if $1+\rho U^{-1}V$ vanish, we cannot conclude anything without 
obtaining higher order expansions, and virtually any auxiliary function may
occur, either connected with higher order regular variation of $\oF_*$ (here
we are talking of third order or even higher order in exceptional cases) or
the third terms in the expansion of $\oG_*$ (or higher order terms in 
exceptional cases).

Cases 2 and 3 are handled in the same way.
\hfill$\qed$

\bigskip

\noindent{\it Example.}
To conclude this example we illustrate the result for a particular
distribution. Consider a Student innovation density 
$$
  f(x)=K_\alpha \Bigl( 1+{x^2\over\alpha}\Big)^{-(\alpha+1)/2} \, ,
  \qquad x\in\RR \, ,
$$
where $K_\alpha$ is the normalizing constant and where $\alpha$ is more than 
$2$.

The Student distribution being symmetric, {\CBa} and {\CBb} are obvious.

To check the assumptions of Resnick and St\u aric\u a's theorem
we use Proposition \PropTail.

It is plain that $\oF$ and $\overline{\M_{-1}F}$ are smoothly varying 
of order at least two, because the second order derivative of $f$ is 
regularly varying. Also, $f$ is continuously differentiable with integrable
derivative; this establish that $F''$ exists, is continuous and integrable.

To prove that $\oF$ is second order regularly varying, we derive a two
terms expansion, writing
$$\eqalign{
  \oF(t)
  &{}=K_\alpha \int_t^\infty x^{-\alpha-1} \alpha^{(\alpha+1)/2}
    (1+\alpha x^{-2})^{-(\alpha+1)/2} \d x \cr
  &{}=\alpha^{(\alpha+1)/2} K_\alpha \int_t^\infty 
    x^{-\alpha-1}-{\alpha(\alpha+1)\over 2} x^{-\alpha-3} + O(x^{-\alpha-5})
    \d x \cr
  &{}=\alpha^{(\alpha+1)/2} K_\alpha \Bigl( {1\over \alpha }t^{-\alpha}
    - {\alpha (\alpha+1)\over 2(\alpha+2)} t^{-\alpha-2}
    + O(t^{-\alpha-4})\Bigr)\, . \cr
  }
$$
Consequently,
$$\eqalign{
  \oF(\lambda t) -\lambda^{-\alpha} \oF(t)
  &{}\sim \alpha^{(\alpha+1)/ 2} K_\alpha {\alpha(\alpha+1)\over 2(\alpha+2)}
  (-\lambda^{-\alpha-2}+\lambda^{-\alpha}) t^{-\alpha-2} \cr
  &{}\sim \lambda^{-\alpha} {\alpha^2(\alpha+1)\over 2(\alpha+2)}
    \oF(t)t^{-2} (1-\lambda^{-2})\, . \cr
  }
$$
Therefore, $\oF$ belongs to $2RV(-\alpha,\Id^{-2})$.

A first order analysis shows that $G_*^\leftarrow(1-u)\asymp u^{1/\alpha}$.
Consequently, choosing $k_n=n^{4\theta /(4+\alpha)}$ with $\theta$ positive
less than $1$ ensures that the condition listed in Proposition 
\PropTail.{\it v}
is satisfied. For such choice, the second assumption in {\Cknb} holds.
We conclude that the distributional assumptions in the Resnick and 
St\u aric\u a theorem are satisfied.

An example of a process where the theorem leads to a fully explicit result is
the AR(1) model, $X_n=\sum_{j\geq 0} r^j Z_{n-j}$ with $|r|$ less than $1$.
In that case 
$$
  \lambda={1+|r|^\alpha\over \alpha^2(1-|r|^\alpha)} \, .
$$

\bigskip


\subsection{Randomly weighted sums.}%
In this subsection, we consider a weighted sum $\sum_{i\in\ZZ}W_iX_i$,
where the weights $W=(W_i)_{i\in\ZZ}$ are random, independent of the $X_i$'s.
We also write $\<W,X\>$ this series.
Clearly, under some assumptions on the weights, the uniformity of Theorem
{\MainTheorem} allows one to obtain an asymptotic expansion for the tail of
the weighted sum given the weights, and then decondition. This
can be achieved with various
integrability hypotheses on $W$, according to the arguments
used in the proof. The one which we provide seems to work well for the
applications which we will study. In particular, it does not add
any moment requirement to the distribution of the $X_i$'s.
In applications, it is often assumed that the weights are nonnegative. This 
is not strictly necessary for deriving tail expansions, but it somewhat 
simplifies the statements and proofs.

We will develop some applications in the next subsections.

Before stating our main result on randomly weighted sums, recall that
$|\cdot|_p$ is the $\ell_p$-norm on sequences. Hence, when the sequence
$W$ is nonnegative, $|W|_p$ is $(\sum_{i\in\ZZ}W_i^p)^{1/p}$. This is always
defined, possibly infinite.

\Theorem{%
  \label{RandomWeights}
  Let $F$ be a continuous distribution function, with tail $\oF$ in 
  $SR_{-\alpha,\omega}$ and such that 
  $\overline{\M_{-1}F}=O(\oF)$ at infinity. Let $m$ be an integer less than 
  $\alpha$ and $\omega$. Let $\gamma$ be a positive number less than $1$ and
  $\omega-m$. Assume that $F^{(m)}$ is bounded.
  Let $X=(X_i)_{i\in\ZZ}$ be a sequence of independent and identically
  distributed random
  variables having distribution $F$. Consider some random weights
  $W=(W_i)_{i\in\ZZ}$, independent of $X$, and such that for any
  $1\leq k\leq j\leq m$,
  $$
    E ( |W|_1^{j-k} \, |W|_k^k\,)
    \II\{\, N_{\alpha,\gamma,\omega}(W)\geq t\,\}
    =o\bigl(t^{-m}\oF(t)\bigr)
    \eqno{\equa{RWCondA}}
  $$
  as $t$ tends to infinity, and
  $$
    EN_{\alpha,\gamma,\omega}(W)^{m+\alpha+\epsilon} 
    < \infty \, .
    \eqno{\equa{RWCondB}}
  $$
  Let $K_W$ be the conditional distribution function of
  $\<W,X\>$ given $W$. Then, 
  $$
    P\{\, \<W,X\>\geq t \,\}
    = \sum_{i\in\ZZ} E \L_{K_W\natural \M_{W_i}F,m}\overline{\M_{W_i}F}
    (t) + o\bigl(t^{-m}\oF(t)\bigr)
  $$
  as $t$ tends to infinity.
}

\bigskip

\Remark In the proof of Theorem {\RandomWeights}, boundedness of $F^{(m)}$
is used only to prove that for any $0\leq j\leq m$, the map 
$(t,w)\mapsto w^{-j}\oF^{(j)}(t/w)$ is bounded on 
$[1,\infty)\times (0,\infty)$.
For the conclusion of the Theorem to hold, it is enough that there exist
$t_0$ and $M$ such that for any $t$ more than $t_0$, for any integer $i$
and for any integer $j$ at most $m$,
$$
  |W_i^{-j}\oF^{(j)}(t/W_i)|\leq M \quad\hbox{a.s.}
$$
In particular, when the $W_i$'s are less than a fixed number, this is 
implied by $\oF$ belonging to $SR_{-\alpha,\omega}$. For example,
the case of the $W_i$'s being Bernoulli random variables occurs in the
analysis of randomly stopped sums.

\bigskip

\Proof Let $R$ be $1/N_{\alpha,\gamma,\omega}(W)$. Define the sequence $c=RW$,
whose elements are $c_i=W_i/N_{\alpha,\gamma,\omega}(W)$. Since 
$N_{\alpha,\gamma,\omega}(\cdot )$ is homogenous of degree $1$, this new random
sequence satisfies $N_{\alpha,\gamma,\omega}(c)=1$. Let $G_c$ be the
conditional distribution function of $\<c,X\>$ given $c$.
Let $\epsilon$ be a positive real number. Let $t_2$ be as in Lemma {\Potter}
as applied to the normalized regularly varying function $\oF$.
Furthemore, let $t_1$ be at least $t_2$, and such that the function 
$\eta(\cdot)$ in Theorem {\MainTheorem} is at most $\epsilon$ on 
$[\, t_1,\infty)$. We apply Theorem {\MainTheorem} conditioning on $W$. So,
on the event $\{\,Rt > t_1\,\}$,
$$
  |\,\oG_c(Rt)-\sum_{i\in\ZZ} \L_{G_c\natural \M_{c_i}F,m}
  \overline{\M_{c_i}F}(Rt)\,| \leq (Rt)^{-m}\oF(Rt)\epsilon \, .
$$
Let $K$ be the (unconditional) distribution function of $\<W,X\>$. 
Clearly, $K(t)$ is the expected value of $G_c(Rt)$. Taking expectation 
with respect to the sequence $W$ in the previous inequality,
$$\displaylines{\quad
  |\,\oK(t)-\sum_{i\in\ZZ}E\L_{G_c\natural \M_{c_i}F,m}
  \overline{\M_{c_i}F}(Rt)\,|
  \hfill\cr\hfill
  \eqalign{
    {}\leq{}&\epsilon t^{-m}ER^{-m}\oF(Rt)\II\{\, Rt\geq t_1\,\}
             +E \oG_c(Rt)\II\{\, Rt\leq t_1\,\} \cr
            &{} + \sum_{i\in\ZZ} E|\,\L_{G_c\natural M_{c_i}F,m}
             \overline{\M_{c_i}F} (Rt)\,| \,\II\{\, Rt\leq t_1\,\} \, . \cr
    }
  \quad\equa{RWeqa}\cr
  }
$$
The equality
$$
  \L_{G_c\natural \M_{c_i}F,m}\overline{\M_{c_i}F}
  = \L_{\M_R(K_W\natural \M_{W_i}F),m} \M_R\overline{\M_{W_i}F}
$$
and Lemma {\LaplaceCharMul} show that the left hand side of {\RWeqa} is
the absolute value of $\oK(t)$ minus the asymptotic expansion given in the 
statement of the theorem.

To bound the right hand side of {\RWeqa}, the Potter bound in Lemma {\Potter}
implies that for $Rt$ and $t$ at least $t_1$,
$$
  \oF(Rt)/\oF(t) \leq R^{-\alpha} (R^{\epsilon}\vee R^{-\epsilon}) \, .
$$
Consequently,
$$
  ER^{-m}\oF(Rt) \II\{\, Rt\geq t_1\,\}
  \leq E R^{-m-\alpha} (R^\epsilon\vee R^{-\epsilon}) \oF (t) \, .
$$
Next, we have the obvious inequality
$$
  EG_c(Rt)\II\{\, Rt\leq t_1\,\}
  \leq \P{Rt\leq t_1}
  = P\{\, N_{\alpha,\gamma,\omega}(W)\geq t/t_1\,\} \, .
$$
Using our integrability assumption on the weights and Markov's inequality, 
this is at most $o\bigl(t^{-m}\oF(t)\bigr)$.

We now bound the third term in the right hand side of \RWeqa. Let $F_*$ be the
distribution function of $|X_i|$ and let $H_c$ be the conditional one 
of $\sum_{i\in\ZZ} c_i|X_i|$ given $c$.
Then, $\mu_{G_c\natural\M_{c_i}F,j}$ is at most $\mu_{H_c,j}$. Consequently,
$$\displaylines{\qquad
  E|\L_{G_c\natural \M_{c_i}F,m}\overline{\M_{c_i}F} (Rt)|
  \,\II\{\, Rt\leq t_1\,\}
  \hfill\cr\hfill
  {}\leq \sum_{0\leq j\leq m} E\mu_{H_c,j}|\D^j\overline{\M_{c_i}F}(Rt)|\,
  \II\{\, Rt\leq t_1\,\} \, .
  \qquad\equa{RWeqb}
  \cr}
$$
We then make use of the following claim which allows us to untangle the weights
and the random  variables $X$. It is a special case of a lemma in Chow
and Teicher (1978, \S 10.3). It controls the moments of a deterministically
weighted sum by that of the random variables and various $\ell_p$-norms
of the weights. Recall the notation $C_k=\sum_{i\in\ZZ} c_i^k$.

\bigskip

\Claim{ 
  Let $p_{k,j}=k-1+(j-k)(j-k+1)/2$.
  For any positive integer $j$ less than $\alpha$,
  $$
    \mu_{H_c,j}\leq \sum_{1\leq k\leq j} 2^{p_{k,j}}\mu_{F_*,k}\mu_{F_*,1}^{j-k}
    C_1^{j-k}C_k \, .
  $$
  }

\Proof The proof is that of Chow and Teicher (1978, \S 10.3) but with the
constant made explicit. For any nonnegative sequence $(a_i)_{i\in\ZZ}$,
$$\eqalign{
  \Bigl(\sum_{i\in\ZZ} a_i\Bigr)^j
  &{}=\sum_{i\in\ZZ} a_i\Bigl(a_i+\sum_{k\in\ZZ\setminus\{i\}} 
   a_k\Bigr)^{j-1} \cr
  &{}\leq 2^{j-1} \Bigl(\sum_{i\in\ZZ} a_i^j + \sum_{i\in\ZZ} a_i
   \Bigl(\sum_{k\in\ZZ\setminus\{ i\}}a_k\Bigr)^{j-1}\Bigr) \, .\cr 
  }
$$
Substituting $c_i|X_i|$ for $a_i$ in this inequality and
taking expectation on both sides with respect to $X$,
$$\eqalign{
  \mu_{H_c,j}
  &{}\leq 2^{j-1} \Bigl( \mu_{F_*,j}C_j+\mu_{F_*,1}\sum_{i\in\ZZ} c_i
     E\Bigl(\sum_{k\in\ZZ\setminus\{ i\}} c_kX_k\Bigr)^{j-1}\Bigr) \cr
  &{}\leq 2^{j-1} \Bigl(\mu_{F_*,j}C_j +\mu_{F_*,1}
     \sum_{i\in\ZZ} c_i\mu_{H_c,j-1}\Bigr) \cr
  &{}=2^{j-1} (\mu_{F_*,j}C_j+\mu_{F_*,1}C_1\mu_{H_c,j-1})\, . \cr
  }
$$
The claim follows by induction.\hfill$\qed$

Continuing the proof of Theorem \RandomWeights, that is, to bound {\RWeqb} 
from above, we have
$$
  \D^j\overline{\M_{c_i}F}(Rt)
  = R^{-j}W_i^{-j} \oF^{(j)}(t/W_i) \, .
$$
Since $\oF^{(j)}$ is regularly varying with index $-\alpha-j$ and bounded, 
since we also can assume that $t$ is at least $1$ say, $\oF^{(j)}(t/w)$
is at most a constant times $(w\wedge 1)^{j+\alpha-\epsilon}$. In particular,
the map $(t,w)\mapsto w^{-j}\oF^{(j)}(t/w)$ is bounded on $[\,1,\infty)\times
(0,\infty)$. This fact, combined with the claim, show that {\RWeqb} is at most
$$
  O(1)\sum_{0\leq j\leq m} \sum_{1\leq k\leq j} E C_1^{j-k}C_k R^{-j}
  \II\{\, Rt\leq t_1\,\} \, ,
$$
that is
$$
  O(1) \sum_{0\leq j\leq m} \sum_{1\leq k\leq j}
  E\Bigl(\sum_{n\in\ZZ} W_n\Bigr)^{j-k}
  \sum_{n\in\ZZ} W_n^k \,\II\{\, Rt\leq t_1\,\} \, .
$$
But one of our integrability assumptions implies that this last expression
is $o\bigl(t^{-m}\oF(t)\bigr)$, which concludes the proof. \hfill$\qed$

\bigskip


\subsection{Randomly stopped sums.}%
Theorem {\RandomWeights} has many applications in applied probability. 
In this subsection, we obtain a tail expansion for randomly stopped 
sums. Randomly stopped sums are a basic model in insurance mathematics,
e.g.\ in modelling total claim size. A discussion of asymptotic behavior
of random sums may be found in Embrechts et al.\ (1997). Some practical
methods for obtaining tail area approximations for compound distributions 
may be found in Beirlant et al.\ (1996). Willmot and Lin (2000) is a good
source of information for compound distributions and may be consulted for
additional references; see also Willekens (1989). We further mention Omey
and Willekens (1986, 1987) who obtained second-order results. We also mention
Geluk (1992, 1996) who provides second order results for subordinated 
probability distributions in the heavy tail case.

Let $N$ be a nonnegative random variable, independent of the $X_i$'s. Define
the sum $S_n=X_1+\cdots + X_n$, with $S_0=0$. Consequently, we agree that
$F^{\star 0}$ is the point mass at $0$. We write $K$ for the distribution
function of $S_N$.

\Theorem%
{\label{RandomlyStoppedSum}
  Assume
  that $F$ is a distribution function on the nonnegative half line, such 
  that $\oF$ 
  belongs to $SR_{-\alpha,\omega}$. Let $m$ be an integer less than $\alpha$ 
  and $\omega$. Let $\gamma$ be a positive number less than $\omega-m-k$
  and $1$. If $N$ has a moment of order more than
  $$
    {\alpha+m\over\gamma} \Bigl({\alpha+\omega\over\alpha}\vee 2\Bigr) + m \, ,
  $$
  then,
  $$
    \oK = E N\L_{F^{\star (N-1)},m} \oF + o(\Id^{-m}\oF) 
  $$
  at infinity.
}

\bigskip

\Proof Let $W$ be the random sequence defined by $W_i=\II\{\, 0<i\leq N\,\}$.
Then $S_N=\<W,X\>$. To apply Theorem {\RandomWeights} with the remark following
it, we need to check its assumptions, which is to check the integrability 
condition on the weights. Since
$$
  {1\over\gamma}\Bigl( {\alpha\over \alpha+\omega}\wedge {1\over 2}\Bigr)^{-1}
  = {1\over\gamma} \Bigl( {\alpha+\omega\over\alpha}\vee 2\Bigr) \, ,
$$
we have
$$
  N_{\alpha,\gamma,\omega}(W)
  =N^{{\ss 1\over\ss\gamma}({\ss\alpha+\omega\over\ss\alpha}\vee 2)} 
  \vee 2^{\ss\alpha\over\ss\alpha+\omega} \, .
$$
Moreover, $|W|_1=N$ and $|W|_k^k=N^k$. So the integrability conditions are
implied by 
$$
  EN^{{\ss m+\alpha+\epsilon\over\ss \gamma} 
  ({\ss\alpha+\omega\over\ss\alpha}\vee 2)}
  <\infty \, ,
$$
and
$$
  E N^m
  \II\{\, N>t^{\gamma({\alpha\over\alpha+\omega}\wedge {1\over 2})}\,\}
  = o\bigl(t^{-m}\oF(t)\bigr) \, .
$$
The first expectation is finite by assumption. For the second one, we
apply the standard trick to prove Markov's inequality. For any $p$, the 
expectation is at most
$$
  E N^{m+p}
  t^{-p\gamma ({\ss\alpha\over\ss \alpha+\omega}\wedge {1\over 2})} \, .
$$
Take $p$ such that the exponent of $t$ is less than $-\alpha-m$ but such
that $EN^{m+p}$ is finite.
Use Proposition 1.3.6 in Bingham, Goldie and Teugels (1989) to conclude
that the assumptions of Theorem {\RandomWeights} hold. To conclude, note
that for $i$ positive and less than $N$, the equalities $\M_{W_i}F=F$ and 
$K\natural \M_{W_i}F=F^{\star (N-1)}$ hold.
\hfill$\qed$

\bigskip

As in Barbe and McCormick (200?), using the Laplace transforms of $X_1$ and $N$
allows one to derive a rather neat expression. Indeed, setting
$$
  \Lambda_X(t) = E e^{-tX} \quad\hbox{ and }\quad\Lambda_N(t)=E e^{-tN} \, ,
$$
equality (2.2.1) in Barbe and McCormick (200?) yields
$$
  EN\L_{F^{\star (N-1)},m} 
  = -\sum_{0\leq j\leq m} {1\over j!} {\d^j\over \d u^j} 
  { \Lambda_N'\bigl(-\log \Lambda_X(u)\bigr)\over \Lambda_X(u)} \Bigr|_{u=0}
  \D^j \, .
  \eqno{\equa{EquaRandomStop}}
$$
Then, the technique explained in Barbe and McCormick (200?) allows for 
efficient computation using computer algebra packages.

Theorem {\RandomlyStoppedSum} has an interesting special case.

\Corollary%
{\label{CorRandomlyStoppedSum}
  Let $F$ be a distribution function satisfying the assumptions of Theorem 
  \RandomlyStoppedSum. If $N$ has a Poisson distribution with parameter $a$,
  then
  $$
    \oK=a \L_{K,m}\oF + o(\Id^{-m}\oF) \, .
  $$
}

\Proof The equality
$$
  \sum_{k>0} k S_{k-1}^j {a^k\over k!} e^{-a}
  = a\sum_{k\geq 1} S_{k-1}^j {a^{k-1}\over (k-1)!} e^{-a}
$$
shows that 
$$
  EN\L_{F^{\star (N-1)},m}=aE\L_{F^{\star N},m}=a\L_{K,m} \, .
$$
The result follows from Theorem \RandomlyStoppedSum.\hfill$\qed$

\bigskip


\subsection{Queueing theory.}
We consider a queue of M/G/1 type. This means that customers arrive
with interarrival time exponentially distributed with mean $\mu$; the service
has general distribution function $B$ with finite mean $\beta$.
We assume in this application that the system is stable, that is
$\beta/\mu$ is less than $1$. Define
$$
  H(t)=\beta^{-1}\int_0^t\oB(s)\d s \, \quad t\geq 0\, .
$$
Set $a=\beta/\mu$.
The Pollaczek-Hincin formula shows that the waiting time has the
compound geometric distribution function
$$
  W=(1-a) \sum_{n\geq 0} a^nH^{\star n} \, .
$$
We refer to Bingham, Goldie and Teugels (1987, p.387) or Cohen (1972)
for this derivation. If $\oB$ is regularly varying at infinity, Cohen (1972)
shows that so are $\oH$ and $\oW$. When $\oB$ statisfies the assumption
of Theorem $\MainTheorem$, higher order expansions can be obtained. We mention
Abate et al.\ (1994), Abate et al.\ (1995), Abate and Whitt (1997) for recent
related work. See also Willekens and Teugels (1992).

\Proposition
{\label{PropQueue}
  Assume that $B$ satisfies the assumptions of Theorem {\RandomlyStoppedSum}.
  Let $\Lambda_H$ be the Laplace transform of $H$. Then
  $$\displaylines{\qquad
    \oW =a(1-a)\sum_{0\leq j\leq m} {1\over j!}{\d^j\over \d u^j}
     \bigl(1-a\Lambda_H(u)\bigr)^{-2}\Big|_{u=0} \oH^j 
     \hfill\cr\hfill
     +o(\Id^{-m}\oH)
     \qquad\cr}
  $$
  at infinity.
}

\bigskip

\Proof Let $S_n$ be a sum of $n$ independent and identically
distributed random variables having distribution function $H$. We set
$S_0=0$. We see that $W$ is the distribution of $S_N$ where
$N$ has a geometric distribution with parameter $a$. The Laplace
transform of $N$ is
$$
  \Lambda_N(u)=(1-a)/(1-ae^{-u}) \, .
$$
Therefore, for a random variable $X$ having distribution $H$,
$$
  {\Lambda_N'\bigl(-\log \Lambda_X(u)\bigr)\over\Lambda_X(u)}
  = {-a(1-a)\over \bigl( 1-a\Lambda_X(u)\bigr)^2 } \, .
$$
Applying formula \EquaRandomStop, we obtain the result.\hfill$\qed$

\bigskip


\subsection{Branching processes}%
Consider an age dependent branching process. Basically,
this refers to a Galton-Watson branching process in which the particles
have a random life time governed by a distribution $F$. The process
starts at time $0$ with one particle; at the end of its life, it generates
an average of $\rho$ particles, which themselves at death will each
generate an average of $\rho$ particle, and so on. We refer to Athreya
and Ney (1972, chapter 4) for a complete description of the process.
It is intuitively clear that if $\rho$ is less than $1$, the so-called
subcritical case, then the process will become extinct.

Let $\nu(t)$ be the expected number of particles alive at time $t$. In the
subcritical case, let also
$N$ be a geometric random variable with parameter $\rho$, that is
$N$ is a nonnegative integer $k$ with probability $(1-\rho)\rho^k$.

\Proposition{
\label{PropBranching}
  Assume that $F$ is a continuous distribution function on the positive
  half line, whose tail $\oF$ belongs to $SR_{-\alpha,\omega}$. Let $m$
  be an integer less than $\alpha$ and $\omega$. Then
  $$
    \nu= \rho^{-1} E N\L_{F^{\star (N-1)},m} \II\{\, N\geq 1\,\}
    \oF  + o(\Id^{-m}\oF)
  $$
  at infinity.
  }

\bigskip

\Proof As before, let us agree that $F^{\star 0}$ is the distribution
function of the point mass at $0$.
Equation (4) in section IV.5 of Athreya and Ney (1972) shows that
$$
  \nu=\sum_{k\geq 0} \rho^k (\overline{F^{\star (k+1)}}-\overline{F^{\star k}})
  \, .
$$
Applying Theorem \RandomlyStoppedSum, we see that
$$
  \sum_{k\geq 0} \rho^k \overline{F^{\star (k+1)}} 
  = \rho^{-1} (1-\rho)^{-1}E( \overline{F^{\star N}}\II\{\, N\geq 1\,\} )
$$
has asymptotic expansion
$$
  \rho^{-1}(1-\rho)^{-1}EN\L_{F^{\star (N-1)},m}\II\{\, N\geq 1\,\} \oF
$$
while 
$$
  \sum_{k\geq 0} \rho^k\overline{F^{\star k}}
$$
has asymptotic expansion
$$
  (1-\rho)^{-1}ENL_{F^{\star (N-1)},m} \II\{\, N\geq 1\,\} \oF \, .
$$
The result follows.\hfill$\qed$

\bigskip

Let us now present an explicit expansion with $m=2$. We write $\sigma_F^2$ the
variance pertaining to the distribution $F$. Direct calculation
yields
$$\displaylines{\quad
  \rho^{-1} EN\L_{F^{\star (N-1)},2} = {1\over 1-\rho}\Id 
  - {2\rho\mu_{F,1}\over (1-\rho)^2} \D
  \hfill\cr\hfill
  {}+ {\rho\over (1-\rho)^3} 
  \bigl( (1-\rho)\sigma_F^2+(1+2\rho)\mu_{F,1}^2\bigr)
  \D^2 \, ,
  \cr}
$$
so that
$$\displaylines{\quad
  \nu
  = {1\over 1-\rho}\oF
  + {2\rho\mu_F\over (1-\rho)^2}F'
  \hfill\cr\hfill
  {}- {\rho\over (1-\rho)^3} 
  \bigl( (1-\rho)\sigma_F^2+(1+2\rho)\mu_{F,1}^2\bigr)
  F'' + o(\Id^{-2}F) \, .
  \cr}
$$
The first order term yields the result in Chover et al.~(1973, p.296) and the
first two terms yield the second-order result result in Grubel (1987). We
refer to Pakes (1975) for related work in a subexponential setting.

Note that if the Laplace transform of $F$ is known, one could use Proposition
{\PropQueue} to obtain an alternative form of the asymptotic expansion.

\bigskip


\subsection{Infinitely divisible distributions.}%
The infinite divisible distributions are those which for any integer
$n$ can be written as $n$-th convolution power of another
distribution.  They are also the limiting distributions of sums of $n$
independent and identically distribution random variables when the
distribution is allowed to change with $n$.  We refer to Feller (1971)
for an introduction to the topic. These distributions are
characterized through their characteristic functions, which are of the
form
$$
  \zeta\in\RR\mapsto\exp\Bigl( i\tau\zeta-{\sigma^2\over 2}\zeta^2+\int\Bigl(
  e^{i\zeta x}-1-i{\zeta x\over 1+x^2}\Bigr) \d \nu(x)\Bigr) \, ,
$$
where $\sigma^2$ is nonnegative, $\tau$ is a real number, $\nu$, the
so-called L\'evy measure, has no mass at the origin and satifies 
$\int x^2/(1+x^2) \d\nu(x)<\infty$. Note that
throughout this subsection, $i$ denotes the imaginary unit; this is the
only part in this paper where complex numbers are used.

For what follows, it is convenient to set
$$
  \onu(t)=\nu(t,\infty) \, .
$$

Let $G_\nu$ be an infinitely divisible distribution function with
L\'evy measure $\nu$.  If $\nu$ has a regularly varying tail, then
$\oG\sim\onu$; see e.g.~Feller (1971); see also Pakes (2004) for a
related first-order result and Embrechts et al.\ (1979) for work in
the subexponential case.  Under a slightly stronger assumption,
Gr\"ubel (1987) proved the two terms expansion
$$
  \oG_\nu \sim \onu - \mu_{G,1}\D\onu
$$
at infinity. With our notation this means $\oG_\nu\sim \L_{G_\nu,1}\onu$.

The following result is then quite natural.

\Proposition
{\label{PropID}
  Let $\nu$ be a measure such that $\onu$ is smoothly varying
  of index $-\alpha$ and order $\omega$ and $\nu (-\infty,-t )
  =O\bigl(\onu(t)\bigr)$ as $t$ tends to infinity. For any 
  integer $m$ less than $\alpha$ and $\omega$,
  $$
    G_\nu=\L_{G_\nu,m}\onu + o(\Id^{-m}\onu ) \, .
  $$
}

In this expansion, the Laplace character $\L_{G_\nu,m}$ involves the
moments of $G_\nu$ of order at most $m$. These are finite when $\onu$
is regularly varying of index less than $-m$ and the right tail of
$\nu$ is dominant, that is under the assumptions of the
Proposition. In practice, those moments can be computed by
differentiating the characteristic function of $G_\nu$. Again,
computer algebra packages are wonderful for doing this type of work.

\bigskip

\Proof The basic idea underlying the proof is a classical representation
of infinitely divisible distributions as convolutions of two distributions,
the first one having light tail, the second one being that of a randomly
stopped sum with heavy tail. This makes it possible to use Theorem 
\RandomlyStoppedSum .

Let $A$ be a positive constant, and let $a$ be the $\nu$-measure of 
$[\,-A,A\,]^c$. Let $\nu_1$ be the measure $\nu(\cdot\cap [\,-A,A\,]^c)$,
and let $F$ be the distribution function pertaining to the probability 
measure $\nu_1/a$. Define
$$
  \tau_1=\tau-\int_{[-A,A]^c} {x\over 1+x^2}\d\nu(x) \, .
$$
The function
$$
  \widehat H(\zeta) 
  =\exp\Bigl( i\zeta\tau_1-{\sigma^2\over 2}\zeta^2+\int_{-A}^A
  e^{i\zeta x}-1-{i\zeta x\over1+x^2} \d\nu(x)\Bigr)
$$
is the characteristic function of an infinitely divisible distribution
function $H$. Write $\widehat F$ for the characteristic function of $F$ and 
$\widehat{G_\nu}$ for that of $G_\nu$. One has
$$
  \widehat{G_\nu}=\widehat H\exp \bigl(a(\widehat F-1)\bigr) \, .
  \eqno{\equa{IDa}}
$$
Let $X=(X_i)_{i\geq 1}$ be a sequence of independent and identically
distributed random variables, all having distribution $F$. They induce
the partial sums $S_n=X_1+\cdots + X_n$, with the usual convention
$S_0=0$. Let $N$ be a Poisson random variable, with mean $a$. Let $K$
be the distribution function of the randomly stopped sum $S_N$. One
can check that the characteristic function of $K$ is
$$
  \widehat K
  =\exp \bigl( a(\widehat F-1)\bigr) \, ,
$$
so that $\widehat{G_\nu}=\widehat H\widehat K$. Consequently, $G_\nu=H\star K$.

It follows from Corollary {\CorRandomlyStoppedSum} that $\oK\sim a\L_{K,m}\oF$.

The remainder of the proof is the only place in section 4 where we rely on 
results which we will establish in the proof of Theorem \MainTheorem.
We write
$$
  \oG(t)=\int_{-\infty}^{t/2} \oK(t-x)\d H(x)
  + \int_{-\infty}^{t/2} \oH(t-x) \d K(x)
  + \oH\,\oK(t/2) \, .
$$
As observed by Gr\"ubel (1987, proof of Theorem 7), the function $\oH$ decays 
exponentially fast to $0$ at infinity. So
$$
  \oG(t)=\int_{-\infty}^{t/2} \oK (t-x)\d H(x) + o\bigl(t^{-m}\oF(t)\bigr) 
  \, .
$$
If follows from Theorem \fixedref{5.4.1}, Lemma \fixedref{6.2.1}
and the asymptotic expansion of $\oK$ that
$$\eqalign{
  \oG(t)
  &{}=\sum_{0\leq j\leq m} {(-1)^j\over j!} a \mu_{K,j}\int_{-\infty}^{t/2}
    \D^j \oF(t-x)\d H(x) + o\bigl( t^{-m}\oF(t)\bigr) \cr
  &{}=\sum_{0\leq j\leq m} {(-1)^j\over j!} a \mu_{K,j}\L_{H,m-j}\D^j\oF(t)
    + o\bigl( t^{-m}\oF(t)\bigr) \, . \cr
  }
$$
Applying Lemma \LaplaceBinomial, we conclude that
$$
  \oG = a\L_{K\star H,m}\oF + o(\Id^{-m}\oF) \, .
$$
Since $K\star H$ is $G_\nu$ and $a\oF$ ultimately coincide with $\onu$,
this concludes the proof.\hfill$\qed$

\bigskip


\subsection{Implicit transient renewal equation and iterative systems.}%
A renewal equation is an integral equation of the form 
$$
  F-aF\star H=K \, ,
$$
where $H$ and $K$ are given distribution functions, $a$ is a real
number and $F$ is the unknown. Such equation is transient if the
absolute value of $a$ is less than $1$.  These equations arise in
applied probability and we refer to Feller (1979) for an introduction
to renewal theory and Bingham, Goldie and Teugels (1989, \S 8.6) for
basic results and references on the tail behaviour of the
solutions. Following Goldie (1991), implicit renewal equations are
equations of the same form, except that $K$ is an integral transform
of $F$. Again, these equations appear in applied probability and
statistics, in connection with the stationary solution of iterative
systems.

The purpose of this subsection is to show that in some cases, one can
derive an asymptotic expansion of the solution by solving linear
equations.  This method is radically different from that of Kesten
(1973) or Goldie (1991); it works under a different set of
assumptions, closer to the one used by Grey (1994).  The appealing
feature of this approach is its simplicity of implementation.  The
idea is very simple: expand the known function on an asymptotic scale;
on this scale, the implicit renewal equation becomes a finite
dimensional linear system, which can be solved by elementary linear
algebraic techniques.  Consequently, we obtain an analogue of the
technique to solve differential equations using formal series
expansions and identifying the coefficients.  Of course, this
technique will only work for a limited type of renewal equation.

To explain the principle, we first look at a simple renewal equation.

\medskip

\noindent{\it Standard renewal equation.} Assume that $H$ and $K$ have moments
of order $m$. Then the renewal equation gives $m$ equations determining the
moments of $F$, that is
$$
  \mu_{F,j}-a\sum_{0\leq i\leq j} {j\choose i} \mu_{F,i}\mu_{H,j-i}
  = \mu_{K,j} \, ,\qquad 1\leq j\leq m \, .
$$
This can be rephrased in a nicer form with Proposition \LaplaceRepresentation,
which combined with the renewal equation imply
$$
  \L_{F,m}-a\L_{F,m}\circ\L_{H,m}=\L_{K,m} \, .
$$
When the absolute value of $a$ is less than $1$, the operator $\Id-a\L_{H,m}$
is not in the ideal generated by $\D$. Therefore, it is invertible in 
$(\RR_m[\,\D\,],\circ)$. This implies that for $|a|$ less than $1$,
$$
  \L_{F,m}=(\Id-a\L_{H,m})^{-1}\circ \L_{K,m} \, .
  \eqno{\equa{RenewalLaplace}}
$$
So, we can calculate the Laplace character and the moments of $F$. 
This calculation is done by manipulating finite dimensional matices.

Iterating the renewal equation yields
$$
  F=\sum_{i\geq 0} a^i K\star H^{\star i} \, .
  \eqno{\equa{RenewalRepF}}
$$
Hence, $F$ converges to $1/(1-a)$ at infinity and its tail is given by
$$
  (1-a)^{-1}-F = \sum_{i\geq 0} a^i \overline{K\star H^{\star i}} \, .
$$
Define the distribution function $G$ by $\oG=1-(1-a)F$. It solves the equation
$$
  \oG - a\overline{G\star H}=(1-a)\oK \, .
  \eqno{\equa{RenewalG}}
$$
Similarly to what we did in subsection 3.2, assume that $K$ and $H$ have
an asymptotic expansion on a $\star$-asymptotic scale $e$. It is 
then conceivable that $\oG$ has
an expansion in that scale. Then {\RenewalG} and Theorem {\MainTheorem}
yield, with notation analogous to that of section 3,
$$
  p_{\oG} - a (\calL_G p_{\oH} + \calL_H p_{\oG})=(1-a)p_{\oK} \, .
$$
That is, $p_{\oG}$ is the solution of a linear system of equations, and can be
made explicit by the formula
$$
  p_{\oG} = (\Id-a\calL_H)^{-1} \bigl( (1-a)p_{\oK}+a\calL_G p_{\oH}\bigr) 
  \, .
  \eqno{\equa{RenewalSolpG}}
$$
In this formula, $\calL_G$ involves the moments of $G$ and those are $1-a$
times the moments of $F$; these can be calculated by solving the
linear system \RenewalLaplace.
Having an expansion for $\oG$ it is then trivial to derive one for $F$.

The only thing missing to make this rigorous is to prove that $G$ has
an asymptotic expansion in the same asymptotic scale $e$ as $H$ and
$K$.  This is trivial because, when $a$ is positive, representation
{\RenewalRepF} shows that $G$ is the distribution function of the sum
of a random variable with distribution $K$ and a randomly stopped sum
having a number of summands distributed according to a geometric
distribution with parameter $a$; and then Theorem
{\RandomlyStoppedSum} applies.  When $a$ is negative, one uses the
same argument but split the sum in {\RenewalRepF} into one where the
index $i$ is odd and one where the index is even.

\medskip

\noindent{\it An implicit renewal equation.} Following Goldie (1991), 
Grey (1994), Grincevi\v cius (1975) and Kesten (1973), consider
the distributional equation involving random variables, $R\equald  Q+MR$
where $(M,Q)$ and $R$ are independent and all random variables are 
nonnegative. We need one important extra assumption, namely that $M$ and
$Q$ are independent. It would be desirable not to make this assumption,
but this would require some interesting generalization of 
Theorem \MainTheorem. To write the implicit renewal
equation equivalent to this distributional identity, let $\convM$ denote the
Mellin-Stieltjes convolution. That is, if $F$ and $H$ are two distribution
functions,
$$
  H\convM F(t)=\int F(t/x)\d H(x) \, .
$$
The connection with iterative systems is that if $(Q_i,M_i)$, $i\geq 1$, are
independent, all distributed as $(Q,M)$, then, under suitable conditions, the
sequence $R_0=0$, $R_n=Q_n+M_nR_{n-1}$ has a distribution converging
to that of $R$, that is to $F$ --- see for instance the survey article
by Diaconis and Freedman (1999). We mention that first-order asymptotics
for random coefficient autoregressive models is consider in Resnick and
Willekens (1991).

Let $H$ and $K$ be the distribution functions of $M$ and $Q$ respectiveley. 
The distributional equation is equivalent to
$$
  F=K\star (H\convM F) \, .
  \eqno{\equa{ImplicitRenewalA}}
$$

Let us assume that $K$ has a tail expansion in a $\star$-asymptotic scale $e$
as in section 3.2. Equation {\ImplicitRenewalA} yields two equations:
one on the Laplace characters, which allows us to identify the moments of $F$,
one on the tail vectors $p_{\oF}$, $p_{\oK}$, which
allows one to find the tail expansion of the solution.

For any integer $k$ for which the $k$-th moments of $K$ and $H$ exist,
the moment equations are
$$\eqalign{
  \mu_{F,k}
  &{}= \sum_{0\leq j\leq k} {k\choose j} \mu_{K,j}\mu_{H\convMss F,k-j} \cr
  &{}=\sum_{0\leq j\leq k} {k\choose j} \mu_{K,j} \mu_{H,k-j}\mu_{F,k-j} \, . 
   \cr
  }
$$
These equations allow one to find the moments $\mu_{F,k}$ by induction on $k$.

If $\oF$ has an expansion in the scale $e$, then {\ImplicitRenewalA} and 
Theorem {\MainTheorem} suggest that
$$
  p_{\oF} = \calL_{H\convMss F} p_{\oK} + 
  \calL_K p_{\overline{H\convMss F}} 
  \, .
$$
Since the variables are nonnegative,
$$
  \overline{H\convM F}(t)
  = \int \oF(t/x) \d H(x) 
  = \int \overline{\M_xF}(t)\d H(x) \, .
$$
Consequently,
$$
  p_{\overline{H\convMss F}} = \int\calM_x \d H(x) \, p_{\oF} \, .
$$
Thus, we obtain the tail expansion of $F$ in the scale $e$, whose coefficients
are given by
$$
  p_{\oF} = \Bigl( \Id - \calL_K\int\calM_x \d H(x)\Bigr)^{-1}
  \calL_{H\convMss F} p_{\oK} \, .
  \eqno{\equa{ExpansionRenewal}}
$$
Again, for this to be justified, we only need to prove that $\oF$ has an
expansion in the scale $e$. Our next result gives a sufficient condition,
and a complete example of application follows its proof.

\Proposition{%
  \label{ImplicitRenewalTheorem}
  Assume that $K$ satistifies the assumptions of Theorem \RandomWeights.
  If $EM_1^{2(\alpha+m+1)}$ is less than $1$, then formula 
  {\ExpansionRenewal} holds.
}

\medskip

Note that by Jensen's inequality, the integrability assumption in Proposition
{\ImplicitRenewalTheorem} implies that $E\log M_1$ is negative, and therefore,
that the implicit renewal equation has a well defined solution.

\bigskip

\Proof Let $Q_i$ (respectively $M_i$), $i\geq 1$, be a sequence of
independent random variables all having the distribution $K$
(respectively $H$).  Let $W_0=1$ and for $k$ at least $1$, let
$W_k=M_1\ldots M_k$. Under the assumption of the Proposition, $R$ has
the same distribution as $\sum_{i\geq 1} Q_iW_{i-1}$. Therefore, we
will derive the proposition from Theorem \RandomWeights. So, we need
to check that {\RWCondA} and {\RWCondB} hold. Referring to {\RWCondA},
we set $p=j-k$.  Clearly $p+k$ is at most $m$. Moreover, define
$$
  \rho={1\over 2}\wedge {\alpha\over \alpha+\omega} \, ,
  \qquad\hbox{ and }\qquad
  q=\Big\lfloor {\alpha+m\over\gamma\rho}+1\Bigr\rfloor \, .
$$
Furthermore, set $s=\gamma\rho q$.
Note that $t^{-s}=o\bigl(t^{-m}\oF(t)\bigr)$ at infinity, for $s$ is more than
$m+\alpha$.

Applying H\"older's inequality, we bound the expectation involved in
{\RWCondA} by
$$
  ( E|W|_1^{2p} |W|_k^{2k} )^{1/2} 
  \bigl( E N_{\alpha,\gamma,\omega}(W)^{2s}\bigr)^{1/2} t^{-s} \, .
  \eqno{\equa{EqPropIRA}}
$$
Now,
\hfuzz=26pt
$$\displaylines{\quad
  |W|_1^{2p}|W|_k^{2k}
  =(2p)!2!\sum_{n_1<\cdots < n_{2p+2}} W_{n_1}\cdots W_{n_{2p}} W_{n_{2p+1}}^k
  W_{n_{2p+2}}^k\qquad\equa{EqPropIRB} 
  \hfill\cr\hfill
  {}+\hbox{other terms where $n_i=n_j$ for some distinct $i$ and $j$.}
  \cr}
$$
\hfuzz=0pt
In all the terms in the right hand side of {\EqPropIRB}, the $M_i$'s are 
raised to a power at most $2(p+k)$, which is at most $2m$. Concerning the 
first sum on the right hand side, 
$$\eqalignno{
  &\hskip -.4in E W_{n_1}\cdots W_{n_{2p}} W_{n_{2p+1}}^k W_{n_{2p+2}}^k \cr
  {}={}& E W_{n_1}^{2p+2k} \Bigl({W_{n_2}\over W_{n_1}}\Bigr)^{2p-1+2k}
         \cdots \Bigl( {W_{n_{2p+1}}\over W_{n_{2p}}}\Bigr)^{2k}
         \Bigl( {W_{n_{2p+2}}\over W_{n_{2p+1}}} \Bigr)^k \cr
  {}={} & \mu_{H,2p+2k}^{n_1} \mu_{H,2p-1+2k}^{n_2-n_1} \cdots
          \mu_{H,2k}^{n_{2p+1}-n_{2p}} \mu_{H,k}^{n_{2p+2}-n_{2p+1}}
        &\equa{EqPropIRC}\cr
  }
$$
By Lyapounov's inequality, for $j$ at most $2m$, the inequality $\mu_{H,j}
\leq \mu_{H,2m}^{j/2m}$ holds. Therefore, {\EqPropIRC} is at most $\mu_{H,2m}$
at the power
$$\displaylines{\qquad
  {1\over 2m} 
  \bigl( (2p+2k)n_1 + (2p-1+2k)(n_2-n_1)+\cdots
  \hfill\cr\hfill
  {}+ 2k(n_{2p+1}-n_{2p})+ k(n_{2p+2}-n_{2p+1})\bigr) \, .
  \quad\cr}
$$
This exponent is at least
$$
  {1\over 2m} \bigl( kn_1+k(n_2-n_1)+\cdots + k(n_{2p+2}-n_{2p+1})\bigr)
  = {kn_{2p+2}\over 2m} \, .
$$
Our moment assumption and Lyapounov's inequality imply that $\mu_{H,2m}$ is 
less than $1$. Consequently, {\EqPropIRC} is at most 
$\mu_{H,2m}^{k n_{2p+2}/2m}$. For $n_{2p+2}$ fixed, there are at most 
$n_{2p+2}^{2p+1}$ integers $n_1,\ldots , n_{2p+1}$ less than $n_{2p+2}$.
Therefore, the first sum in the right hand side of {\EqPropIRB} is at most
$$
  \sum_{n\geq 1} n^{2p+1}(\mu_{H,m}^{k/2m})^n \, ,
$$
which is finite.

The other sums in the right hand side of {\EqPropIRB} are similarly shown to
be finite. 

So, to check {\RWCondA} it remains to prove that 
$N_{\alpha,\gamma,\omega}(W)^{2s}$ has finite expectation. It suffices to show
that $E|W|_{\gamma\rho}^{2s}$ and $E|W|_\infty^{2s}$ are finite. 

We have
$$
  E|W|_{\gamma\rho}^{2s} 
  = E\Bigl( \sum_{i\geq 1}W_i^{\gamma\rho}\Bigr)^{2q}
  = \sum_{i_1,\ldots , i_{2q}} W_{i_1}^{\gamma\rho}\cdots 
  W_{i_{2q}}^{\gamma\rho}
  \, .
$$
Again, this sum involves $M_i$'s raised at power at most $2q\gamma\rho$,
that is $2s$, which is less than $2(\alpha+m+\gamma\rho)$, which is less
than $2(\alpha+m+1)$. By assumption,
$\mu_{H,2(\alpha+m+1)}$ is less than $1$, and the same argument as before
shows that $E|W|_{\gamma\rho}^{2s}$ is finite.

Finally, note that
$$\eqalign{
  \P{ |W|_\infty^{2s}\geq t }
  &{}\leq \sum_{i\geq 1} \P{ W_i\geq t^{1/2s} } \cr
  &{}\leq \sum_{i\geq 1} t^{-(\alpha+m+1)/s} EW_i^{2(\alpha+m+1)}\cr
  &{}\leq t^{-(\alpha+m+1)/s} (1-\mu_{H,2(\alpha+m+1)})^{-1} \, .\cr
  }
$$
Since $(\alpha+m+1)/s$ is more than $1$, this shows that $|W|_\infty^{2s}$ has
finite expectation.

Finally, it is now straigtfoward to check that {\RWCondB} holds.\hfill$\qed$

\bigskip

Let us now give an explicit example. Its purpose is to show that it is now
easy to obtain the tail expansion, at least when the assumptions of
Proposition {\ImplicitRenewalTheorem} are satisfied. So for our example, take
$H$ to be the exponential distribution function of mean $\theta$,
$$
  H(t) = 1-e^{-t/\theta} \, , \qquad t\geq 0 \, ,
$$
and $K$ to be the Pareto distribution
$$
  K(t) = 1-(1+t)^{-\alpha} \, , \qquad t\geq 0 \, .
$$

To check the assumptions of Proposition \ImplicitRenewalTheorem,
we calculate
$$
  EM^\lambda 
  = \int x^\lambda \theta^{-1}e^{-x/\theta} \d x
  = \theta^\lambda \Gamma (1+\lambda ) \, .
$$
Therefore, the condition of Proposition {\ImplicitRenewalTheorem} is
simply
$$
  \theta \leq \Gamma\bigl( 2(\alpha+m)+3)\bigr)^{-1/ 2(\alpha+m+1)}
  \, .
$$

It is natural to consider the asymptotic scale $e_i(t)=t^{-\alpha-i}$,
with $i$ positive and less than $\alpha$. We will derive only two 
terms, again not because
of the difficulty of getting more, but for the space of writing the
coefficients. So we assume $\alpha$ larger that $1$ and derive a two terms
expansion for the solution of \ImplicitRenewalA.

In the chosen scale, $p_{\oK}=(1~ {-\alpha})^{\rm t}$.
As we have seen in section 3.2, the matrix
representing the derivative is defined by 
$\calD e_i = -(\alpha+i)e_{i+1}$, that is
$$
  \calD=\pmatrix{ 0 & 0 \cr -\alpha & 0 \cr} \, .
$$
The matrix representing the multiplication is
$$
  \calM_c = c^\alpha \pmatrix{ 1 & 0 \cr 0 & c \cr} \, .
$$

The first moment of the Pareto distribution is $\mu_{K,1}=1/(\alpha-1)$
while that of the exponential distribution is $\mu_{H,1} = \theta$.
So the moment equation is
$$
  \mu_{F,1} 
  = \mu_{K,1}+ \mu_{H,1}\mu_{F,1} 
  = (\alpha-1)^{-1} + \theta\mu_{F,1} \, .
$$
Therefore,
$$
  \mu_{F,1} = 1/(1-\theta)(\alpha-1) \, .
$$
We then evaluate the matrices involved in \ExpansionRenewal. So,
$$
  \calL_K = \pmatrix{ 1 & 0 \cr \alpha/(\alpha-1) & 1 \cr }
$$
and 
$$\eqalign{
  \int \calM_x \d H(x) 
  &{}=\int\diag(x^\alpha,x^{\alpha+1}) \theta^{-1} e^{-x/\theta} \d x
    \cr
  &{}= \theta^\alpha\Gamma(\alpha+1)\pmatrix{1 & 0 \cr 0 & \theta(\alpha+1)\cr}
   \cr
  }
$$
After some elementary calculation,
$$
  \calL_{H\convM F}
  =\pmatrix{1 & 0 \cr\noalign{\vskip 2pt} 
            \theta\alpha/(1-\theta)(\alpha-1) & 1\cr}
$$
Now, applying formula \ExpansionRenewal, we proved that $\oF$ has an 
asymptotic expansion $\oF\sim p_{\oF,0}e_0+p_{\oF,1}e_1$ with
$$
  p_{\oF}=
\pmatrix{
  {\ds 1\over\ds 1-\theta^\alpha\Gamma(\alpha+1)} \cr
  \noalign{\vskip 3pt}
  \alpha {\ds\theta^\alpha \Gamma(\alpha+1)(\theta-\alpha+\theta\alpha)
  +\alpha-1-\theta\alpha\over\ds
  (1-\alpha)(1-\theta)(1-\theta^\alpha\Gamma(\alpha+1))
  (1-\theta^{\alpha+1}\Gamma(\alpha+2)) } \cr}
  \, .
$$

\bigskip

To conclude this section, we mention that there are other parts of mathematics
where regular variation plays an important
role is differential equations; see for instance Mari\'c and Tomi\'c (1977),
Omey (1981) and Mari\'c (2000). Since the asymptotically smooth functions
of fixed order can be differentiated a certain number of times, they
form a natural class to use in differential equations. The general 
philosophy of this paper could be applied to some of these equations. Moreover,
because of Theorem {\MainTheorem} it is possible to obtain asymptotic expansion
for some integro-differential equations involving convolutions, Mellin 
transforms, multiplication by a function, nonlinearity due to taking powers 
of the unknown function and similar features for which we can find a stable
asymptotic scale.

\bigskip


\section{Proof in the positive case.}

In this section, we prove Theorem {\MainTheorem} under some extra assumptions.
Throughout this section, we suppose that $F^{(k)}$ is bounded and Lebesgue
integrable over the positive half line. More importantly, we also suppose that
both the $c_i$'s and the 
$X_i$'s are nonnegative. Since $\ZZ$ and $\NN^*$ are in bijection,
there is no loss of generality to index our sequence $(c_i)$ by
$\NN^*$, which we will do throughout the proof.
In the first subsection, we derive an expression
for convolutions which will be suitable for our analysis. In the
second subsection, we outline the proof. The third subsection
contains some basic facts about regularly varying functions and tail 
estimates.
Subsection 4 contains a key estimate. Then we prove some basic
lemmas (subsection 5) and finally conclude the proof by induction.

There will be many results during the proof which state that
some term $A$ say is at most some term $B$. In some cases, nothing
prevents a priori $B$ to be infinite. However, at the end, we will
see that all the upper bounds are indeed finite and prove the
theorem.

During the proof, $M$ denotes a generic constant which may 
change from place to place. This constant `depends' only on $F$ and $\omega$.
Also several of our lemmas below have conclusions which hold for
all sufficiently large reals. The bounds above which these
conclusions hold are denoted by $t_i$. These $t_i$'s only depend on $F$ 
and $\omega$.

Since we consider nonnegative $c_i$'s thoughout this section, we
assume without any loss of generality that the sequence $(c_n)_{n\geq
1}$ is nonincreasing. It is also understood that except in subsection
5.1, all distributions in the current section are supported by the
nonnegative half line.  Recall the notation $C_s$ for the sum of the
$c_i^s$, that is, in this section, $C_s=\sum_{i\geq 1} c_i^s$.

Note also that since $\omega$ is at least $1$, the function $\oF$ is normalized
regularly varying.

\bigskip

\subsection{Decomposition of the convolution into integral and
multiplication operators.}%
Let $K$ be a distribution function and $h$ be a function integrable with 
respect to the measure $\d K$. For any $\eta$ between $0$ and $1$, define 
the operator
$$
  \T_{K,\eta}h(t) = \int_{-\infty}^{\eta t} h(t-x)\d K(x)\, .
$$
Recall that for any positive real number $c$, we defined the multiplication
operator on functions
$$
  \M_c h(t)=h(t/c) \, .
$$
It is then natural to define, when possible,
$$
  \M_0 h(t)=\cases{\lim_{s\to+\infty} h(s) & if $t\geq 0$, \cr
                 \lim_{s\to-\infty} h(s) & if $t<0$.\cr }
$$
Note that if a random variable $X$ has distribution function $F$, then
$cX$ has distribution function $\M_cF$.  Moreover, for any $c$
nonnegative, $\overline{\M_cF}=\M_c\oF$.

These operators allow us to write convolutions and their derivative in
a way that will be suitable for our analysis.

\Proposition%
{\label{ConvolInOperators}%
  Let $F$ and $G$ be two distribution functions. For any $\eta$ between $0$ 
  and $1$,
  $$
    \overline{F\star G}
    =\T_{G,1-\eta}\oF + \T_{F,\eta}\oG + \M_{1/\eta}\oF \M_{1/(1-\eta)}\oG \, .
  $$
  Let $k$ be a positive integer. If $F$ and $G$ are $k$ times differentiable,
  \hfuzz=4pt
  $$
    \overline{F\star G}^{(k)}
    = \T_{G,1-\eta}\oF^{(k)} + \T_{F,\eta}\oG^{(k)}
    -\sum_{1\leq i\leq k-1}\M_{1/\eta}\oF^{(i)}\M_{1/(1-\eta)}\oG^{(k-i)} \, .
  $$
  \hfuzz=0pt
}

When $k=1$, the sum $\sum_{1\leq i\leq 0}$ in the proposition must be read
as $0$. Equivalently,
$$
  (F\star G)'=\T_{G,1-\eta}F' + \T_{F,\eta}G' \, .
$$

The proof of the proposition will be based on a lemma describing the
behavior of the operator $\T_{K,\eta}$ and $\M_c$ under differentiation.

\Lemma
{\label{LemmaConvolInOperators}
  If $K$ is a differentiable distribution function and if $h$ is a 
  differentiable function, then
  $$
    (\T_{K,\eta}h)'=\T_{K,\eta}h'+\eta \M_{1/(1-\eta)}h \M_{1/\eta}K'
  $$
  and
  $$
    (\M_ch)'=c^{-1}\M_ch' \, .
  $$
}

\Proof If $K$ is differentiable, then
$$
  \T_{K,\eta}h(t)=\int_{-\infty}^{\eta t} h(t-x)K'(x)\d x \, .
$$
The chain rule yields
$$
  {\d\over\d t} \T_{K,\eta}h(t)
  = \int_{-\infty}^{\eta t}h'(t-x)K'(x)\d x 
  +\eta h\bigl(t(1-\eta)\bigr)K'(\eta t)
$$
which is the first statement. The second one is immediate.\hfill$\qed$

\bigskip

\noindent{\bf Proof} (of the Proposition). The tail of the convolution $F*G$ 
is
$$\eqalignno{
  \overline{F\star G}(t)
  &{}=\int_{-\infty}^\infty \oF (t-y)\d G(y) 
  &\equa{ConvolInOperatorsEqA}\cr
  &{}= \int_{-\infty}^{t(1-\eta)}\oF(t-y) \d G(y) 
    +\int_{t(1-\eta)}^\infty\oF(t-y)\d G(y)
    \, .\cr
  }
$$
We integrate by parts and make a change of variable to obtain
$$\displaylines{\qquad
  \int_{t(1-\eta)}^\infty\oF(t-y)\d G(y)
  \hfill\cr\hfill
  \eqalign{
    &{}= \bigl[ \oF(t-y)\bigl( G(y)-1\bigr)\bigr]_{t(1-\eta)}^\infty
     +\int_{t(1-\eta)}^\infty\oG(y)\d \oF(t-y) \cr
    &{}=\oF(t\eta)\oG\bigl(t(1-\eta)\bigr) 
     + \int_{-\infty}^{t\eta}\oG(t-y)\d F(y) \, . \cr
    }
  \qquad\cr
  }
$$
Combined with \preveq , we obtain the first assertion of the proposition.

To prove the second assertion, we proceed by induction, starting to prove the
result for $k=1$. Using the lemma and differentiating the expression for
$\overline{F\star G}$, we see that
$$\eqalign{
  \overline{F\star G}{}'
  &{}=\T_{G,1-\eta}\oF' + \T_{F,\eta}\oG' 
   + (1-\eta)\M_{1/\eta}\oF \M_{1/(1-\eta)}G'
   \cr
  \noalign{\vskip2pt}
  &\qquad {}+\eta \M_{1/(1-\eta)}\oG \M_{1/\eta}F' 
   + \eta \M_{1/\eta}\oF' \M_{1/(1-\eta)}\oG \cr
  \noalign{\vskip 2pt}
  &\qquad {}+(1-\eta) \M_{1/\eta}\oF \M_{1/(1-\eta)}\oG' \, .\cr
  }
$$
Since $\oF'=-F'$, we obtain $\overline{F\star G}{}'=\T_{G,1-\eta}\oF'
+\T_{F,\eta}\oG'$.

Assume now that the relation holds for any integer up to $k-1$. Using the
lemma,
$$\eqalign{
  (\overline{F\star G}^{(k-1)})'
    &{}=\T_{G,1-\eta}\oF^{(k)} + \T_{F,\eta}\oG^{(k)} \cr
    &\qquad {}+ (1-\eta) \M_{1/\eta}\oF^{(k-1)} \M_{1/(1-\eta)}G' \cr
    &\qquad {}+\eta \M_{1/(1-\eta)} \oG^{(k-1)}\M_{1/\eta}F' \cr
    &\qquad {}-\sum_{1\leq i\leq k-2}\eta \M_{1/\eta}\oF^{(i+1)}
     \M_{1/(1-\eta)}\oG^{(k-1-i)} \cr
    &\qquad {}-\sum_{1\leq i\leq k-2} (1-\eta) \M_{1/\eta}\oF^{(i)}
     \M_{1/(1-\eta)}\oG^{(k-i)} \, ,\cr
  }
$$
which after collecting the terms give the proper 
expression for $\overline{F\star G}^{(k)}$.\hfill$\qed$

\bigskip

\subsection{Organizing the proof.}%
The proof goes essentially by induction on $k$, and for fixed $k$ by induction 
on $n$. We write $G_n$ the distribution funciton 
of $\sum_{1\leq i\leq n} c_iX_i$ and
$G$ that of $\< c,X\>$. We also write $F_i$ for $M_{c_i}F$, that is for the
distribution function of $c_iX_i$. Hence, $G_n=F_1\star \cdots \star F_n$.

We also define a paramater $\rho$ in $(0,1)$, which will be chosen
at the end of the proof. We set $d_n=c_n^\rho$.

Applying Proposition \ConvolInOperators, we see that
$$\displaylines{\qquad
  \oG_n^{(k)} = \T_{G_{n-1},1-d_n}\oF_n^{(k)} + \T_{F_n,d_n}\oG_{n-1}^{(k)}
  \hfill\cr\noalign{\vskip 3pt}\hfill
  {}- \sum_{1\leq i\leq k-1} \M_{1/d_n}\oF_n^{(i)} \M_{1/(1-d_n)} 
  \oG_{n-1}^{(k-i)}
  \, .
  \qquad\cr}
$$
Let us write $\AmoGnk$ the $m$ terms asymptotic expansion
of $\oGnk$ --- or, more precisely, for the time being, the candidate
for this asymptotic expansion ---
$$
  \AmoGnk=\sum_{1\leq i\leq n} \L_{G_n\natural F_i,m}\oFik \, .
$$
The `$\calA_m$' in this formula is not an operator; it is aimed as a short hand
for `$m$-terms approximation of' and merely to help the memory. Sometimes
we will omit the subscript $m$, writing simply $\calA$.

We will approximate $\T_{G_{n-1},1-d_n}$ by $\L_{G_{n-1}}$ 
(Theorem \fixedref{5.4.1})
and $\T_{F_n,d_n}$ by $\L_{F_n}$. Thus, we expect an approximation of $\oGnk$
by $\L_{G_{n-1}}\oFnk+\L_{F_n}\oG_{n-1}^{(k)}$. Since $\oG_{n-1}^{(k)}$
is close to its asymptotic expansion, we also hope to approximate 
$\L_{F_n}\oG_{n-1}^{(k)}$ by $\L_{F_n}\calA \oG_{n-1}^{(k)}$. But
$$
  \L_{F_n}\calA \oG_{n-1}^{(k)}
  = \sum_{1\leq i\leq n-1} \L_{F_n} \L_{G_{n-1}\natural F_i} \oFik \, .
$$
Now, $\L_{F_n}\L_{G_{n-1}\natural F_i}$ is not 
$\L_{F_n\star (G_{n-1}\natural F_i)}$ (the latter being $\L_{G_n\natural F_i}$); 
however, these two operators are equal up to differential operators 
involving $\D^{m+1+i}$'s for $i$ nonnegative. Since $\oF_i^{(k+l)}\asymp \Id^{-l}\oFik$ 
when $\oF_i^{(k+l)}$ is regularly varying, we expect that
$$
  \L_{F_n}\calA \oG_{n-1}^{(k)} 
  \sim \sum_{1\leq i\leq n-1} \L_{G_n\natural F_i}\oFik
  \, .
$$
Hence, we should obtain 
$$
  \oGnk\sim \L_{G_{n-1}}\oFnk 
  +\sum_{1\leq i\leq n-1}\L_{G_n\natural F_i}\oFik \, , 
$$
which is equal to $\calA\oGnk$. Thus, assuming that $\oG_{n-1}^{(k)}\sim \calA
\oG_{n-1}^{(k)}$, we can expect to prove $\oG_n^{(k)}\sim \calA\oGnk$, and
build this way a proof by induction on $n$. In any case, this 
derivation shows that
$$\eqalignno{%
  \hskip -15pt\oGnk-\calA\oGnk
  &{}=(\T_{G_{n-1},1-d_n}-\L_{G_{n-1}})\oFnk \cr
  \noalign{\vskip 2pt}
  &\quad {}+ \T_{F_n,d_n}(\oG_{n-1}^{(k)}-\calA\oG_{n-1}^{(k)}) \cr
  \noalign{\vskip 3pt}
  &\quad {}+\Bigl(\T_{F_n,d_n}\calA \oG_{n-1}^{(k)}
   -\sum_{1\leq i\leq n-1} \L_{G_n\natural F_i}\oFik \Bigr) \cr
  &\quad {}-\sum_{1\leq j\leq k-1} \M_{1/d_n}\oF_n^{(j)}
   \M_{1/1-d_n}\oG_{n-1}^{(k-j)} \, .~
  &\equa{InductionEquality}\cr
  }
$$

In the right hand side of this equality, we control $\oFnk$ by assumptions
on $\oF$ and the difference $\T_{G_{n-1},1-d_n}-\L_{G_{n-1}}$ by Theorem 
\fixedref{5.4.1}.
Since $\calA\oG_{n-1}^{(k)}$ is expressed in terms of functions $\oF_i$'s,
hence, ultimately in terms of $\oF$, it is an explicit term on which
we can work. The only term on which we do not have a rather direct
control is $\T_{F_n,d_n}(\oG_{n-1}^{(k)}-\calA\oG_{n-1}^{(k)})$. Up
to a shift of index $n$, it is related to the left hand side of \preveq ,
provided we can control the operator norm of $\T_{F_n,d_n}$. This will
allow us to prove the asymptotic expansion
$\oGnk\sim\AmoGnk$ is valid uniformly in $n$.

For the induction and the estimates to be written in a manageable form,
we use a pseudo-semi-norm to control the remainder terms. For any
function $h$ on the nonnegative half line, we define
$$
  |h|_{m,\tau}=\sup_{t\geq\tau} t^m |h(t)|/\oF(t) \, .
$$
The choice of this pseudo-semi-norm is motivated by the fact that
$\lim_{\tau\to\infty} |h|_{m,\tau}=0$ is equivalent to 
$h(t)=o\big(t^{-m}\oF(t)\big)$ as $t$ tends to infinity. Therefore, 
proving an asymptotic equivalence up to a $t^{-m}\oF(t)$ term 
amounts to prove that the pseudo-semi-norm $|\cdot|_{m,\tau}$ of some 
function tends to $0$ as $\tau$ tends to infinity.

\bigskip

\subsection{Regular variation and basic tail estimates.}%
In this subsection, we collect some facts on regular variation
and on the distribution of $\< c,X\>$.
Our first lemma sandwiches the tail of $\cX$.

\Lemma{
  \label{TrivialTailBound}
  If $C_\rho\leq 1$, then
  $$
    \M_{c_1}\oF\leq \oG_n \leq \oG \leq \sumn \M_{c_n^{1-\rho}}\oF \, .
  $$ 
  Moreover, let $\epsilon$ be a positive real number. There exists $t_3$
  such that for any $t$ larger than $t_3$, 
  $$
    \oG(t)\leq C_{\alpha(1-\rho)-\epsilon}\oF(t) \, .
  $$
}

\Proof The first two bounds follow from the inequalities
$$
  c_1X_1 \leq c_1X_1+\cdots + c_nX_n\leq \cX \, .
$$
Next, if $c_n^{1-\rho}X_n\leq t$ for all $n$, then
$$
  \sumn c_n X_n \leq t\sumn c_n^\rho = t C_\rho \leq t \, .
$$
Therefore,
$$
  \P{\cX > t} 
 \leq P\bigl\{\,\cup_{n\geq 1}\{\, c_n^{1-\rho}X_n > t\,\}\,\bigr\} 
  = \sumn \P{ c_n^{1-\rho} X_n > t} \, , 
$$
which is the last inequality in the first statement. This
inequality asserts that $\oG(t)$ is
at most $\sum_{n\geq 1} \oF(t/c_n^{1-\rho})$. Since $C_\rho$ is at most
$1$, so are all the $c_i$'s.  Apply Lemma {\Potter} with $g=\oF$, 
$\delta=\epsilon/(1-\rho)$ and $\lambda=c_n^{\rho-1}$ to obtain
the second statement.
\hfill$\qed$

\bigskip

This lemma will be useful when combined with the next analytical result.

\Lemma
{
  \label{StochasticOrderBound}%
  Let $g$, $h$ be two nonincreasing functions of bounded variations,
  such that $g\leq h$ and 
  $$
    \limt g(t) =\limt h(t) =0 \, .
  $$
  Let $f$ be a nonnegative nondecreasing and right continuous function 
  on $\RR^+$. Then
  $$
    -\int f \d g\leq -\int f \d h \, .
  $$
}

\Proof Define the left continuous inverse of $f$ by
$$
  f^\leftarrow (u)=\inf\{\, x \, :\, f(x)\geq u\,\}\, .
$$
Observe that $f(x)\geq u$ if and only if $f^\leftarrow (u)\leq x$. 
We use the Fubini-Tonnelli theorem to write
$$\eqalign{
  \int_0^\infty f\d g 
  &{}=\int_0^\infty \int_0^\infty \II\{\,u\leq f(x)\,\} \d u\, \d g(x) \cr
  &{}=\int_0^\infty \int_0^\infty \II\{\, x\geq f^\leftarrow (u)\,\} \d g(x) 
   \d u \cr
  &{}=-\int_0^\infty g\circ f^\leftarrow (u) \d u \, .\cr
  }
$$
The result follows from the inequality $g\circ f^\leftarrow\leq h\circ 
f^\leftarrow$.\hfill$\qed$

\bigskip

Since $\oG$, $\oGn$, $\oF_i$ and $\oF$ are of the same order, so are their 
truncated means $\int_t^\infty x^k \d K(x)$ for $K=\oG$, $\oGn$, $\oF_i$ 
and $\oF$. Order is related to asymptotic behavior, and because we 
will need to have estimates that are uniform
over $n$, we need a quantitative statement comparing truncated means. The
following will do for our purpose and may be useful in other contexts.

\Lemma 
{
  \label{TruncatedMomentBoundF}%
  Let $\oF$ be a normalized regularly varying function of index $-\alpha$. Let
  $\delta$ be a positive number. There exists $t_4$ such that for 
  any $0<u\leq v$, any $t\geq t_4$, any integer $k<\alpha$,
  $$
    \int_{tv}^\infty x^k \d \M_uF(x) 
    \leq {2\alpha\over\alpha-k} v^k \Bigl({u\over v}\Bigr)^{\alpha-\delta} 
    t^k \oF(t) \, .
  $$
}

\Proof A change of variable yields
$$
  \int_{tv}^\infty x^k \d \M_uF(x) 
  = u^k\int_{tv/u}^\infty x^k \d F(x) \, .
$$
By Karamata's theorem for Stieltjes integrals (Bingham, Goldie, Teugels, 
1989, Theorem 1.6.5), there exists a $t_4'$ such that for $t$ more than $t_4'$
$$
  \int_t^\infty x^k\d F(x)
  \leq {2\alpha\over \alpha-k} t^k\oF (t) \, .
$$
This $t_4'$ may be taken independent of $k$ since there is only a finite number
of integers less than $\alpha$. Since $v/u\geq 1$, this implies that for
$t$ at least $t_4'$,
$$
  \int_{tv}^\infty x^k \d \M_uF(x)
  \leq {2\alpha\over \alpha-k} v^k t^k \oF(tv/u) \, .
$$
By Lemma \Potter, if $t\geq t_2$, we have 
$\oF (tv/u)\leq (u/v)^{\alpha-\delta}\oF (t)$.
Take $t_4=t_2\vee t_4'$ to conclude.\hfill$\qed$

\Lemma
{ 
  \label{TruncatedMomentBoundGn}%
  Let $\oF$ be a normalized regularly varying function
  with index $-\alpha$. Let $\epsilon$ be a positive number. There exists
  $t_5$ such that for any $t\geq t_5$, any $n\geq 1$, any sequence
  $(c_n)_{n\geq 1}$ with $C_\rho\leq 1$, and any
  $k<\alpha$,
  $$
    \int_t^\infty x^k \d G_n(x)
    \leq \int_t^\infty x^k \d G(x)
    \leq {2\alpha\over \alpha-k}\, C_{\alpha(1-\rho)-\epsilon} t^k\oF (t) 
    \, .
  $$
}

\Proof Combining Lemmas {\TrivialTailBound} and \StochasticOrderBound, we have
$$
  \int_t^\infty x^k\d G_n(x)
  \leq \int_t^\infty x^k \d G(x)
  \leq \sumn \int_t^\infty x^k \d \M_{c_n^{1-\rho}} F(x) \, ,
$$
where we used that $C_\rho\leq 1$ in the last inequality. Since $C_\rho$ is 
at most $1$, so are all the $c_n$'s. We apply Lemma 
{\TruncatedMomentBoundF} with $\delta=\epsilon/(1-\rho)$ to obtain
$$
  \int_t^\infty x^k \d \M_{c_n^{1-\rho}} F(x)
  \leq {2\alpha\over \alpha-k} c_n^{(1-\rho)(\alpha-\delta)} t^k \oF (t)
$$
for $t$ at least $t_4$\hfill$\qed$

\bigskip

We will also use the following bound on the moment of $G$.

\Lemma
{\label{MomentBoundG}%
  Assume $C_\rho$ is at most $1$. For any $s$ at least $1$,
  $$
    \mu_{G,s}\leq C_{s(1-\rho)}\mu_{F,s} \, . 
  $$
}

\Proof
Write
$$
  \mu_{G,s}=s\int_0^\infty x^{s-1}\oG(x) \d x \, .
$$
Apply Lemma {\TrivialTailBound} to obtain
that the moment of order $s$ of $G$ is at most
$$\eqalignno{
  \sum_{n\geq 1} \int_0^\infty s x^{s-1} \M_{c_n^{1-\rho}}\oF(x) \d x
  &{}= \sum_{n\geq 1} \int_0^\infty (c_n^{1-\rho})^s x^s \d F(x) \cr
  &{}= C_{s(1-\rho)}\mu_{F,s} \, .
  &\qed\cr
  }
$$
 
\bigskip

\subsection{The fundamental estimate.}%
The origin of our asymptotic expansions is the following estimate,
whose proof is almost deceptively simple.

\Theorem
{\label{ApproxTByL}%
  Let $m$ be a positive integer and let $r$ be in $[\, 0,1)$. 
  Furthermore, let $K$ be a distribution function on the nonnegative 
  half line, whose $m$-th
  moment is finite. If $h$ is smoothly varying of order $m+r$, then
  $$\displaylines{\quad
    |(\T_{K,\eta}-\L_{K,m})h|(t)
    \leq \sum_{0\leq j\leq m} {|h^{(j)}(t)|\over j!}
    \int_{\eta t}^\infty x^j \d K(x) 
  \hfill\cr\hfill
    {}+{|h^{(m)}(t)|\over t^r m!} 
    \int_0^{\eta t}\oDelta_{t,x/t}^r(h^{(m)})x^{m+r}\d K(x) \, .
  \quad\cr}
  $$
}

\Proof The expression of $\T_{K,\eta}$ and Proposition {\Taylor} show that
$$\displaylines{\qquad
    \Bigl| \T_{K,\eta}h(t)-\int_0^{\eta t}\sum_{0\leq j\leq m} (-1)^j 
    {x^j\over j!} h^{(j)}(t)\d K(x)\Bigr|
  \hfill\cr\hfill
    {}\leq {|h^{(m)}(t)|\over t^r m!} \int_0^{\eta t}
    \oDelta_{t,x/t}^r(h^{(m)}) 
    x^{m+r} \d K(x) \, .
  \qquad\global\advance\equanumber by 1\preveq
  \cr}
$$
Since
$$
  \int_0^{\eta t} (-1)^j {x^j\over j!} h^{(j)}(t) \d K(x)
  = {(-1)^j\over j!} \Bigl( \mu_{K,j}-\int_{\eta t}^\infty x^j \d K(x)\Bigr)
  \D^j h(t) \, ,
$$
the result follows from the triangle inequality.\hfill$\qed$

\bigskip

We could also bound {\preveq} by the simpler estimate
$$
  t^{-r}|h^{(m)}(t)|\oDelta_{t,\eta}^r(h^{(m)})\mu_{K,m+r} \, .
$$
Unfortunately, this is not good enough for our purpose.

Let us now explain why the estimate in the previous Theorem is useful.
When $h$ is smoothly varying and has index $-\beta$, the term of 
smallest order in $\L_{K,m}h$ is given by $h^{(m)}$, which is in 
$RV_{-\beta-m}$. If $\oK$ is regularly varying with index $-\alpha$,
Karamata's theorem for Stieltjes integrals (Bingham, Goldie and
Teugels, 1989, Theorem 1.6.5) shows that the function $t\mapsto
\int_{\eta t}^\infty x^i \d K(x)$ is in $RV_{-\alpha+i}$.
Thus, $h^{(i)}(t)\int_{\eta t}^\infty x^i\d K(x)$ is regularly varying
with index $-\alpha-\beta$; this index is smaller than $-\beta-m$
if $m<\alpha$. Thus, Theorem {\ApproxTByL}
proves that when $h$ is smoothly varying, $\T_{K,\eta}h\sim \L_{K,m} h$, 
and, moreover, gives
an estimate of the error term of this asymptotic expansion.

\bigskip

\subsection{Basic lemmas.}%
The goal of this subsection is to give estimates for the first three
terms of the right hand side of \InductionEquality.

\hfuzz=2pt
Our first lemma takes care of the term 
$(\T_{G_{n-1},1-d_n}-\L_{G_{n-1},m})\oF_n^{(k)}$. It asserts that
it is of smaller order than $t^{-m-k}\oF(t)$, uniformly in $n$ and
in some sequences $(c_n)_{n\geq 1}$. We write $v(t)$ a function
such that $0\leq v(t)\leq t$ for all nonnegative $t$, and $\limt v(t)
=+\infty$ while $\limt v(t)/t=0$. We assume moreover that $v(t)/t$ is
ultimately nonincreasing. For instance, $v(t)=t^\kappa$ with $0<\kappa <1$
will do.

\hfuzz=0pt

\Lemma
{\label{BasicLemmaA}%
  Let $\epsilon$ be a positive number.
  Assume that $\oF$ is in 
  $SR_{-\alpha,m+r+k}$ with $m$ an integer, $r$ in $[\,0,1)$ and $m+r$
  smaller than $\alpha$. 
  There exist $M$ and $t_6$ depending only on $F$, $m$, $k$ and $\epsilon$
  such that whenever $C_\rho\vee C_{(1-\rho)(m+r)}\leq 1$,
  and whenever $d_n\leq 1/2$, for any $t$ at least $t_6$,
  $$\displaylines{\quad
    |(\T_{G_{n-1},1-d_n}-\L_{G_{n-1},m})\oF_n^{(k)}|_{k+m,t}
    \hfill\cr\noalign{\vskip 3pt}\hfill
    {}\leq M c_n^{\alpha-(\alpha+k+m)\rho-\epsilon} 
    \Bigl( \oDelta_{t,v(t)/t}^r(\oF^{(k+m)})+v(t)^{m+r}\oF\circ v(t)\Bigr) \, .
    \quad\cr}
  $$
}

\Proof Applying Theorem \ApproxTByL, we obtain for $t$ at least $1$,
$$\displaylines{
  |(\T_{G_{n-1},1-d_n}-\L_{G_{n-1},m})\oF_n^{(k)}|(t)
  \hfill\cr\noalign{\vskip 3pt}\quad
  {}\leq{}
   \sum_{0\leq j\leq m} |\oF_n^{(k+j)}(t)| 
     \int_{t(1-d_n)}^\infty x^j\d G_{n-1}(x) 
   \hfill\cr\noalign{\vskip 3pt}\hfill
     {}+\, |\oF_n^{(k+m)}(t)| 
     \int_0^{t(1-d_n)} \oDelta^r_{t,x/t}(\oF_n^{(k+m)}) x^{m+r} 
     \d G_{n-1}(x) \, .
  \quad\global\advance\equanumber by1\preveq\cr
  }
$$
Since the $c_i$'s are all at most $1$, it suffices to prove
the result for $\epsilon$ small enough. Thus, we can assume that $m+r$ 
is at most $\alpha-\epsilon$. Let $\delta=\epsilon/(2-\rho)$. We can assume
that $\epsilon$ is small enough so that $(1-\rho)\alpha-\delta$ is at least
$(1-\rho)(m+r)$. Then,
the sum $C_{\alpha (1-\rho)-\delta}$ is at most $C_{(1-\rho)(m+r)}$,
hence at most $1$.
Combining Lemmas {\BoundDkMF} and \TruncatedMomentBoundGn, we see that for
$t(1-d_n)\geq t_1\vee t_5$,
$$\eqalign{
  |\oF_n^{(k+j)}(t)| \int_{t(1-d_n)}^\infty x^j \d G_{n-1}
  & {}\leq M c_n^{\alpha-\delta} t^{-k-j}\oF(t) t^j\oF(t) \cr
  \noalign{\vskip 3pt}
  & {}\leq M c_n^{\alpha-\delta} t^{-k} \oF(t)^2 \, . \cr
  }
$$
Applying Lemmas \BoundDkMF, \DeltaDM, {\TrivialTailBound} and 
{\StochasticOrderBound}, we also have, for $t$ at least $t_1$ and $t_2$,
$$\displaylines{\quad
    |\oF_n^{(k+m)}(t)|\int_0^{t(1-d_n)} \oDelta_{t,x/t}^r(\oF_n^{(k+m)})
    x^{m+r}\d G_{n-1}(x)
  \hfill\cr\hfill
    {}\leq M c_n^{\alpha-\delta} t^{-k-m} \oF(t) 
    \int_0^{t(1-d_n)} \oDelta_{t/c_n,x/t}^r(\oF^{(k+m)}) x^{m+r} \d G(x) \, .
  \quad\cr
  }
$$
We split the integral in this upper bound as one for $x$ between $0$ and
$v(t)$, and one between $v(t)$ and $t(1-d_n)$.  We then have the easy bound
$$
  \int_0^{v(t)}\oDelta_{t/c_n,x/t}^r(\oF^{(k+m)})x^{r+m} \d G(x)
  \leq \oDelta_{t,v(t)/t}^r (\oF^{(k+m)})\mu_{G,m+r} \, .
$$
We use Lemma {\MomentBoundG} to bound $\mu_{G,m+r}$ by
$C_{(m+r)(1-\rho)}\mu_{F,m+r}$, which is then at most $\mu_{F,m+r}$ under
the assumptions of the lemma.

To bound the integral for $x$ between $v(t)$ and $t(1-d_n)$, we first
bound $\oDelta_{t/c_n,x/t}^r(\oF^{(k+m)})$ by 
$\oDelta_{t/c_n,1-d_n}^r (\oF^{(k+m)})$. The latter is at most
\hfuzz=1pt
$$
  \oDelta_{t/c_n,1/2}^r(\oF^{(k+m)}) 
  + \sup_{s\geq t/c_n}\sup_{1/2<|y|\leq 1-d_n} \hskip -4pt y^{-r}
  \biggl( \Bigl|{\oF^{(k+m)}\bigl( s(1-y)\bigr)\over \oF^{(k+m)}(s)}\Bigr|
         + 1 \biggr) .
$$
\hfuzz=0pt
If $s$ is at least $t/c_n$ and $y$ at most $1-d_n$, then $s(1-y)$ is
at least $td_n/c_n$, and therefore at least $t$. Thus, by the standard
Potter bounds, there exists $t_6'$ such that for $s$ at least $t_6'$,
$$\eqalign{
  \Bigl| {\oF^{(k+m)}\bigl( s(1-y)\bigr)\over \oF^{(k+m)}(s)} \Bigr|
  & {}\leq 2 \bigl( (1-y)^{-\alpha-k-m-\delta}\vee 
    (1-y)^{-\alpha-k-m+\delta} \bigr) \cr
  & {}\leq 2 d_n^{-\alpha-k-m-\delta} \, .\cr
  }
$$
Consequently,
$$\displaylines{\quad
    \int_{v(t)}^{t(1-d_n)} \oDelta_{t/c_n,x/t}^r (\oF^{(k+m)}) x^{m+r} \d G(x)
  \hfill\cr\hfill
    {}\leq M \bigl( \oDelta_{t/c_n,1/2}^r(\oF^{(k+m)})
    +d_n^{-\alpha-k-m-\delta}+ 1\bigr) \int_{v(t)}^\infty x^{m+r}\d G(x) \, .
  \cr}
$$
Since $d_n$ is smaller than $1/2$, we can simplify this upper bound into
$$
  M d_n^{-\alpha-k-m-\epsilon}\int_{v(t)}^\infty x^{m+r}\d G(x) \, .
$$
Applying Lemma {\TruncatedMomentBoundGn}, for $v(t)$ at least $t_5$,
this last quantity is at most
$$
  M d_n^{-\alpha-k-m-\delta}v(t)^{m+r}\oF\circ v(t) \, .
$$
It follows that
$$\displaylines{\qquad
    \int_0^{t(1-d_n)} \oDelta_{t/c_n,x/t}^r(\oF^{(k+m)}) x^{m+r}\d G(x)
  \hfill\cr\hfill
    {}\leq M\oDelta_{t,v(t)/t}^r (\oF^{(k+m)}) + M d_n^{-\alpha-k-m-\delta}
    v(t)^{m+r} \oF\circ v(t) \, . 
  \qquad\cr}
$$
Combined with the other bounds, we obtain that the right hand side of {\preveq}
is at most
$$\displaylines{\qquad
  M c_n^{\alpha-\delta} t^{-k-m}\oF(t)\Bigl( t^m\oF(t)
  +\oDelta_{t,v(t)/t}^r(\oF^{(k+m)})
  \hfill\cr\hfill 
  {}+ d_n^{-\alpha-k-m-\delta}v(t)^{m+r}\oF\circ v(t) \Bigr) \, .
  \qquad\cr}
$$
Since $m+r$ is smaller than $\alpha$, the function $s\mapsto s^{m+r}\oF(s)$
is ultimately nonincreasing (Bingham, Goldie and Teugels, 1989, Theorem 1.5.3).
Thus, for $t$ large enough, $t^m\oF(t)$ is smaller than 
$2v(t)^{m+r}\oF\circ v(t)$.
The proof is completed in noting that $d_n^{-\alpha-k-m-\delta}$ is at least
$1$.\hfill$\qed$

\bigskip

Our next lemma shows that when $d_n$ is small, that is $n$ is large,
then $T_{F_n,d_n}$ is very close to be a contraction.

\Lemma%
{\label{TNearContraction}%
  There exists some positive numbers
  $M$ and $t_7$ such that for any $t$ at least $t_7$, any positive $n$,
  any $k$, any $d_n$ smaller than $1/2$, any function $h$,
  $$
    |\T_{F_n,d_n}h|_{k,t} \leq |h|_{k,t(1-d_n)}(1+Md_n) \, .
  $$
}

\Proof A change of variable shows that
$$
  \T_{F_n,d_n}h(t)
  =\int_0^{td_n/c_n} h(t-c_nx) \d F(x) \, .
$$
If $x$ is at most $td_n/c_n$, then $t-c_nx$ is at least $t(1-d_n)$.
Consequently,
$$
  |h(t-c_nx)|
  \leq |h|_{k,t(1-d_n)} {\oF(t-c_nx)\over (t-c_nx)^k}
  \leq |h|_{k,t(1-d_n)} {\oF\bigl(t(1-d_n)\bigr)\over t^k(1-d_n)^k} \, .
$$
This implies
$$
  {t^k |\T_{F_n,d_n}h(t)|\over \oF(t)} 
  \leq |h|_{k,t(1-d_n)} {\oF\bigl(t(1-d_n)\bigr)\over\oF (t)} (1-d_n)^{-k}
  \, . 
$$
Since $d_n$ is at most $1/2$, we can apply Lemma \Potter, and then bound
$(1-d_n)^{-\alpha-k-\delta}$ by $1+Md_n$.\hfill$\qed$

\Remark
{In proving Lemma \TNearContraction, our improved Potter bound is 
  crucial. With the usual one we could only get $A(1+Md_n)$ for some $A$ 
  greater than $1$. This is not good enough for the inductions to come. 
  Extensive attempts to prove our theorem under weaker assumptions 
  suggest that the behavior of the function $c(\cdot )$ in the Karamta 
  representation (see Bingham, Goldie and Teugels, 1989, Theorem 1.3.1) 
  of $\oF$ cannot be arbitrary for our asymptotic expansion to hold. 
  In particular, the asymptotic behavior of $c(\cdot)$ plays a key 
  role in the asymptotic expansion. This is already the case for a two
  terms expansion where the second term involves derivative; one should
  indeed remark that the absolute continuity of the function $c(\cdot)$
  is equivalent to that of $\oF$.
}

\bigskip

In particular, in connection with the second term in \InductionEquality,
we obtain the following estimate, considering 
$h=\oG_{n-1}^{(k)}-\calA_m\oG_{n-1}^{(k)}$ in Lemma \TNearContraction.

\Corollary
{\label{CorToTNearContraction}%
  There exist positive numbers
  $M$ and $t_8$ such that for any $t$ at least $t_8$ and any sequence
  $(c_n)_{n\geq 1}$ with $(d_n)_{n\geq 1}$ uniformly bounded by $1/2$,
  any $k$ and any $n$ at least $2$,
  $$\displaylines{\qquad
    |\T_{F_n,d_n} (\oG_{n-1}^{(k)}-\calA_m \oG_{n-1}^{(k)})|_{k+m,t}
    \hfill\cr\noalign{\vskip 3pt}\hfill
    {}\leq |\oG_{n-1}^{(k)}-\calA_m\oG_{n-1}^{(k)}|_{k+m,t(1-d_n)} (1+Md_n) \, .
    \qquad\cr
    }
  $$
}

This estimate shows the need to introduce the family
of operators $\T_{K,\eta}$ and not use $T_{K,1/2}$ as in Barbe and McCormick
(2005). Indeed, when doing the induction later, we see that the value
of $t$ in the norm of the right hand side drops by a factor $1-d_n$. When
applying this bound inductively, we will obtain something like
$$
  |\oG_2^{(k)}-\calA\oG_2^{(k)}|_{k+m,t\Pi_{1\leq i\leq n}(1-d_i)} \, .
$$
If $\prod_{i\geq 1}(1-d_i)=0$, we will not be able to permute the limits
as $n$ tends to infinity and $t$ tends to infinity. It is therefore 
essential that the series $\sum_{n\geq 1} d_n$ converges.

\bigskip

Our last lemma in this section will handle the third term in the right
hand side of \InductionEquality.

\Lemma%
{\label{BasicLemmaC}%
  Let $\epsilon$ be a positive number at most $\alpha-1$. 
  Assume that $F$ is in $SR_{-\alpha,m+k+r}$.
  There exist positive $M$ and $t_9$
  such that for any sequence $(c_n)_{n\geq 1}$ with $C_\rho$ at most $1$,
  any $t$ at least $t_9$ and any $n$ at least $2$,
  $$\displaylines{\qquad
      \Bigl|\T_{F_n,d_n}\calA_m\oG_{n-1}^{(k)}
      -\sum_{1\leq i\leq n-1}\L_{G_n\natural F_i} \oF_i^{(k)}\Bigr|_{k+m,t}
    \hfill\cr\hfill
      \leq M\Bigl( c_n^{\alpha (1-\rho)-\epsilon}
      +c_n^r\oDelta_{t,d_n}^r(\oF^{(k+m)}) \Bigr) v(t)^{m+r}\oF\circ v(t)
    \hfill\cr\noalign{\vskip 3pt}\qquad\qquad\qquad
      {}+ Mc_n^r \oDelta_{t,c_nv(t)/t}^r(\oF^{(k+m)}) \, .
    \hfill\cr
    }
  $$
}

The proof of this lemma requires an auxilliary result which we state
as a claim. It compares $T_{K,\eta}\L_H$ with $L_{K\star H}$ (cf.\
subsection \fixedref{5.2}).

\medskip

\Claim
{ For any nonnegative $t$,
  $$\displaylines{
      |\T_{K,\eta}\L_H h-\L_{K\star H} h|(t)
    \hfill\cr\noalign{\vskip 3pt}\hfill
    \eqalign{
      {}\leq & \sum_{0\leq s\leq m} {|h^{(s)}(t)|\over s!} 
               \sum_{0\leq j\leq s} {s\choose j}
               \mu_{H,j} \int_{\eta t}^\infty x^{s-j}\d K(x)\cr
             & {}+{|h^{(m)}(t)|\over t^r m!}
               \sum_{0\leq j\leq m} {m\choose j} \mu_{H,j} 
               \int_0^{\eta t}\oDelta_{t,x/t}^r(h^{(m)})x^{m+r-j} \d K(x)\cr
    }%
   \cr}
  $$
}

\Proof 
By linearity of $\T_{K,\eta}$ we have
$$
  \T_{K,\eta}\L_H 
  =\sum_{0\leq j\leq m} {(-1)^j\over j!} \mu_{H,j}\T_{K,\eta}\D^j \, .
$$
Applying Theorem {\ApproxTByL} to $h^{(j)}$ in the form of a bound
for $|\T_{K,\eta}h^{(j)}-\L_{K,m-j}h^{(j)}|$, we obtain
$$\displaylines{\quad
  \Bigl| T_{K,\eta}\L_{H,m}h-\sum_{0\leq j\leq m} {(-1)^j\over j!}
         \mu_{H,j}\L_{K,m-j}h^{(j)} \Bigr|(t)
  \hfill\cr\qquad
  {}\leq \sum_{0\leq j\leq m} {\mu_{H,j}\over j!} 
  \biggl(\, \sum_{0\leq l\leq m-j}
  {|h^{(j+l)}(t)|\over l!}\int_{\eta t}^\infty x^l \d K(x) 
  \hfill\cr\hfill
  {}+ {|h^{(m)}(t)|\over t^r (m-j)!}
  \int_0^{\eta t}\oDelta_{t,x/t}^r(h^{(m)})x^{m-j+r}\d K(x)
  \biggr) \, .
  \quad\cr
  }
$$
Using Lemma {\LaplaceBinomial}, the left hand side is the absolute value of 
$(\T_{K,\eta}\L_H-\L_{K\star H})h$ evaluated at $t$. Setting $s=l+j$,
the right hand side is at most
$$\displaylines{\quad
  \sum_{s,j}\II\{\, j\leq s\leq m\,\} {|h^{(s)}(t)|\over s!}
  {s\choose j} \mu_{H,j}\int_{\eta t}^\infty x^{s-j}\d K(x) 
  \hfill\cr\hfill
  {}+{|h^{(m)}(t)|\over t^r m!}\, \Bigl(
  \sum_{0\leq j\leq m} {m\choose j} \mu_{H,j}
  \int_0^{\eta t}\oDelta_{t,x/t}^r(h^{(m)}) x^{m-j+r}\d K(x) \Bigr)
  \, .\cr
  }
$$
This is our claim.\hfill$\qed$

\noindent{\bf Proof} (of Lemma \BasicLemmaC). First the triangle
inequality implies the pointwise inequality
$$\displaylines{\qquad
  \Bigl|\T_{F_n,d_n}\calA \oG_{n-1}^{(k)}
   -\sum_{1\leq i\leq n-1}\L_{G_n\natural F_i} \oF_i^{(k)}\Bigr|
  \hfill\cr\hfill
  {}\leq \sum_{1\leq i\leq n-1} |\T_{F_n,d_n}\L_{G_{n-1}\natural F_i}
  \oF_i^{(k)} - \L_{G_n\natural F_i}\oF_i^{(k)}| \, .
  \qquad
  \equadef{PBasicLemmaCEqA}\cr
  }
$$
Since $F_n\star (G_{n-1}\natural F_i)=G_n\natural F_i$, the claim
yields for $t$ at least $1$,
$$\displaylines{\quad
    |\T_{F_n,d_n}\L_{G_{n-1}\natural F_i}\oF_i^{(k)}
    -\L_{G_n\natural F_i}\oF_i^{(k)}|(t)
  \hfill\cr\noalign{\vskip 2pt}\hfill
    {}\leq \sum_{0\leq s\leq m} {|\oF_i^{(k+s)}(t)|\over s!}
    \sum_{0\leq j\leq s} {s\choose j} \mu_{G_{n-1}\natural F_i,j}
    \int_{td_n}^\infty x^{s-j}\d F_n(x)
  \cr\qquad\qquad\qquad
    {}+{|\oF_i^{(k+m)}(t)|\over m!}
    \sum_{0\leq j\leq m} {m\choose j} \mu_{G_{n-1}\natural F_i,j}
   \hfill\equadef{PBasicLemmaCEqB}\cr\hfill
    \times \int_0^{d_nt}\oDelta_{t,x/t}^r(\oF_i^{(k+m)}) x^{m+r-j}\d F_n(x) 
    \, .
  \cr
  }
$$
Since we assume $C_\rho$ to be at most $1$, all the $c_n$'s are at most
$1$. From Lemma {\BoundDkMF}, we deduce that for $t$ at least $t_1$,
$$
  |\oF_i^{(k+s)}(t)|\leq M c_i^{\alpha-\epsilon}t^{-k-s}\oF(t) \, .
$$
Moreover, $\mu_{G_{n-1}\natural F_i,j}\leq \mu_{G,j}$. But $\mu_{G,j}$
is a fixed polynomial in the $C_l$'s and $\mu_{F,l}$'s, $1\leq l\leq j$. 
Since all these $C_l$'s are at most $C_\rho$, we conclude that 
$\mu_{G_{n-1}\natural F_i,j}$ is at most
some fixed constant $M$; this constant does not depend on $n$ and
sequences $(c_n)_{n\geq 1}$ with $C_\rho\leq 1$.

Lemma {\TruncatedMomentBoundF} implies that for any $t$ at least $t_4$,
$$
  \int_{td_n}^\infty x^{s-j}\d F_n(x)
  \leq M d_n^{s-j}c_n^{\alpha (1-\rho)-\epsilon} t^{s-j}\oF(t) \, .
$$
It follows that the sum in the first term of {\preveq} is at most
$$\displaylines{
    \sum_{0\leq s\leq m} Mc_i^{\alpha-\epsilon}t^{-k-s}\oF(t)
    \sum_{0\leq j\leq s} d_n^{s-j} c_n^{\alpha (1-\rho)-\epsilon}
    t^{s-j}\oF(t)
  \hfill\cr\hfill
    {}\leq M t^{-k} \oF(t)^2 c_i^{\alpha-\epsilon} 
    c_n^{\alpha(1-\rho)-\epsilon} \, .
  \cr
  }
$$
Similarly, for the second term in \preveq, we use the same bounds 
and Lemma {\DeltaDM} to obtain that it is at most
$$
  M c_i^{\alpha-\epsilon} t^{-k-m}\oF(t) 
  \max_{0\leq j\leq m} \int_0^{d_nt}\oDelta_{t/c_i,x/t}^r(\oF^{(k+m)}) 
  x^{m+r-j}\d F_n(x) \, .
$$
Using a change of variables, we rewrite the integral terms in this bound as
$$
  c_n^{m+r-j}\int_0^{td_n/c_n} \oDelta_{t/c_i,c_nx/t}^r(\oF^{(k+m)})x^{m+r-j}
  \d F(x) \, .
$$
As in the proof of Lemma {\BasicLemmaA}, we split the integral in this
bound as one for $x$ between $0$ and $v(t)$, and one between $v(t)$ and
$td_n/c_n$. We have
$$\displaylines{\qquad
  \int_0^{v(t)} \oDelta_{t/c_i,c_nx/t}^r(\oF^{(k+m)}) x^{m+r-j}\d F(x)
  \hfill\cr\hfill
  {}\leq \oDelta_{t,c_nv(t)/t}^r(\oF^{(k+m)})\mu_{F,m+r-j} \, ,
  \qquad\cr}
$$
and using Lemma {\TruncatedMomentBoundF} (with $u=v=1$ in it), 
for $v(t)$ at least $t_1$,
$$\displaylines{\qquad
    \int_{v(t)}^{td_n/c_n} \oDelta_{t/c_i,c_nx/t}^r(\oF^{(k+m)}) x^{m+r-j} 
    \d F(x)
  \hfill\cr\hfill
    {}\leq M\oDelta_{t/c_i,d_n}^r (\oF^{(k+m)}) v(t)^{m+r-j}\oF\circ v(t) \, .
  \qquad\cr}
$$
We conclude that for $t$ large enough {\PBasicLemmaCEqB} is at most
$$\displaylines{\quad
    M t^{-k-m}\oF(t) c_i^{\alpha-\epsilon} 
    \Bigl( t^m\oF(t)c_n^{\alpha (1-\rho)-\epsilon} 
    +c_n^r\oDelta_{t,c_nv(t)/t}^r(\oF^{(k+m)})
  \hfill\cr\hfill
    {}+ c_n^r\oDelta_{t,d_n}^r(\oF^{(k+m)}) v(t)^{m+r}\oF\circ v(t) 
    \Bigr) \, .
  \quad\cr}
$$
Since $m+r$ is smaller than $\alpha$, the function $s\mapsto s^{m+r}\oF(s)$ 
is asymptotically equivalent to a nonincreasing function (Bingham, Goldie
and Teugels, Theorem 1.5.3). Hence, for $v(t)$ large enough, $t^m\oF(t)$
is at most $2v(t)^{m+r}\oF\circ v(t)$. Therefore, {\PBasicLemmaCEqB} is
at most 
$$\displaylines{
    M t^{-k-m}\oF(t) c_i^{\alpha-\epsilon} 
    \Bigl( \bigl( c_n^{\alpha(1-\rho)-\epsilon} 
    + c_n^r\oDelta_{t,d_n}^r(\oF^{(k+m)})\bigr) v(t)^{m+r}\oF\circ v(t)
  \hfill\cr\hfill
    {}+c_n^r\oDelta_{t,c_nv(t)/t}^r(\oF^{(k+m)})\Bigr) \, .
  \cr}
$$
We then use {\PBasicLemmaCEqA} to finish.\hfill$\qed$

\bigskip

\subsection{Inductions.}%
Define
$$
  \gamma_{m,n}^{(k)}(t)=|\oGnk-\AmoGnk|_{k+m,t} \, .
$$
Our proof of Theorem {\MainTheorem} consists in bounding 
$\gamma_{m,n}^{(k)}(t)$ uniformly in $n$. This is done by induction on
$n$, and we first need to settle the case $n=2$.


\Lemma%
{\label{GammaTwo}%
  Take $\rho$ less than $\alpha/(\alpha+k)$.
  There exists a function $\eta_0$ with limit $0$ at infinity, such that 
  whenever $d_1$ and $d_2$ are at most $1/2$ and $k+m<\omega$,
  $$
    \gamma_{m,2}^{(k)}(t)\leq \eta_0(t)\, .
  $$
}

\Proof Proposition {\ConvolInOperators} shows that
$$
  \oG_2=\T_{F_1,1-d_2}\oF_2 + \T_{F_2,d_2}\oF_1 + \M_{1/(1-d_2)}\oF_1
  \M_{1/d_2}\oF_2 \, ,
  \eqno{\equa{PGammaTwoEqA}}
$$
while for $k$ at least $1$,
$$\displaylines{\qquad
    \oG_2^{(k)}=\T_{F_1,1-d_2}\oF_2^{(k)} +\T_{F_2,d_2}\oF_1^{(k)} 
  \hfill\cr\hfill
    - \sum_{1\leq j\leq k-1} \M_{1/(1-d_2)}\oF_1^{(j)}\M_{1/d_2}\oF_2^{(k-j)} 
    \, .
    \qquad\equa{PGammaTwoEqB}\cr
  }
$$
Let $r$ be a nonnegative number less than $1$ and $\omega-k-m$.
We use a slight modification of Lemma {\BasicLemmaA} with $n$ equal $2$ to 
obtain
$$\displaylines{\quad
  |(\T_{F_1,1-d_2}-\L_{F_1,m})\oF_2^{(k)}|_{k+m,t}
  \hfill\cr\noalign{\vskip 2pt}\hfill
  \leq
  M c_2^{\alpha-\rho(\alpha+k+m)-\epsilon} 
  \bigl( \oDelta_{t,v(t)/t}^r(\oF^{(k+m)}) + v(t)^{m+r} \oF\circ v(t)\bigr)
  \, .
  \quad\cr}
$$
Permuting $F_1$ and $F_2$ yields an upper bound for the error committed
in approximating $\T_{F_2,d_2}\oF_1^{(k)}$ by $\L_{F_2,m}\oF_1^{(k)}$.

Combined with {\PGammaTwoEqB}, this proves the lemma when $k$ is $1$.

When $k$ vanishes, since $m$ is smaller than $\alpha$, 
$$
  \M_{1/d_2}\oF_2(t) \M_{1/(1-d_2)}\oF_1(t) 
  = o\bigl(t^{-m}\oF(t)\bigr)
$$
as $t$ tends to infinity. Since $\M_{1/d_2}\oF_2=\M_{c_2^{1-\rho}}\oF$ this 
is uniform in $d_1$ and $d_2$ at most
$1/2$. Combined with {\PGammaTwoEqA}, this proves the result for $k$ vanishing.

When $k$ is larger than $1$, let $\epsilon$ be less than 
$\alpha (1-\rho)-k\rho$. We use Lemma {\BoundDkMF} to obtain
$$\displaylines{\qquad
  |\M_{1/d_2}\oF_2^{(i)} \M_{1/(1-d_2)}\oF_1^{(k-i)}|(t)
  \hfill\cr\noalign{\vskip 2pt}\hfill
  \eqalign{
    & {}= \bigl|c_2^{-i}\oF^{(i)}(td_2/c_2)c_1^{i-k}
      \oF^{(k-i)}\bigl( t(1-d_2)/c_1\bigr) \bigr| \cr
    &{}\leq M t^{-i} c_2^{\alpha(1-\rho)-i\rho-\epsilon} \oF (t) t^{-k+i}
      c_1^{\alpha-\epsilon}\oF(t) \cr
  }
  \qquad\cr
  }
$$
which is at most $M t^{-k}\oF(t)^2$. We conclude as before.\hfill$\qed$

\bigskip

The next step is to show that $\gamma_{m,n}^{(k)}$ can be bounded by
induction. For this purpose, we need another bound that will control
the sum involving the multiplication operators in \InductionEquality.

\Lemma
{\label{BoundMSum}%
  Let $k$ be at least $2$.
  Let $\epsilon$ be a positive real number. Assume that all $d_i$'s
  are smaller than $1/2$. There exist $t_9$ and $M$ such that for any $t$
  at least $t_9$, any positive integers $j$, $k$ with $j\leq k-1$,
  $$\displaylines{\qquad
      \sum_{i\geq 2} |\M_{1/d_i}\oF_i^{(j)}\M_{1/(1-d_i)}
      \oG_{i-1}^{(k-j)}|_{k+m,t}
    \hfill\cr\hfill
      {}\leq M t^m\oF(t) \bigl( C_{\alpha-\epsilon/(1-\rho)}+\sup_{i\geq 2}
      \gamma_{0,i}^{(k-j)}(t/2)\bigr) C_{\alpha-\rho(\alpha+j)-\epsilon}
      \, . \cr
    }
  $$
}

\Proof Define $\delta=\epsilon/(1-\rho)$. Under our assumptions, 
all the $c_i$'s are at most $1$. Consequently,
$d_i/c_i$ is at least $1$ and Lemma {\BoundDkMF} (taking $c$ to be $1$ in 
that lemma) implies that for $s$ at least $t_1$,
$$
  |\M_{1/d_i}\oF_i^{(j)}(s)|
  = c_i^{-j}|\oF^{(j)}(sd_i/c_i)|
  \leq c_i^{-j}M (sd_i/c_i)^{-j} \oF(sd_i/c_i) \, .
$$
Applying Lemma \Potter, using again the fact that $d_i/c_i$ is at least
$1$, we see that for $s$ at least $t_1$ and $t_2$,
$$
  |\M_{1/d_i} \oF_i^{(j)}(s)|
  \leq M c_i^{-j+(j+\alpha-\delta)(1-\rho)} s^{-j}\oF(s) \, .
$$
Since $\calA_0\oG_i^{(k-j)}=\sum_{1\leq l\leq i}\oF_l^{(k-j)}$, the 
triangle inequality yields
$$
  |\oG_i^{(k-j)}(s)|
  \leq \sum_{1\leq l\leq i} |\oF_l^{(k-j)}(s)| + \gamma_{0,i}^{(k-j)}(s)
  s^{-k+j}\oF(s) \, .
$$
Using again Lemma {\BoundDkMF} and the fact that the $c_i$'s are at most
$1$, we see that for $s$ at least $t_1$,
$$\eqalign{
  |\oG_i^{(k-j)}(s)|
  &{}\leq M\sum_{1\leq l\leq i} c_l^{\alpha-\delta} s^{-k+j}\oF(s)
    + \gamma_{0,i}^{(k-j)}(s) s^{-k+j} \oF(s) \cr
  &{}\leq M s^{-k+j} \oF(s) 
    \bigl( C_{\alpha-\delta}+\sup_{n\geq 2}\gamma_{0,n}^{(k-j)}(s)\bigr)
    \, . \cr
  }
$$
Consequently, since $d_i$ is at most $1/2$, we obtain that for $t$ large 
enough,
$$\displaylines{\quad
  \sum_{i\geq 2} \sup_{s\geq t} s^{k+m}
  |\M_{1/d_i}\oF_i^{(j)}(s)\M_{1/(1-d_i)}\oG_{i-1}^{(k-j)}(s)|/\oF(s)
  \hfill\cr\hfill
  {}\leq M\sup_{s\geq t} s^m \oF(s/2)
  \bigl( C_{\alpha-\delta}+\sup_{n\geq 2}\gamma_{0,n}^{(k-j)}(s/2)\bigr)
  C_{(j+\alpha-\delta)(1-\rho)-j} \, ,
  \quad\cr
  }
$$
the factor $1/2$ in $s/2$ coming from the multiplication operator acting
on $\oG_{i-1}^{(k-1)}$ and our assumption that $d_i$ is at
most $1/2$.
One more application of Potter's bound allows us to replace $\oF (s/2)$ by
$\oF(s)$ in the upper bound. Since $\gamma_{0,i}^{(k-j)}(\cdot )$ is 
nonincreasing and $s\mapsto s^m\oF(s)$ is asymptotically equivalent to an
ultimately nonincreasing function  (Bingham, Goldie and Teugels, 1989, 
Theorem 1.5.3), we can remove the supremum in $s$
in the above bound (and, of course, we increase $M$ in doing 
that).\hfill$\qed$

\bigskip

Our next lemma shows that $\gamma_{m,n}^{(k)}(\cdot )$ can be bounded 
by induction on $k$, whenever $k$ is at least $2$. We define
$$
  \eta_\sigma(t)
  = \bigl( v(t)/t\bigr)^\gamma + v(t)^\sigma\oF\circ v(t) \, .
$$

\Lemma%
{\label{InductionBound}%
  Let $\epsilon$ be a positive real number, and let $k$ be a positive
  integer. There exist positive numbers $M$
  and $t_{10}$ such that for any integer $n$ at least $2$ and any sequence
  $(c_n)_{n\geq 1}$ with $C_{\gamma\rho}\vee 
  C_{\alpha-(\alpha+k+m)\rho-\epsilon}\leq 1$ and $d_n\leq 1/2$,
  $$\eqalign{
    \gamma_{m,n}^{(k)}(t)\leq{}
    & M\Bigl( 1+\max_{1\leq j\leq k-1}\sup_{i\geq 2} 
      \gamma_{0,i}^{(k-j)}(t/M)\Bigr) \eta_m (t/M) \cr
    & {}+M\gamma_{m,2}^{(k)}(t/M) \, . 
    \cr}
  $$
}%
In the statement of this lemma, we agree that $\max_{1\leq j\leq 0}$
is $0$.

\bigskip

\Proof For $u$ positive, the inequality 
$\oDelta_{t,u}^r(h)\leq u^{\gamma-r}\oDelta_{t,u}^\gamma (h)$ 
holds. We will take $r$ to be $0$ in this inequality and apply it in
the bounds of the previous subsection. 
Under the assumption of Theorem \MainTheorem, we also have
$\limsup_{t\to\infty}\oDelta_{t,1/2}^0(\oF^{(k+m)})$
finite.
Using {\InductionEquality}, Lemma \BasicLemmaA, 
Corollary {\CorToTNearContraction} and \BasicLemmaC, we obtain that
there exists $M$ and $t_{10}'$ such that for any $t$ at least $t_{10}'$ and
any $n$ at least $3$,
$$\eqalign{
  \gamma_{m,n}^{(k)}(t)\leq{}
  & M c_n^{\alpha-(\alpha+k+m)\rho -\epsilon}
    \Bigl( \bigl( v(t)/t\bigr)^\gamma+ v(t)^m\oF\circ v(t) \Bigr) \cr
  &{}+ \gamma_{m,n-1}^{(k)}\bigl( t(1-d_n)\bigr) (1+Md_n) \cr
  \noalign{\vskip 3pt}
  &{}+ M(c_n^{\alpha(1-\rho)-\epsilon}+d_n^\gamma) v(t)^m\oF\circ v(t)
    + Mc_n^\gamma \bigl( v(t)/t\bigr)^\gamma \cr
  \noalign{\vskip 3pt}
  &{} +\sum_{1\leq j\leq k-1} 
    |\M_{1/d_n}\oF_n^{(j)}\M_{1/(1-d_n)}\oG_{n-1}^{(k-j)}|_{k+m,t} \, . \cr
  }
$$
Collecting the terms and using that all the $c_n$'s are at most $1$,
$$\eqalignno{
  \gamma_{m,n}^{(k)}(t)\leq{}
  & M(c_n^{\alpha-(\alpha+k+m)\rho-\epsilon}+c_n^{\gamma\rho})\eta_m(t) 
    \cr
  \noalign{\vskip 2pt}
  &{}+\gamma_{m,n-1}^{(k)}\bigl( t(1-d_n)\bigr) (1+Md_n) 
  &\equa{InductionBoundEqA}\cr
  \noalign{\vskip 2pt}
  &{}+\sum_{1\leq j\leq k-1} 
    |\M_{1/d_n}\oF_n^{(j)}\M_{1/(1-d_n)}\oG_{n-1}^{(k-j)}|_{k+m,t} \, . \cr
  }
$$
Taking $t_{10}'$ large enough, we can assume that $\eta_m (\cdot )$ is 
nonincreasing on $[\,t_{10}',\infty)$. Dropping the subscript $m$ and 
superscript $k$ temporarily, inequality {\InductionBoundEqA} has the form
$$
  \gamma_n(t)\leq a_n(t)+\gamma_{n-1}\bigl( t(1-d_n)\bigr) (1+Md_n) \, ,
$$
where $a_n(\cdot)$ is the nonnegative and nonincreasing function
on $[\, t_{10}',\infty)$ given by
$$\displaylines{\qquad
  a_n(t)= M(c_n^{\alpha-(\alpha+k+m)\rho-\epsilon}+c_n^{\rho\gamma})
  \eta_m(t)
  \hfill\cr\hfill
  + \sum_{1\leq j\leq k-1}
  |\M_{1/d_n}\oF^{(j)}_n\M_{1/(1-d_n)}\oG_{n-1}^{(k-j)}|_{k+m,t}
  \, .
  \qquad\cr}
$$

For $i+1$ at most $n$, define $A_{i,n}=\prod_{i+1\leq j\leq n} (1-d_j)$
and $B_{i,n}=\prod_{i+1\leq j\leq n}(1+Md_j)$. We also set $A_{n,n}=B_{n,n}=1$.
By induction, inequality {\InductionBoundEqA} implies
$$
  \gamma_n(t)\leq \sum_{3\leq i\leq n} B_{i,n}a_i(tA_{i,n}) 
  + B_{2,n}\gamma_2(tA_{2,n}) \, .
$$
Set $A=\prod_{j\geq 1} (1-d_j)$ and $B=\prod_{j\geq 1}(1+Md_j)$. Since all the
$a_n(\cdot)$'s are nonnegative and nonincreasing, we have
$$
  \gamma_n(t)\leq \sum_{i\geq 3} B a_i(tA)+B\gamma_2(tA)\, .
$$
The inequality $\log (1-x)\geq -2x$ for $x$ nonnegative and at most $1/2$
implies
$$
  A\geq \exp\Bigl( -2\sum_{j\geq 1} d_j\Bigr) = \exp (-2C_\rho ) \, ,
$$
while the the inequality $\log (1+x)\leq x$ implies
$$
  B\leq \exp\Bigl( \sum_{j\geq 1} Md_j\Bigr) = \exp(MC_\rho) \, .
$$
Consequently,
$$
  \gamma_n(t)\leq e^{MC_\rho}\Bigl( \sum_{i\geq 3} a_i(te^{-2C_\rho})
  +\gamma_2(te^{-2C_\rho})\Bigr) \, .
$$
We also have
$$\displaylines{\qquad
    \sum_{i\geq 3} a_i(s) 
    \leq M (C_{\alpha-(\alpha+k+m)\rho-\epsilon}+C_{\rho\gamma})\eta_m(s)
  \hfill\cr\hfill
    {}+\sum_{1\leq j\leq k-1} \sum_{i\geq 3} 
    |M_{1/d_i}\oF_i^{(j)}\M_{1/(1-d_i)}\oG_{i-1}^{(j-k)}|_{k+m,s} \, .
  \qquad\cr}
$$
Applying Lemma \BoundMSum, this bound is at most
$$\displaylines{
    M(C_{\alpha-(\alpha+k+m)\rho-\epsilon}+C_{\rho\gamma})\eta_m(s)
  \hfill\cr\noalign{\vskip 2pt}\hfill
    {}+ M s^m\oF(s) \bigl( C_{\alpha-\epsilon/(1-\rho)}+
    \max_{1\leq j\leq k-1}\sup_{i\geq 2} \gamma_{0,i}^{(k-j)}(s/2)\bigr)
    C_{\alpha-\rho(\alpha+k)-\epsilon} \, .
  \cr}
$$
Since $s^m\oF(s)$ is at most $\eta_m (s)$ for $s$ at least $1$, this 
bound yields
\hfuzz=3pt
$$\displaylines{%
    \gamma_n(t)\leq M\Bigl( C_{\alpha-(\alpha+k+m)\rho-\epsilon}
    + C_{\gamma\rho} +1+\max_{1\leq j\leq k-1} \sup_{i\geq 2} 
    \gamma_{0,i}^{(k-j)} (te^{-2C_\rho}/2) \Bigr)
  \hfill\cr\hfill
    \times \eta_m (te^{-2C_\rho}) + e^{MC_\rho}\gamma_2(te^{-2C_\rho}) \, .
  \qquad\cr}
$$
\hfuzz=0pt
To conclude the proof, if $C_{\rho\gamma}$ is at most $1$, so are
all the $c_i$'s and $C_\rho$ is at most $C_{\rho\gamma}$.\hfill$\qed$

\bigskip

We need an analogue of Lemma {\InductionBound} when $k$ vanishes. For this
purpose, we state the analogue of Lemma {\BoundMSum} in this case.

\bigskip

\Lemma
{\label{BoundMSumkZero}
  Assume that $C_\rho$ is at most $1$ and that all $d_i$'s are at most
  $1/2$. For $t$ larger than some $t_{11}$, the following inequality holds:
  $$
    \sum_{i\geq 1} |M_{1/d_i}\oF_i\M_{1/(1-d_i)}\oG_{i-1}|_{k+m,t}
    \leq C_{\alpha(1-\rho)-\epsilon}^2 t^m\oF(t) \, .
  $$
}

\Proof
Set $\delta=\epsilon/(1-\rho)$.
Apply Lemma {\Potter} to obtain that for $s$ at least $t_2$,
$$
  \M_{1/d_i}\oF_i(s)
  =\oF(sd_i/c_i)
  \leq c_i^{(1-\rho)(\alpha-\delta)}\oF(s) \, .
$$
Since $d_i$ is at most $1/2$, Lemma {\TrivialTailBound} implies
$$
  \M_{1/(1-\d_i)}\oG_{i-1}(s) \leq \sum_{n\geq 1} \oF(s/2c_n^{1-\rho}) \, ,
$$
which, by Lemma {\Potter} is at most $M\sum_{n\geq 1} 
c_n^{(1-\rho)(\alpha-\delta)}\oF(s)$ for $s$ large enough. This implies
the conclusion.\hfill$\qed$

\bigskip

The next lemma is Theorem {\MainTheorem} stated in an other way.

\Lemma
{\label{FinalLemmaPositive}%
  Let $\epsilon$ be a positive number at most $1$. There
  exists a function $\eta (\cdot)$ with limit $0$ at infinity, such that,
  for any $k$ and $m$ with $k+m$ smaller than $\omega$ and $m$ smaller than
  $\alpha$, any $n$ at least $2$, and 
  any $t$ at least $1$ say, any nonnegative sequence $(c_i)_{i\geq 1}$ 
  with $\sup_{i\geq 1} d_i\leq 1/2$ and both $C_{\gamma\rho}$ and 
  $C_{\alpha-(\alpha+k+m)\rho-\epsilon}$ at most $1$,
  $$
    \gamma_{m,n}^{(k)}(t) \leq \eta (t) \, .
  $$
}

\Proof
For $k=1$, Lemma {\InductionBound} implies
$$
  \gamma_{m,n}^{(1)}(t)\leq M \eta_m (t/M) + M\gamma_{m,2}^{(1)}(t/M) \, .
$$
The result follows then from the assumptions and Lemma \GammaTwo.

By induction on $k$, Lemmas {\InductionBound} and {\BoundMSum} imply the 
result for $k$ at least $1$.

When $k$ vanishes, Lemma {\BoundMSumkZero} and the same arguments as 
those in the proof of Lemma {\InductionBound} show that
$$
  \gamma_{m,n}^{(0)}(t)\leq M (C_{(1-\rho)(\alpha-\epsilon)}+1)
  \eta_m(t/M) + M\gamma_{m,2}^{(0)} + MC_{\alpha (1-\rho)-\epsilon}^2
  t^m\oF(t)
$$
(compare with \InductionBoundEqA). The result follows.\hfill$\qed$

\bigskip

\subsection{Conclusion.}%
To obtain Theorem \MainTheorem, it mostly remains to show that as $n$
tends to infinity, $\gamma_{m,n}(t)$ converges to 
$|\oG^{(k)}-\calA_m\oG^{(k)}|_{k+m,t}$. This is achieved in two steps,
one consisting in proving that the sequence of approximations converges,
the other one in proving that the sequence of functions $G_n^{(k)}$
converges to $G^{(k)}$.

\bigskip

\Lemma
{\label{ContinuityAm}%
  Assume that $C_{1-\rho}$ and $|c|_\infty$ are at most $1$. 
  Then there exists $t_{12}$ and
  $M$ such that for any $t$ at least $t_{12}$,
  $$
    |\calA_m\oG_n^{(k)}-\calA_m\oG^{(k)}|_{k,t}
    \leq M \sum_{i\geq n+1}c_i^{1-\rho}
  $$%
}%
\finetune{\vskip -.1in}
\Remark For our purpose, it is enough that $\calA_m\oG_n^{(k)}$ converges
pointwise to $\calA_m\oG^{(k)}$. In order to simulate properly tail behavior,
it would be desirable to have convergence in $|\cdot|_{k+m,t}$ norm.

\bigskip

\Proof From the definition of the approximation and since both
$\mu_{G_n\natural F_i,0}$ and $\mu_{G\natural F_i,0}$ equal $1$ for $i$ 
at most $n$,
$$\eqalignno{
  \calA_m\oG_n^{(k)}-\calA_m\oG^{(k)}
  &{}= \sum_{0\leq j\leq m} {(-1)^j\over j!}\sum_{1\leq i\leq n}
    (\mu_{G_n\natural F_i,j}-\mu_{G\natural F_i,j}) \D^{j+k}\oF_i \cr
  &{}\qquad {}-\sum_{0\leq j\leq m} {(-1)^j\over j!} \sum_{i\geq n+1}
    \mu_{G\natural F_i,j}\D^{j+k}\oF_i \, .
  &\equa{BoundForContinuityAm}
  \cr
  }
$$
Write $Y=\sum_{\matrix{\scriptstyle1\leq j\leq n\cr\noalign{\vskip -3pt}
\scriptstyle j\not= i\hfill\cr}} c_jX_j$ and $R=\sum_{\matrix{\scriptstyle 
j\geq n+1\cr\noalign{\vskip -3pt}\scriptstyle j\not= i\hfill\cr}} c_jX_j$.
Then, for $j$ at least $1$,
$$
  \mu_{G\natural F_i,j}-\mu_{G_n\natural F_i,j}
  =E\bigl( (Y+R)^j-Y^j\bigr)
  =\sum_{1\leq l\leq j} {j\choose l} E R^l E Y^{j-l} \, .
$$
Applying Lemma {\MomentBoundG}, we see that for any $l$ at least $1$,
$$
  E R^l
  \leq \sum_{i\geq n+1} c_i^{l(1-\rho)}\mu_{F,l}
  \leq \sum_{i\geq n+1} c_i^{1-\rho}\mu_{F,l} \, ,
$$
while for $l$ less than $j$, by the same token and under the assumption of the
lemma,
$$
  EY^{j-l}
  \leq \mu_{G,j-l}
  \leq C_{(j-l)(1-\rho)}\mu_{F,j-l}
  \leq \mu_{F,j-l} \, .
$$
Consequently, for $j$ at least $1$,
$$\eqalign{ 0\leq \mu_{G\natural F_i,j}-\mu_{G_n\natural F_i,j}
  & {}\leq \sum_{1\leq l\leq j} {j\choose l} \mu_{F,l}\mu_{F,j-l}
    \sum_{i\geq n+1} c_i^{1-\rho} \cr
  & {}\leq\mu_{F\star F,j} \sum_{i\geq n+1} c_i^{1-\rho} \, . \cr
  }
$$
For $t$ at least $t_1$, Lemma {\BoundDkMF} implies
$$
  |\D^{j+k}\oF_i(t)|\leq Mc_i^{\alpha-\epsilon} t^{-k-j} \oF(t) \, .
$$
So, the first double sum in {\BoundForContinuityAm} is at most
$$
  M t^{-k} \oF(t) \sum_{i\geq n+1} c_i^{1-\rho} \, .
$$
To bound the second double sum, we use Lemma {\MomentBoundG} to obtain,
for $j$ at least $1$.
$$
  \mu_{G\natural F_i,j} 
  \leq \mu_{G,j}
  \leq C_{1-\rho} \mu_{F,j} \, .
$$
Hence, the second double sum in {\BoundForContinuityAm} is at most
$$
  M\sum_{i\geq n+1} c_i^{\alpha-\epsilon} t^{-k} \oF(t) \, .
$$
This concludes the proof.\hfill$\qed$

\bigskip

The next lemma is the only place where the boundedness and continuity of 
$F^{(k)}$ is used. Its proof is an adaptation of that of Proposition 9.1.6 
in Dudley (1989).

\Lemma
{\label{GnkCvToGk}%
  Assume that $F$ is $k$-times continuously differentiable
  on $(0,\infty )$, and that $F^{(k)}$ is bounded and in $L^1(\d x)$. Then,
  $\limn G_n^{(k)}=G^{(k)}$ pointwise.
}

\bigskip

\Proof
Write 
$$
  G_n(t)=\int_0^t F_1(t-y)\d G_n\natural F_1(y) \, .
$$
Since $F^{(k)}$ exists and is in $L^1(\d x)$, so is $F_1$. Then
$$
  G_n^{(k)}(t)=\int_0^t F_1^{(k)}(t-y)\d G_n\natural F_1(y) \, .
$$
But the sequence $(G_n\natural F_1)_{n\geq 1}$ converges weakly* to
the continuous distribution function $G\natural F_1$ and since $F_1^{(k)}$ 
is continuous and bounded, $G_n^{(k)}(t)$ converges to 
$$
  \int_0^t F_1^{(k)}(t-y)\d G\natural F_1(y)
  = \bigl( F_1\star (G\natural F_1)\bigr)^{(k)}(t)
  = G^{(k)}(t) \, .
  \eqno{\qed}
$$

To obtain Theorem {\MainTheorem} when the $c_i$'s are nonnegative, it 
remains to do some rewriting of Lemma {\FinalLemmaPositive}. For this purpose,
define the sets of nonnegative sequences
$$
  \calC_{\alpha,\omega,\gamma}
  =\{\, c\in [\,0,\infty)^{\NN^*} \, : \, N_{\alpha,\gamma,\omega}(c)\leq 1 
  \,\}\, .
$$
The following is exactly Theorem {\MainTheorem} but for the positive case
and everywhere smooth distribution function.

\Proposition
{\label{FinalPropositionPositive}
  Let $\oF$ in $SR_{-\alpha,\omega}$ and let $k$, $m$ be two nonnegative
  integers with $k+m$ smaller than $\omega$ and $m$ smaller than $\alpha$. 
  Let $\gamma$ be a positive
  real number smaller than both $1$ and $\omega-k-m$.
  Assume that $F^{(k)}$ is countinuous and bounded.
  There exists a function
  $\eta(\cdot)$ which tends to $0$ at infinity and a real number $t_{13}$
  such that for any $t$ at least $t_{13}$ and any sequence $c$ in
  $\calC_{\alpha,\omega,\gamma}$
  $$
    |\oG_c^{(k)}-\calA_m\oG_c^{(k)}|_{k+m,t}\leq \eta (t) \, .
  $$
}

\Proof
For $\epsilon$ positive, define 
$\rho=(1/2)\wedge \bigl( \alpha/(\alpha+\omega)\bigr)$. Let $c$
be a sequence in $\calC_{\alpha,\omega,\gamma}$. By definition of 
$\calC_{\alpha,\omega,\gamma}$ the series
$C_{\rho\gamma}$ is at most $1$. Moreover, by taking $\epsilon$ small enough,
$C_{\alpha-\epsilon}$ is at most $C_{\rho\gamma}$, hence at most $1$.

Since $k+m+\gamma<\omega$, we also have $\rho<\alpha/(\alpha+\gamma+k+m)$,
that is $\alpha-\rho(\alpha+k+m)>\rho\gamma$.
Hence, for $\epsilon$ small enough $C_{\alpha-(\alpha+k+m)\rho-\epsilon }$
is at most $1$.

To conclude the proof, lemmas {\ContinuityAm} and {\GnkCvToGk} show that
$$
  \oG^{(k)}-\calA_m\oG^{(k)} = \limn \oG_n^{(k)}-\calA_m \oG_n^{(k)}
$$
pointwise. But for any $n$ at least $2$ and any $t$ at least $t_1$, 
Lemma {\FinalLemmaPositive} yields
$$
  |\oG_n^{(k)}(t)-\calA_m\oG_n^{(k)}(t)|\leq \eta (t) t^{-m-k}\oF(t)\, .
$$
Hence the same inequality holds with $\oG^{(k)}$.\hfill$\qed$

\bigskip


\section{Removing the sign restriction on the random variables.}

In section 5, we assumed that the random variables are positive.
Our goal in this section is to remove this assumption. But we will keep
the assumption that $F^{(k)}$ and $\M_{-1}F^{(k)}$ exist, are bounded 
and continuous on the whole positive half line.
The basic
argument consists in conditioning. To be specific, define the
random set of all indices corresponding to a nonpositive random variable,
that is
$$
  I=\{\, i\in \NN^* \, :\, X_i\leq 0\,\} \,.
$$
Write $H_I$ for the distribution function of $\sum_{i\in I}c_i X_i$
and $G_I$ for that of $\sum_{i\in \NN^*\setminus I} c_iX_i$. Let $F_+$
(resp. $F_-$) be the distribution function of $X_1$ say, given that
$X_1$ is positive (resp. nonpositive), that is
$$\eqalignno{
  \overline{F_+}= \oF/\oF(0) &\quad\hbox{ on}\quad (0,\infty) \cr
\noalign{\noindent and} 
  F_-=F/F(0) &\quad\hbox{ on}\quad (-\infty,0\,]\, .\cr
}
$$
We write $F_{+,i}$ (resp. $F_{-,i}$) for $\M_{c_i}F_+$ (resp. $\M_{c_i}F_-$).
Note that there is no ambiguity in this notation for when $c$ is nonnegative,
$(\M_cF)_+=\M_c(F_+)$ and $(\M_c F)_-=\M_c(F_-)$.
Given $I$, the results of section 5 can be applied to both $\M_{-1}H_I$ and 
$G_I$, using respectively the distribution functions $M_{-1}F_-$ and $F_+$ 
instead of $F$. The distribution function of the whole series $\<c,X\>$ 
is $EH_I\star G_I$. The identity
$$
  \overline{H_I\star G_I}(t)=\int_{-\infty}^0 \oG_I(t-x)\d H_I(x) \, ,
$$
suggests the relevance to the operator
$$
  \U_Hh(t) = \int_{-\infty}^0 h(t-x) \d H(x) \, ,
$$
in term of which $\overline{H_I\star G_I}=\U_{H_I}\oG_I$. We then need
to prove the analogue of some of the results of section 5 on $\T_{F,\eta}$,
but now for the operator $\U_H$. The first one which we will prove implies
$\D^k \U_H\oG = \U_H\oG^{(k)}$. It is then clear how the proof of the
asymptotic expansion for the distribution function $K$ of $\<c,X\>$ and 
its derivatives will go. Indeed, we will have
$$
  K^{(k)}
  = E \U_{H_I}\oG_I^{(k)}
  = E\U_{H_I}\calA_m\oG_I^{(k)} + E\U_{H_I}(\oG_I^{(k)}-\calA_m\oG_I^{(k)}) 
  \, .
  \eqno{\equa{KEqual}}
$$
We will see that $\U_H$ is a contraction for the right norm. Hence, 
equality {\KEqual} combined with the results of section 5 gives the asymptotic
expansion
$$
  K^{(k)}\sim E\U_{H_I}\calA_m\oG_I^{(k)} \, .
$$
We will be able to approximate $\U_{H_I}$ by $\L_{H_I}$. Because the error
term is bounded uniformly with respect to the sequence $(c_i)_{i\geq 1}$,
we will be able to permute expectation and asymptotic expansions.

\bigskip

\subsection{Elementary properties of \poorBold{$\U_H$}.}%
Since we are interested in expansions of derivatives, our first elementary
result deals with the composition of the derivative and $\U_H$. It asserts
that those two operators commute whenever acting on sufficiently regular
functions.

\Lemma
{\label{DUCommute}
  Let $h$ be a function on the nonnegative half line with Lebesgue
  integrable $k$-th derivative, and such that for any nonnegative $l$ at 
  most $k$,
  $$
    \lim_{t\to\infty} h^{(l)}(t) = 0 \, .
  $$
  Then, $\D^k\U_H h = \U_H \D^k h$ almost everywhere.
}

\bigskip

\Proof By induction it suffices to prove the result for $k$ equal $1$. We write
$$
  \U_H h(t) = -\int\int h'(y-x) \II\{\, y>t\,\} \d y\d H(x) \, .
$$
The function $U_Hh'$ is Lebesgue integrable for
$$
  \int \Bigl| \int h'(y-x)\d H(x) \Bigr| \d y
  \leq \int\int |h'(y-x)| \d y \d H(x)
  = |h'|_{L^1(\d x)} \, .
$$
Applying Fubini's theorem,
$$
  \U_H h(t)
  = -\int_t^\infty \int_{-\infty}^0 h'(y-x)\d H(x) \d y
  = -\int_t^\infty \U_H h'(y) \d y \, .
$$
Since we proved that $\U_Hh'$ is Lebesgue integrable, it is almost everywhere
the derivative of $\U_Hh$.\hfill$\qed$

\bigskip

In particular, if $F^{(k)}$ is continuous and bounded and $i$ is not in $I$, 
then the analogue to the expression obtained in the proof of Lemma \GnkCvToGk,
is
$$
  \oG_I^{(k)}=\int_0^t \oF_{+,i}^{(k)}(t-y) \d G_I\natural F_{+,i}(y) \, .
  \eqno{\equa{ExprGI}}
$$
Thus, $\oG_I^{(k)}$ is bounded and integrable with respect to the 
Lebesgue measure. Moreover, since $\oF$ is smoothly varying with negative
exponent $-\alpha$, its $k$-th derivative tends to $0$ at infinity.
Then {\ExprGI} shows that $\oG_I^{(l)}$ tends to $0$ at infinity for any 
nonnegative $l$ at most $k$. Consequently,
$$
  \overline{H_I\star G_I}^{(k)} = \U_{H_I}\oG^{(k)}_I \, .
$$

Our next result allows us to replace $\oG_I^{(k)}$ by its asymptotic expansion
when looking for an expansion for $\U_{H_I}\oG_I^{(k)}$.

\Lemma
{\label{UContraction}%
  For any nonnegative $t$, the operator $\U_H$ 
  is a contraction with respect to the norms $|\cdot |_{p,t}$.
}

\bigskip

\Proof Let $h$ be a function whose $|\cdot |_{p,t}$-norm is finite. Since
$t\mapsto t^p/\oF(t)$ is nondecreasing on the nonnegative half line
$$\eqalign{
  \Bigl| {t^p\over \oF(t)} \U_Hh(t)\Bigr|
  &{}\leq \int_{-\infty}^0 {t^p\over \oF(t)} |h(t-x)| \d H(x) \cr
  &{}\leq \int_{-\infty}^0 {(t-x)^p\over \oF(t-x)} |h(t-x)| \d H(x) \, .\cr
  }
$$
The integrand is at most $|h|_{p,t}$; so is the integral, for $H$ is
a distribution function.\hfill$\qed$

\bigskip

\subsection{Basic expansion of \poorBold{$\U_H$}.}%
We now show that $\U_H$ has an expansion $\L_{H,m}$, and therefore behaves 
similarly to $\T_{F,\eta}$. Indeed, the next
lemma should be compared to Theorem \ApproxTByL.

\Lemma
{\label{ApproxUByL}%
  Let $m$ be a positive integer and let $r$ be in $[\, 0,1)$. Let $H$ be a 
  distribution funciton on the nonpositive half line with finite $m$-th 
  moment. If $h$ is smoothly varying of order $m+r$, then
  $$\eqalign{
    |(\U_H-\L_{H,m})h|(t)
    &{}\leq \sum_{0\leq j\leq m} {|h^{(j)}(t)|\over j!} \int_{-\infty}^{-t/2}
      |x|^j \d H(x) \cr
    &\hskip 18pt {}+{|h^{(m)}(t)|\over t^r m!} \int_{-t/2}^0 
      \oDelta_{t,|x|/t}^r (h^{(m)}) |x|^{m+r} \d H(x)\cr
    &\hskip 18pt {}+ H(-t/2) \sup_{s\geq t} |h(s)| \, . \cr
  }
  $$
}

\Proof We first have the inequality
$$\eqalign{
  \Bigl| \U_H h(t) -\int_{-t/2}^0 h(t-x) \d H(x) \Bigr|
  &{}\leq \int_{-\infty}^{-t/2} |h(t-x)| \d H(x) \cr
  \noalign{\vskip 2pt}
  &{}\leq H(-t/2) \sup_{s\geq t} |h(s)| \, . \cr
  }
$$
Applying Proposition \Taylor, we have
$$\displaylines{\qquad
    \Bigl| \int_{-t/2}^0 \Bigl( h(t-x) -\sum_{0\leq j\leq m} {(-1)^j\over j!}
    x^j h^{(j)} (t)\Bigr) \d H(x) \Bigr|
  \hfill\cr\noalign{\vskip 3pt}\hfill
    {}\leq \int_{-t/2}^0 {|x|^{m+r}\over t^r} {|h^{(m)}(t)|\over m!}
     \oDelta_{t,|x|/t}^r (h^{(m)}) \d H(x) \, . 
  \qquad\cr
  }
$$
The result follows.\hfill$\qed$

\bigskip

\subsection{A technical lemma.}%
Looking at equality \KEqual, we need to approximate
$\U_{H_I}\calA_m\oG^{(k)}_I$. The following result, which is the
analogue to Lemma {\BasicLemmaC} will do. Since $\M_{-1}H_I$ has the
same properties as the distribution function $G$ studied in section 5,
it satisfies the assumptions of the next lemma. Note that even if
$\M_{-1}F$ is not regularly varying, the assumption
$\overline{\M_{-1}F}=O(\oF)$ and Lemma \TrivialTailBound --- for its
part which does not require regular variation --- ensure that $H_I$
still satisfies the assumption of the following Lemma.

\Lemma
{\label{BasicLemmaU}%
  Assume that there exist constants $M_0$ and $t_{14}$ such that for any 
  integer $k$ less than $\alpha$ and any $t$ at least $t_{14}$,
  $$
    \int_{-\infty}^{-t} |x|^k\d H(x) \leq M_0 t^k F(-t)
  $$
  \hfuzz=1pt
  and for any nonnegative integer $j$ at most $m$, the moments 
  $\mu_{\M_{-1}H,j+r}$ are at most $M_0$. Then, there exist a positive $M$ 
  and a $t_{15}$ such that for any $G=\star_{i\geq 1}F_i$ as in section 5 
  (that is the distribution function of an infinite weighted sum of 
  nonnegative random 
  variables with nonnegative weights with the assumption of section 5 
  satisfied), for any $t$ at least $t_{15}$,
  $$\displaylines{\quad
      | \U_H\calA_m\oG^{(k)}-\sum_{i\geq 1} \L_{H\star G\natural F_i,m}
        \oF_i^{(k)} |_{k+m,t} 
    \hfill\cr\hfill
      {}\leq M M_0 C_{\alpha-\epsilon} \Bigl( \oDelta_{t,v(t)/t}^r(\oF^{(k+m)})
      + \oDelta_{t,1/2}^r (\oF^{(k+m)}) v(t)^{m+r} F\bigl(-v(t)\bigr) 
      \cr\hfill
      {}+t^m H(-t/2) \Bigr) \, .
    \cr}
  $$
  \hfuzz=0pt
}

The proof uses an auxiliary result which we state as a claim, very much as
we did in the proof of Lemma \BasicLemmaC.

\Claim
{ For any nonnegative $t$, any distribution function $H$ supported on the 
  nonpositive half line and $K$ supported on the nonnegative half line,
  $$
    \displaylines{%
     |\U_H\L_K h- \L_{H\star K}h|(t)
     \hfill\cr\noalign{\vskip 3pt}\hfill
     \eqalign{
      {}\leq{}
      & \sum_{0\leq s\leq m} {|h^{(s)}(t)|\over s!} \sum_{0\leq j\leq s}
        {s\choose j} \mu_{K,j}\int_{-\infty}^{-t/2} |x|^{s-j} \d H(x) \cr
      & {}+ {|h^{(m)}(t)|\over t^rm!} \sum_{0\leq j\leq m} {m\choose j}
        \mu_{K,j} \int_{-t/2}^0 \oDelta_{t,|x|/t}^r (h^{(m)}) |x|^{m-j+r}
        \d H(x) \cr
      & {} + H(-t/2)\sum_{0\leq j\leq m} {\mu_{K,j}\over j!}
       \sup_{s\geq t} |h^{(j)}(s)| \, . \cr
     }
    \cr}
  $$
}

\Proof By linearity of $\U_H$,
$$
  \U_H\L_K=\sum_{0\leq j\leq m} {(-1)^j\over j!} \mu_{K,j} \U_H \D^j \, .
$$
Applying Lemma {\ApproxUByL} to bound $(\U_H-\L_{H,m-j})\D^jh$, we obtain
$$\displaylines{\quad
  |\U_H\L_K h - \sum_{0\leq j\leq m} {(-1)^j\over j!} \mu_{K,j}\L_{H,m-j}
  h^{(j)}|(t)
  \hfill\cr\noalign{\vskip 3pt}\hfill
  \eqalign{
  {}\leq{}
  & \sum_{0\leq j\leq m} {\mu_{K,j}\over j!} \sum_{0\leq l\leq m-j}
    {|h^{(l+j)}(t)|\over l!} \int_{-\infty}^{-t/2} |x|^l\d H(x) \cr
  &{}+ \sum_{0\leq j\leq m} {\mu_{K,j}\over j!} {|h^{(m)}(t)|\over t^r (m-j)!}
    \int_{-t/2}^0 \oDelta_{t,|x|/t}^r (h^{(m)}) |x|^{m-j+r} \d H(x) \cr
  &{}+\sum_{0\leq j\leq m} {\mu_{K,j}\over j!} H(-t/2) \sup_{s\geq t}
    |h^{(j)}(s)| \, .\cr
  }
  \cr
  }
$$
Again, Lemma {\LaplaceBinomial} shows that the left hand side of this 
inequality is the absolute value of $(\U_H\L_K-\L_{H\star K})h$ evaluated
at $t$. Setting $s=l+j$, the right hand side is exactly the upper bound
given in the claim.\hfill$\qed$

\bigskip

\noindent{\bf Proof} (of Lemma \BasicLemmaU). The triangle inequality
implies the pointwise inequality
$$\displaylines{\qquad
    \Bigl| \U_H\calA_m\oG^{(k)}-\sum_{i\geq 1} \L_{H\star G\natural F_i,m}
    \oF_i^{(k)}\Bigr|
  \hfill\cr\noalign{\vskip 3pt}\hfill
    {}\leq \sum_{i\geq 1} |\U_H\L_{G\natural F_i,m}\oF_i^{(k)}
    -\L_{H\star G\natural F_i,m}\oF_i^{(k)} | \, .
  \qquad\equadef{BasicLemmaUEqA}\cr
  }
$$
The claim yields for any positive $t$,
$$\displaylines{\quad
  |\U_H\L_{G\natural F_i}\oF_i^{(k)} -\L_{H\star G\natural F_i}\oF_i^{(k)}|(t)
  \hfill\cr\noalign{\vskip 3pt}\qquad
    {}\leq \sum_{0\leq s\leq m} {|\oF_i^{(k+s)}(t)|\over s!} 
    \sum_{0\leq j\leq s}{s\choose j} \mu_{G\natural F_i,j} 
    \int_{-\infty}^{-t/2} |x|^{s-j} \d H(x)
  \hfill\cr\qquad\qquad
    {}+ {|\oF_i^{(k+m)}(t)|\over t^r m!} \sum_{0\leq j\leq m} {m\choose j}
    \mu_{G\natural F_i,j} 
  \hfill\cr\hfill
    {}\times \int_{-t/2}^0 
    \oDelta_{t,|x|/t}^r(\oF_i^{(k+m)}) |x|^{m-j+r} \d H(x)
  \cr\qquad\qquad
    {}+H(-t/2) \sum_{0\leq j\leq m} {\mu_{G\natural F_i,j}\over j!}
    \sup_{s\geq t} |\oF_i^{(k+j)} (s)| \, .
  \hfill\equadef{BasicLemmaUEqB}\cr
  }
$$
We use the same estimates as in the proof of Lemma \BasicLemmaC. So, Lemma
{\BoundDkMF} yields
$$
  |\oF_i^{(k+s)}(t)| \leq M c_i^{\alpha-\epsilon} t^{-k-s}\oF(t) \, ,
$$
and we also have the moment bound
$$
  \mu_{G\natural F_i,j}\leq \mu_{G,j}\leq M \, .
$$
In \BasicLemmaUEqB, the term
$$
  \int_{-t/2}^0 \oDelta_{t,|x|/t}^r(\oF_i^{(k+m)})  |x|^{m-j+r} \d H(x)
$$
is rewritten as a sum of an integral over $[\,-t/2,-v(t))$ plus an integral
over $[\,-v(t),0\,]$. In this decomposition, the second integral is at most
$$
  \oDelta_{t,v(t)/t}^r (\oF^{(k+m)})\mu_{\M_{-1}H,m-j+r} \, ,
$$
while the first one is at most
$$
  \oDelta_{t,1/2}^r(\oF^{(k+m)})\int_{-\infty}^{-v(t)} |x|^{m-j+r}\d H(x) \, .
$$
By assumption, for $t$ at least $2t_{14}$, the inequality
$$
  \int_{-\infty}^{-t/2} |x|^{s-j}\d H(x)
  \leq M_0 (t/2)^{s-j} F(-t/2)
$$
holds true. It follows that for $t$ more than some $t_{14}'$, the right hand 
side of {\BasicLemmaUEqB} is at most
$$\displaylines{\quad
    M c_i^{\alpha-\epsilon} t^{-k} \oF(t) F(-t/2)  
    + M c_i^{\alpha-\epsilon} t^{-k-m-r} \oF(t) 
    \Bigl( \oDelta_{t,v(t)/t}^r(\oF^{(k+m)})
  \hfill\cr\noalign{\vskip 2pt}\hfill
     {}+\oDelta_{t,1/2}^r(\oF^{(k+m)}) v(t)^{m+r} F\bigl(-v(t)\bigr) 
    \Bigr)
  \cr\noalign{\vskip 2pt}
    \phantom{\quad M c_i^{\alpha-\epsilon} t^{(-k)} \oF(t) F(-t/2) }
    + MH(-t/2) c_i^{\alpha-\epsilon} t^{-k} \oF(t) \, .
  \hfill\cr}
$$
Consequently, for $t$ at least $2t_{14}\wedge t_{14}'$, {\BasicLemmaUEqA} is 
at most
$$\displaylines{\quad
  MC_{\alpha-\epsilon} t^{-k-m}\oF(t) 
  \Bigl( t^m F(-t/2)+\oDelta_{t,v(t)/t}^r(\oF^{(k+m)}) 
  \hfill\cr\hfill
  {}+ \oDelta_{t,1/2}^r(\oF^{(k+m)}) v(t)^{m+r} F\bigl( -v(t)\bigr) 
  + t^m H(-t/2) \Bigr) \, .
  \cr}
$$
For $t$ large enough, $t^m F(-t/2)$ is at most 
$v(t)^{m+r}F\bigl( -v(t)\bigr)$. This concludes the proof of the 
lemma.\hfill$\qed$

\bigskip

\subsection{Conditional expansion and removing conditioning.}%
For any $I$, the distribution functions $\M_{-1}H_I$ and $G_I$ have
asymptotic expansions given by the results of section 5. In
particular, $H_I$ satisfies the assumptions of Lemma {\BasicLemmaU}
with a constant $M_0$ independent of both $I$ and the sequence
$c=(c_i)_{i\geq 1}$ which $N_{\alpha,\omega,\gamma}(c)$ at most $1$.

The result of section 5 shows that for some function $\eta(\cdot)$ which tends
to $0$ at infinity, for any $t$ at least some $t_{13}$, for any set $I$ and
any sequence $c$ with $N_{\alpha,\omega,\gamma}(c)$ at most $1$,
$$
  |\oG_I^{(k)}-\calA_m\oG_I^{(k)}|_{k+m,t} \leq \eta(t) \, , 
$$
where
$$
  \calA_m\oG_I^{(k)} 
  = \sum_{i\in \NN^*\setminus I}\L_{G_I\natural F_{+,i}}\oF_{+,i}^{(k)} \, .
$$

Next, Lemma {\DUCommute} shows that $\overline{H_I\star G_I}^{(k)}
=\U_{H_I}\oG_I^{(k)}$, which justifies equality \KEqual. Then, Lemma 
{\UContraction} and the result of section 5 yield
$$
  |K^{(k)}-E\U_{H_I}\calA_m\oG_I^{(k)}|_{k+m,t} \leq \eta(t) \, .
$$
Using Lemma \BasicLemmaU, we conclude that there exists a function 
$\eta^*(\cdot)$ with limit $0$ at infinity such that
$$
  \Bigl|K^{(k)}-E\sum_{i\in\NN^*\setminus I} \L_{H_I\star G_I\natural F_{+,i}}
  \oF_{+,i}^{(k)} \Bigr|_{k+m,t} \leq\eta^*(t)\, .
$$

To calculate the expectation involved in the inequality, we rewrite it as
$$
  \sum_{i\geq 1} E\bigl( \II\{\, X_i> 0\,\} \L_{H_I\star G_I\star F_{+,i},m}
  \bigr) \oF_{+,i} \, ,
$$
or, equivalently, since the $X_i$'s are independent and identically 
distributed, as
$$
  \sum_{i\geq 1} E(\L_{H_I\star G_I\natural F_{+,i},m} \, | \, X_i>0 )
  \oF(0) \oF_{+,i}^{(k)} \, .
$$
Conditionally upon $X_i$ being positive, $H_I\star G_I\natural
F_{i,+}$ is the distribution function of $\sum_{ \matrix{\ss j\geq
1\cr\noalign{\vskip -4pt}\ss j\not = i\cr}}c_iX_i$ given $I$. Thus,
the $X_i$'s being indpendent,
$$
  E( \L_{H_I\star G_I\natural F_{+,i}}\, |\, X_i>0 ) = \L_{K\natural F_i} \, .
$$
Next, $\oF(0) F_{+,i}^{(k)}=\oF_i^{(k)}$ on the positive half line, and so 
the expectation that we wanted to calculate is simply 
$\sum_{i\geq 1}\L_{K\natural F_i}\oF_i^{(k)}$.
This proves Theorem {\MainTheorem} when the $X_i$'s may assume arbitrary sign
and the constants $c_i$'s are nonnegative and $F$ is smooth.

\bigskip


\section{Removing the sign restriction on the constants.}

So far, we proved Theorem {\MainTheorem} assuming that the $c_i$'s are 
nonnegative. To drop this requirement is now rather easy, but unfortunately
requires some checking which is very much in the flavor of either section
5 or 6. So we will give the details only for part of the proof. In this 
section, we keep assuming that $\oF^{(k)}$ and $\overline{\M_{-1}F}{}^{(k)}$
exist, are bounded and continuous on the positive half line.

We define the set $J$ of all indices pertaining to a positive constant $c_i$,
that is
$$
  J=\{\, i\in \NN^* \, : \, c_i >0\,\} \, .
$$
Define $G$ (resp. $H$) to be the distribution function of 
\smash{$\sum_{i\in J}c_i X_i$} 
(resp. \smash{$\sum_{i\in \NN^*\setminus J} c_i X_i$}). The tail expansion of 
\smash{$\oG^{(k)}$} is given by the result of section 6. That 
of \smash{$\oH^{(k)}$} follows by the same token, since 
$$
  \sum_{i\in \NN^*\setminus J} c_i X_i
  =\sum_{i\in\NN^*\setminus J}(-c_i)(-X_i) \, .
$$
Now, the distribution function $K$ of $\<c,X\>$ is the convolution $G\star H$.

To obtain an asymptotic expansion of $\overline{G\star H}^{(k)}$, we
have two strategies. One is to decompose the distribution functions by
conditioning with respect to signs, as we did in section 6; another
one is to go more along the line of section 5. Both routes are about
equal length, and we go for the latter.

Our starting point is Proposition {\ConvolInOperators} which asserts that
$$
  \oK = \T_{G,1/2}\oH + \T_{H,1/2}\oG + \M_2 \oG \M_2 \oH
$$
and for a positive integer $k$,
$$
  \oK^{(k)} = \T_{G,1/2}\oH^{(k)} + \T_{H,1/2}\oG^{(k)} + 
  \sum_{1\leq i\leq k-1} \M_2 \oG^{(i)} \M_2 \oH^{(k-i)} \, .
  \eqno{\equa{KEq}}
$$

\bigskip

\subsection{Neglecting terms involving the multiplication operators.}%
Starting with \KEq, the purpose of this subsection is to show that
an expansion of $K^{(k)}$ can be obtained from the ones for 
$\T_{G,1/2}\oH^{(k)}$ and $\T_{H,1/2}\oG^{(k)}$, or, equivalently, that the
terms $\M_2\oG^{(i)}\M_2\oH^{(k-i)}$ can be neglected.

Our first lemma is an analogue of Lemma \UContraction. The operators $\M_c$
are not contraction, but for our problems, they behave very much as if
they were bounded, which is good enough.

\Lemma%
{\label{MBounded}%
  There exists $t_{16}$ such that for any positive real number $c$ and any 
  positive $t$ with  $(t/c)\wedge t$ at least $t_{16}$,
  $$
    |\M_c h|_{p,t} \leq (c^{p+\alpha+1}\vee c^{p+\alpha-1}) |h|_{p,t/2} \, .
  $$
}

\Proof We simply write
$$
  {s^p\over \oF(s)} h(s/c)
  = {(s/c)^p h(s/c)\over \oF(s/c)} \, c^p {\oF(s/c)\over \oF(s)}\, ,
$$
from which we deduce
$$
  |\M_ch|_{p,t} 
  \leq c^p |h|_{p,t/c} \sup_{s\geq t} \oF(s/c)/\oF(s) \, .
$$
We conclude in using Lemma \Potter.\hfill$\qed$

\bigskip

Our next lemma states a rather obvious property of the norm $|\cdot |_{p,t}$.

\Lemma
{\label{NormProduct}%
  For any positive real $p$, $q$ and $t$,
  $$
    |fg|_{p+q+m,t} \leq |f|_{p,t} |g|_{q,t}\, \sup_{s\geq t} s^m\oF(s) \, .
  $$
}

\Proof It follows from the identity
$$
  {s^{p+q+m}\over \oF(s) } f(s) g(s)
  =  {s^p\over \oF(s)} f(s) \, {s^q\over \oF(s)} g(s) \, s^m \oF(s) \, .
  \eqno{\qed}
$$

\bigskip

Our first two lemmas in this subsection imply the following bound.

\Lemma
{\label{NeglectMTerms}%
  Assume that the hypotheses of Theorem {\MainTheorem} hold. Suppose also 
  that $F^{(k)}$ exists, is continuous and bounded. Then, there exists a 
  function $\eta^*(\cdot )$ with limit $0$ at infinity and a positive number 
  $t_{17}$ such that for any $t$ at least $t_{17}$ and any sequence $c$ with 
  $N_{\alpha,\omega,\gamma}(c)$ at most $1$,
  $$
    \Bigl| \sum_{1\leq i\leq k-1} \M_2\oG^{(i)}\M_2\oH^{(k-i)}\Bigr|_{k+m,t}
    \leq \eta^*(t) \, .
  $$
}

\Proof The triangle inequality, Lemmas {\MBounded} and {\NormProduct}
yield for $t$ at least $1$,
$$\displaylines{\quad
    \Bigl| \sum_{1\leq i\leq k-1} \M_2\oG^{(i)}\M_2\oH^{(k-i)}\Bigr|_{k+m,t}
  \hfill\cr\hfill
    {}\leq \sum_{1\leq i\leq k-1} 2^{k+m+\alpha+1} |\oG^{(i)}|_{i,t/2}
    |\oH^{(k-i)}|_{k-i,t/2}\sup_{s\geq t/2}s^m\oF(s) \, .
  \quad\cr}
$$
From the result of section 6,
$$
  |\oG^{(i)}-\calA_0\oG^{(i)}|_{i,t/2}\leq \eta (t/2) \, .
$$
Since 
$$
  \calA_0\oG^{(i)} = \sum_{n\in J}\oF_n^{(i)} \, ,
$$
the triangle inequality and lemma {\BoundDkMF} imply that for $t$ at least
$t_1$,
$$
  |\calA_0\oG^{(i)}(t)| \leq M \sum_{n\in J} c_n^{\alpha-\epsilon} t^{-i}\oF(t)
  \, .
$$
Consequently,
$$
  |\oG^{(i)}|_{i,t/2}
  \leq  \eta(t/2) + M C_{\alpha-\epsilon} \, .
$$
We have a similar bound for $|H^{(k-i)}|_{k+m,t/2}$. Since $s^m\oF(s)=o(1)$
as $s$ tends to infinity, the result follows.
\hfill$\qed$

\bigskip

\subsection{Substituting \poorBold{$\oH^{(k)}$} and \poorBold{$\oG^{(k)}$}
by their expansions.}%
Looking at equality \KEq, we would like to use Theorem {\ApproxTByL} to 
conclude that $\T_{G,1/2}\oH^{(k)}$ has asymptotic expansion 
$\L_{G_m}\oH^{(k)}$. This is not possible since we do not know if $\oH^{(k)}$
is smoothly varying --- the problem is not so much about the regular variation
part of the smoothly varying condition, but more about the continuity of
$\oH^{(k+m)}$. The trick is to replace first $\oH^{(k)}$ by its asymptotic
expansion. Since the expansion involves explicitely $\oF$ and its derivatives,
we will be able to study the action of $\T_{G,1/2}$ on the expansion.
The first step is to show that $\T_{G,1/2}$ behaves like a bounded 
operator in our problem (though it is not.)

\Lemma
{\label{TBounded}%
  For $t$ at least $t_1$,
  $$
    |\T_{G,1/2}h|_{p,t} \leq 2^{p+\alpha+1} |h|_{p,t/2} \, .
  $$%
}

\Proof It follows the proof of Lemma \TNearContraction.\hfill$\qed$

\bigskip

This lemma and the result of section 6 imply
$$
  |\T_{G,1/2}\oH^{(k)}-\T_{G,1/2}\calA_m\oH^{(k)}|_{k+m,t}
  \leq 2^{k+m+\alpha+1} \eta (t/2) \, .
$$
Of course a similar inequality holds if we permute $G$ and $H$.

\bigskip

To conclude the proof of Theorem {\MainTheorem} when $F$ is smooth, we write
$$
  \T_{G,1/2}\calA_m\oH^{(k)} 
  = \sum_{n\in\NN^*\setminus J} \T_{G,1/2}\L_{H\natural F_n,m}\oF_n^{(k)} \, .
$$
Since $\L_{H\natural F_n,m}\oF_n^{(k)}$ involves derivatives of $\oF_n$, it
then suffices to study the expansion of $\T_{G,1/2}\oF_n^{(k)}$ and have a
good error bound in this expansion. This is now straightfoward. Indeed, write
$$
  \T_{G,1/2}\oF_n^{(k)}
  = \int_{-\infty}^{-t/2} \oF_n(t-x) \d G(x) 
  + \int_{-t/2}^{t/2} \oF_n^{(k)}(t-x) \d G(x) \, .
$$
We then bound
$$\eqalign{
  \Bigl| \int_{-\infty}^{-t/2} \oF_n^{(k)} (t-x) \d G(x)\Bigr|
  &{}\leq G(-t/2) \sup_{s\geq t/2} |\oF_n^{(k)}(s)| \cr
  &{}\leq M |c_n|^k G(-t/2) t^{-k}F(-t/2) \, , \cr
  }
$$
and we expand
$$
  \int_{-t/2}^{t/2} \oF_n^{(k)}(t-x) \d G(x)
$$
into $\L_{G,m}\oF^{(k)}$ as we did in sections 5 and 6.

This completes the proof of Theorem {\MainTheorem} when $F^{(k)}$ is
continuous and bounded.

\bigskip


\section{Removing the smoothness assumption.}

We now want to remove the assumption that $\oF^{(k)}$ 
exist, is continuous and bounded on the whole real line.

In Theorem \MainTheorem, the membership of $\oF$ and $\overline{\M_{-1}F}$
to $SR_{-\alpha,\omega}$ ensures that there exists some $A$ more than $2$ 
say such that $\oF^{(k)}$, $\overline{\M_{-1}F}{}^{(k)}$ exist, are continuous
and bounded on $(A,\infty)$.
Thus, we can write $F=(1-p)F_0+pF_1$ where $F_0$ is concentrated on $(-A-1,
A+1)$ and $F_1$ is concentrated on $(-A,A)^c$, and such that $F_1^{(k)}$
and $\overline{\M_{-1}F_1}{}^{(k)}$
exist, are continuous and bounded on the whole real
line. Let $\epsilon_i$ be a sequence of independent Bernoulli random variable
with mean $p$. Let $I$ be the random set of integer $i$ for which $\epsilon_1$
is $1$. Given $I$, let $G_I$ be $\star_{i\in\ZZ\setminus I}\M_{c_i}F$
and $H_I$ be $\star_{i\in I}\M_{c_i}F_0$. Then $\star_{i\in\ZZ}\M_{c_i}F_1$
is $EG_I\star H_I$. Since the $c_i$'s are summable, the support of the
distributions pertaining to $G_I$ is included in $|c|_1[\, -A-1,A+1\,]$
no matter what $I$ is. it follows that on $2|c|_1 (A+1,\infty)$,
$$
  \overline{G_I\star H_I}=\T_{G_I}\overline{H_I} \, .
$$
The expansion for $\overline{G_I\star H_I}$ follows from that of $H_I$ using
the same arguments that we used in sections 6 and 7. This concludes the proof
of Theorem \MainTheorem.

\bigskip

{\noindent\bf Appendix. {\bftt Maple} code.}
The goal in writing the following {\tt Maple} code was to see to which
extend the tail calculus described in section 3.2 could be
automatized.  With this code, the user enters a distribution function
having an expansion on some power of $t^{-1}$ and the integer $m$
occuring in Theorem \MainTheorem.  The program then generates the
coefficients of the asymptotic expansion of $G_c$ (with the notation
of Theorem \MainTheorem), assuming that the $c_i$'s are positive.

The example which we ran here is that of the Burr distribution in subsection
3.2. 

The tail $\oF$ is {\tt oF} (for 'overlined $F$'). One needs to specify 
the parameter $\alpha$, which is named {\tt palpha} and $m$.

\verbatim@
restart; with(LinearAlgebra):
pgamma:=10: ptau:=3/2:
oF:=(1+(t^(ptau))/beta)^(-pgamma); 
palpha:=ptau*pgamma; m:=4;

@

\noindent We then expand the tail of $F$. We write the tail expansion as
a polynomial of $x=1/t$.

\verbatim@
equiF:=subs(t=1/x,convert(asympt( oF,t,palpha+m),
                          polynom)):

@

\noindent We build the set of powers of $x$ involved in this expansion.

\verbatim@
index_set:={}:
for i from 1 by 1 to nops(equiF) do
  a:=op(i,equiF): 
  index_set:=index_set union {limit(log(a)/log(x),
                                    x=infinity)}:
end do:
power_list:=sort(convert(index_set,'list')):

@

\noindent We complete this list by the augmentation procedure 
described in section 3.2.

\verbatim@
for i in op(index_set) do
  for j from i by 1 while  j <= palpha+m do
    index_set:=index_set union {j};
  end do;
end do;
index_list:=sort(convert(index_set,'list')):

@

\noindent Then we calculate the dimension of the vector space we 
will work with.

\verbatim@
vdim:=nops(index_list):

@

We calculate the vector $p_{\oF}$. Unfortunately, with {\tt Maple}, it is
easy to work with polynomials, but it has no useful command to work with 
monomials of noninteger degree. That induces the following code where we
obtained the monomials one by one.

\verbatim@
pF:=Matrix(vdim,1,readonly=false):
equiF_tmp:=equiF:
for i from 1 by 1 to nops(power_list) do
  a:=op(i,equiF_tmp):
  xpower:=limit(log(a)/log(x),x=infinity):
  member(xpower,index_list,'k'):
  pF[k,1]:=a/x^xpower:
end do:

@

The next step is to build the matrices $\calD$ and $\calM_c$, which
we call {\tt Dmat} and {\tt Mcmat} in the code.

\hfuzz=3pt

\verbatim@
Dmat:=Matrix(vdim,readonly=false):
for i from 1 by 1 to vdim-1 do
  for j from i+1 by 1 to vdim do
    if evalb(index_list[j]-index_list[i]=1) 
      then Dmat[j,i]:=-index_list[i]; end if;
end do;end do;
Mcmat:=Matrix(vdim,readonly=false):
for i from 1 by 1 to vdim do 
  Mcmat[i,i]:=c^index_list[i]: 
end do:

@

\hfuzz=0pt

We construct the Laplace character of $F$ in its matrix form $\calL_F$.
The matrix is {\tt Lmat}.

\verbatim@
Lmat:=Matrix(vdim,readonly=false):
temp_mat:=Matrix(vdim,readonly=false):
for i from 1 by 1 to vdim do temp_mat[i,i]:=1: end do:
mu[0]:=1:
for k from 0 by 1 to m do
  Lmat:=(Lmat+((-1)^k*temp_mat*(c^k)*mu[k]/k!)):
  temp_mat:=(temp_mat.Dmat):
end do:

@

Its inverse is

\verbatim@
Lmat_inv:=MatrixInverse(Lmat,method='subs'):

@

We write the expansion for $G$ as
$\calL_G\sum_{i\in\ZZ}\calL_F^{-1}p_{\overline{\M_{c_i}F}}$.
We calculate the generic summands $\calL_F^{-1}p_{\overline{\M_cF}}$.

\verbatim@
Msum:=Lmat_inv.Mcmat.pF:

@

To sum these summands amounts to substitute $c^p$ by $C_p$ in this
sum, we we do now.

\verbatim@
for i from 1 by 1 to vdim do
  Msum[i,1]:=expand(collect(Msum[i,1],c),c):
end do:
for i from 1 by 1 to vdim do
  a:=Msum[i,1]:
  for j from 1 by 1 to nops(a) do
    if is ( c in indets(op(j,a)) )
    then
      p:=limit(log(op(j,a))/log(c),c=infinity):
      Msum[i,1]:=subs(c^p=C[p],Msum[i,1]):
    end if:
  end do:
end do:

@

To calculate the Laplace character of $G_c$, we first obtain its moments.
The algorithm is that described in Barbe and McCormick (2005). These
moments are coded as {\tt muG[k]}.

\verbatim@
for k from 0 by 1 to m do
  Q[k]:=t^k*c^k*mu[k]/k!:
end do:
P1:=add(Q[k],k=0..m):
P2:=convert(series(ln(P1),t=0,m+1),polynom):
P3:=add(C[k]*coeff(P2,c,k),k=0..m):
P4:=convert(series(exp(P3),t=0,m+1),polynom):
for k from 0 by 1 to m do:
  muG[k]:=k!*coeff(P4,t,k):
end do:

@

We construct the Laplace character of $G$.

\verbatim@
LGmat:=Matrix(vdim,readonly=false):
temp_mat:=Matrix(vdim,readonly=false):
for i from 1 by 1 to vdim do temp_mat[i,i]:=1: end do:
for k from 0 by 1 to m do
  LGmat:=(LGmat+((-1)^k*temp_mat*muG[k]/k!)):
  temp_mat:=(temp_mat.Dmat):
end do:

@

And finally obtain the coefficient of the tail, that is $p_{\oG}$.

\verbatim@
tail:=LGmat.Msum:

@

In our example for the Burr distribution, these coefficients are messy.
The remaining code makes them nicer looking; at least looking good
enough so that they can be written as we did.

The first step in the simplification is to substitute the centered moments
for the noncentered ones. In the code, we write {\tt kappa[k]} the $k$-the
centered moment. We express it as a function of the noncentered ones.

\verbatim@
for k from 1 by 1 to 5 do
  kappa[k]:=sum('(-1)^j*mu[j]*(mu[1]^(k-j))
                 *binomial(k,j)','j'=0..k):
  kappa[k]:=eval(kappa[k]):
end do:

@

Then we solve inductively for the noncentered moments, thereby expressing
them as functions of the centered ones.

\verbatim@
for k from 2 by 1 to 5 do 
  mu[k]:=[solve(kappa[k]=s[k],mu[k])][1];
end do:

@

Finally we do the substitution and arrange the terms by powers of $\beta$.

\verbatim@
for i from 1 to vdim do
  tail[i,1]:=collect(collect(simplify(
                        expand(tail[i,1])),beta),mu[1]);
end do;

@


\noindent{\bf Acknowledgements.} Philippe Barbe thanks
Peter Haskell and Konstantin Mischaikow for changing his views
on algebraic constructions through their marvelous teaching.

\bigskip


\noindent{\bf References}
\medskip

{\leftskip=\parindent
 \parindent=-\parindent
 \par

J.~Abate, G.L.~Choudhury, D.M.~Lucantoni, W.~Whitt (1995). Asymptotic
analysis of tail probabilities based on computation of moments, {\sl Ann.\
Appl.\ Probab.}, 5, 983--1007.

J.~Abate, G.L.~Choudhury, W.~Whitt (1994). Waiting-time tail probabilities
in queues with long-tail service-time distribution, {\sl Queueing Systems
Theory Appl.}, 16, 311--338.

J.~Abate, W.~Whitt (1997). Asymptotics for M/G/1 low-priority waiting-time
tail probability, {\sl Queueing Systems Theory Appl.}, 25, 173--233.

K.B.~Athreya, P.~Ney (1972). {\sl Branching Processes}, Springer.

Ph.~Barbe, W.P.~McCormick (2005). Asymptotic expansions of con\-volutions
of regularly varying distributions, {\sl J.\ Austr.\ Math.\ Soc.}, to appear.

Ph.~Barbe, W.P.~McCormick (200?). Tail calculus with remainder, applications
to tail expansions for infinite order moving averages, randomly stopped sums,
and related topics, submitted.

J.~Beirlant, J.L.~Teugels, P.~Vynckier (1996). {\sl Practical Analysis of 
Extreme Values}, Leuven University Press, Leuven, Belgium.

N.H.~Bingham, C.M.~Goldie, J.L.~Teugels (1989) {\sl Regular Variation},
2nd ed., Cambridge

R.~Bojanic, J.~Karamata (1963). On slowly varying functions and asymptotic
relations, Math.\ Research Center Tech.\ Report, 432, Madison, Wisconsin.

A.A.~Borovkov, K.A.~Borovkov (2003). On large deviation probabilities
for random walks, I, regularly varying distribution tails, {\sl Theory Probab.
Appl.}, 46, 193--213.

L.~Breiman (1968). {\sl Probability}, Addison-Wesley.

P.J.~Brockwell, R.A.~Davis (1991). {\sl Time Series: Theory and
Methods}, 2nd ed., Springer.

M.~Broniatowski, A.~Fuchs (1995). Tauberian theorems, Chernoff inequality,
and the tail behavior of finite convolutions of distribution functions,
{\sl Adv.\ Math.}, 116, 12--33.

J.~Chover, P.~Ney, S.~Wainger (1973). Functions of probability measures, 
{\sl J.\ Analyse Math.}, 26, 255--302.

Y.S.~Chow, H.~Teicher (1978). {\sl Probability Theory, Independence, 
Interchangeability, Martingales}, Springer.

J.W.~Cohen (1972). On the tail of the stationary waiting-time distribution
and limit theorem for M/G/1 queue, {\sl Ann.\ Inst.\ H.\ Poincar\'e},{\sl B},
8, 255--263. 

R.A.~Davis, M.~Rosenblatt (1991). Parameter estimation for some time series
models without contiguity, {\sl Statist. Prob. Lett.}, 11, 515--521.

J.~Delsarte (1938). Sur une extension de la formule de Taylor, {\sl Journ.
Math. Pures et Appl.}, 28(3), 213--231.

P.~Diaconis, D.~Freedman (1999). Iterated random functions, {\sl SIAM Review},
41, 45--70.

R.M.~Dudley (1989). {\sl Real Analysis and Probability}, Chapman \& Hall.

P.~Embrechts, C.M.~Goldie, N.~Veraverbeke (1979). Subexponentiality
and infinite divisibility, {\sl Z.\ Wahrsch.\ verw.\ Geb.}, 49, 335--347. 

P.~Embrechts, C.~Kl\"uppelberg, T.~Mikosch (1997). {\sl Modelling 
Extremal Events}, Springer.

W.~Feller (1971). {\sl An Introduction to Probability Theory and Its
Applications}, vol.~2, Wiley.

J.L.~Geluk (1992). Second order tail behaviour of a subordinated
probability distribution, {\sl Stoch.\ Proc.\ Appl.}, 40, 325--337.

J.L.~Geluk (1994). Asymptotic behaviour of the convolution tail of 
distributions each having a first or second order regularly varying tail,
{\sl Analysis}, 14, 163--183.

J.L.~Geluk (1996). Tails of subordinated laws: the regular varying case,
{\sl Stoch.\ Proc.\ Appl.}, 61, 147--161.

J.L.~Geluk, L.~de Haan, S.~Resnick, C.~Starica (1997). Second-order regular 
variation, convolution and the central limit theorem, {\sl Stoch.\ Proc.\  
Appl.}, 69, 139--159.

J.~Geluk, L.~Peng, C.G.~De Vries (2000). Convolutions of heavy-tailed 
random variables and applications to portfolio diversification and 
MA($1$) time series, {\sl Adv.\ Appl.\ Prob.}, 32, 1011--1026.

C.M.~Goldie (1991). Implicit renewal theory and tails of solutions of random
equations, {\sl Ann.\ Appl.\ Probab.}, 1, 126--166.

D.R.~Grey (1994). Regular variation in the tail behavior of solutions
of random difference equations, {\sl Ann.\ Probab.}, 4, 169--183.

A.K.~Grincevi\v cius (1975). On limit distribution for a random walk on the 
line, {\sl Lithuanian Mat.\ J.}, 15, 580--589.

R.~Gr\"ubel (1987). On subordinated distributions and generalized renewal
measures, {\sl Ann.\ Probab.}, 15, 394--415.

P.~Hall, I.~Weissman (1997). On the estimation of extreme tail probabilities,
{\sl Ann. Statist.}, 25, 1311--1326.

B.M.~Hill (1973). A simple general approach to inference about the tail
of a distribution, {\sl Ann.\ Statist.}, 3, 1163--1174.

H.~Kesten (1973). Random difference equations and renewal theory for product
of random matrices, {\sl Acta. Math.}, 131, 207--248.

L.~Mattner (2004). Cumulants are universal homomorphisms into Hausdorff
groups, {\sl Probab.\ Theor.\ Rel.\ Fields}, 130, 151--166.

V.~Mari\'c (2000). {\sl Regular Variation and Differential Equations},
{\sl Lecture Notes in Mathematics}, 1726, Springer.

V.~Mari\'c, M.~Tomi\'c (1977). Regular variation and asymptotic properties
of solutions of nonlinear differential equations, {\sl Publ.\ Inst.\
Math.\ (Beograd)}, 21 (35), 119--129.

A.~Nijenhuis, H.~Wilf (1978). {\sl Combinatorial Algorithms},
Academic Press, second edition.

F.W.J.~Olver (1974). {\sl Asymptotic and Special Functions}, Academic Press.

E.~Omey (1981). Regular variation and its applications to second order
linear differential equations, {\sl Bull.\ Soc.\ Math.\ Belg.}, 32, 207--229.

E.~Omey (1988). Asymptotic properties of convolution products of functions,
{\sl Publ.\ Inst.\ Math.\ (Beograd) (N.S.)}, 43, 41--57.

E.~Omey, E.~Willekens (1986). Second order behavior of the tail of a 
subordinated probability distribution, {\sl Stoch.\ Proc.\ Appl.}, 21,
339--353.

E.~Omey, E.~Willekens (1987). Second order behavior of distributions
subordinated to a distribution with finite mean, {\sl Comm.\ Statist.\ Stoch.\
Models}, 3, 311--342.

A.G.~Pakes (1975). On the tails of waiting-time distributions, {\sl J.\ Appl.\
Probab.}, 12, 555-564.

A.G.~Pakes (2004). Convolution equivalence and infinite divisibility,
{\sl J.\ Appl.\ Probab.}, 41, 407--424.

S.~Resnick (1986). Point processes, regular variation and weak convergence,
{\sl Adv.\ Appl.\  Probab.}, 18, 66--138.

S.~Resnick (1987). {\sl Extreme Values, Regular Variation, and Point
Processes}, Springer.

S.~Resnick, C.~St\u aric\u a (1997). Asymptotic behavior of Hill's
estimator for autoregressive data, heavy tails and highly volatile
phenomena, {\sl Comm.\ Statist., Stoch.\ Models}, 13, 703--721.

S.~Resnick, E.~Willekens (1991). Moving averages with random coefficients
and random coefficients autoregressive models, {\sl Comm.\ Statist.\ Stoch.\
Models}, 7, 511--525.

T.~Rolski, H.~Schmidli, V.~Schmidt, J.~Teugels (1999). {\sl Stochastic
Processes for Insurance and Finance}, Wiley.

B.~Solomyak (1995). On the random series $\sum\pm\lambda^n$ (an 
Erd\"os problem), {\sl Ann. Math.}, 142, 611--625.

D.~Stanton, D.~White (1986). {\sl Constructive Combinatorics},
Springer.

E.~Willekens (1989). Asymptotic approximation of compound distributions
and some applications, {\sl Bull.\ Soc.\ Math.\ Belg., ser.\ B}, 41, 55-61.

E.~Willekens, J.L. Teugels (1992). Asymptotic expansions for waiting time
probabilities in an M/G/1 queue with long-tailed service time, 
{\sl Queueing Systems Theory Appl.}, 10, 295--311.

G.E.~Willmot, X.S.~Lin (2000). {\sl Lundberg Approximations for Compound
Distributions with Insurance Applications}, {\sl Lecture Notes in Statistics},
156, Springer.

}

\vskip .22in
\halign{#\hfil&\hskip 40pt #\hfill\cr
  Ph.\ Barbe            & W.P.\ McCormick\cr
  90 rue de Vaugirard   & Dept.\ of Statistics \cr
  75006 PARIS           & University of Georgia \cr
  FRANCE                & Athens, GA 30602 \cr
                        & USA \cr
                        & bill@stat.uga.edu \cr}

\bye